\documentclass[11pt]{amsart}
\usepackage[dvipsnames,usenames]{color}
\usepackage{hyperref}
\usepackage{float}
\usepackage{graphicx, import}
\graphicspath{ {./images/} }
\usepackage[latin1]{inputenc}
\usepackage{amsmath}
\usepackage{amsfonts}
\usepackage{amssymb}
\usepackage{amsthm}
\usepackage{amscd}
\usepackage{verbatim}
\usepackage{subfigure}
\usepackage{caption}
\usepackage{pinlabel}
\usepackage{stmaryrd}
\usepackage{enumerate, enumitem}
\usepackage{todonotes}
\usepackage{bm}
\usepackage{thmtools}
\usepackage{thm-restate}
\usepackage{lipsum}
\usepackage{setspace}
\usepackage{mathtools}
\usepackage[abs]{overpic}
\usepackage{mathdots}
\usepackage{tikz}
\usetikzlibrary{arrows}
\usetikzlibrary{decorations.pathreplacing}
\usepackage{verbatim}
\usetikzlibrary{cd}
\tikzset{taar/.style={double, double equal sign distance, -implies}}
\tikzset{amar/.style={->, dotted}}
\tikzset{dmar/.style={->, dashed}}
\tikzset{aar/.style={->, very thick}}
\usepackage[section]{placeins}

\newtheorem{theorem}{Theorem}[section]
\newtheorem{Theorem}{Theorem}
\newtheorem{lemma}[theorem]{Lemma}
\newtheorem{proposition}[theorem]{Proposition}

\newtheorem{corollary}[theorem]{Corollary}

\newtheorem{Conjecture}[Theorem]{Conjecture}

\theoremstyle{definition}
\newtheorem{definition}[theorem]{Definition}

\theoremstyle{remark}
\newtheorem{remark}[theorem]{Remark}

\def\F{\mathbb{F}}
\def\N{\mathbb{N}}
\def\C{\mathbb{C}}
\def\Q{\mathbb{Q}}
\def\R{\mathbb{R}}
\def\Z{\mathbb{Z}}

\def\talpha{\mathbb{T}_{\bm{\alpha}}}
\def\tbeta{\mathbb{T}_{\bm{\beta}}}

\def\Sym{\mathrm{Sym}}

\def\co{\colon }

\def\HF {\mathit{HF}}
\def\HFK {\mathit{HFK}}
\def\HFL{\mathit{HFL}}
\newcommand\HFKhat{\widehat{\HFK}}
\newcommand\HFLhat{\widehat{\HFL}}

\newcommand\HFhat{\widehat{\HF}}

\newcommand{\M}{M^{\operatorname{inv}}}
\newcommand{\Lo}{L_0^{\operatorname{inv}}}
\newcommand{\Li}{L_1^{\operatorname{inv}}}

\author[A. Parikh]{Aakash Parikh}
\email{ap1792@math.rutgers.edu}
\address{Rutgers University, New Brunswick, NJ, USA}
\thanks{The author was partially supported by NSF CAREER Grant DMS-2019396.}

\numberwithin{equation}{section}

\title{Localization and the Floer homology of strongly invertible knots}

\begin{document}

\begin{abstract}  
We establish two spectral sequences in knot Floer homology associated to a directed strongly invertible knot K. The first is from the knot Floer homology of K to the Heegaard Floer homology of the three sphere. The second is from the singular knot Floer homology of a singular knot associated to K to the knot Floer homology of the quotient knot of K. 
\end{abstract}
\maketitle
\tableofcontents
	\section{Introduction}
\subsection{Main results}\label{mainresults}
A \emph{symmetric knot} $(K,\tau)$ is a knot $K\subset S^3$ along with a finite order diffeomorphism $\tau\colon (S^3,K)\to (S^3,K)$. A symmetric knot $K$ is \emph{strongly invertible} if $\tau$ is order 2, orientation preserving, and $|\textrm{Fix}(\tau)\cap K|=2$. The study of strongly invertible knots was initiated by Sakuma in \cite{Sak}. 

The purpose of this paper is to establish two new spectral sequences in knot Floer homology for strongly invertible knots. Throughout $\F\coloneq \F_2$ denotes the field of two elements.

\begin{theorem}\label{Mainthm1}
    Given $K$ a strongly invertible knot, there is a spectral sequence with $E_1$ page equal to $\HFKhat(K,0)\otimes \F[\theta,\theta^{-1}]$ and $E_\infty$ page isomorphic to $\HFhat(S^3)\otimes \F[\theta,\theta^{-1}]$. Every page of this spectral sequence is an invariant of $K.$
\end{theorem}

From Theorem \ref{Mainthm1} we extract a numerical invariant $s_\tau$ of strongly invertible knots constructed similarly to many concordance invariants including Rasmussen's $s$ invariant \cite{rasmussen2004khovanov}, and Hendricks-Lipshitz-Sarkar's $q_\tau(K)$ and $d_\tau(K,m)$ \cite{Hendricks_2016}, and conjecture that this invariant is an equivariant concordance invariant.

\begin{figure}[htb!]
 \begin{center}
 \includegraphics[scale=.3]{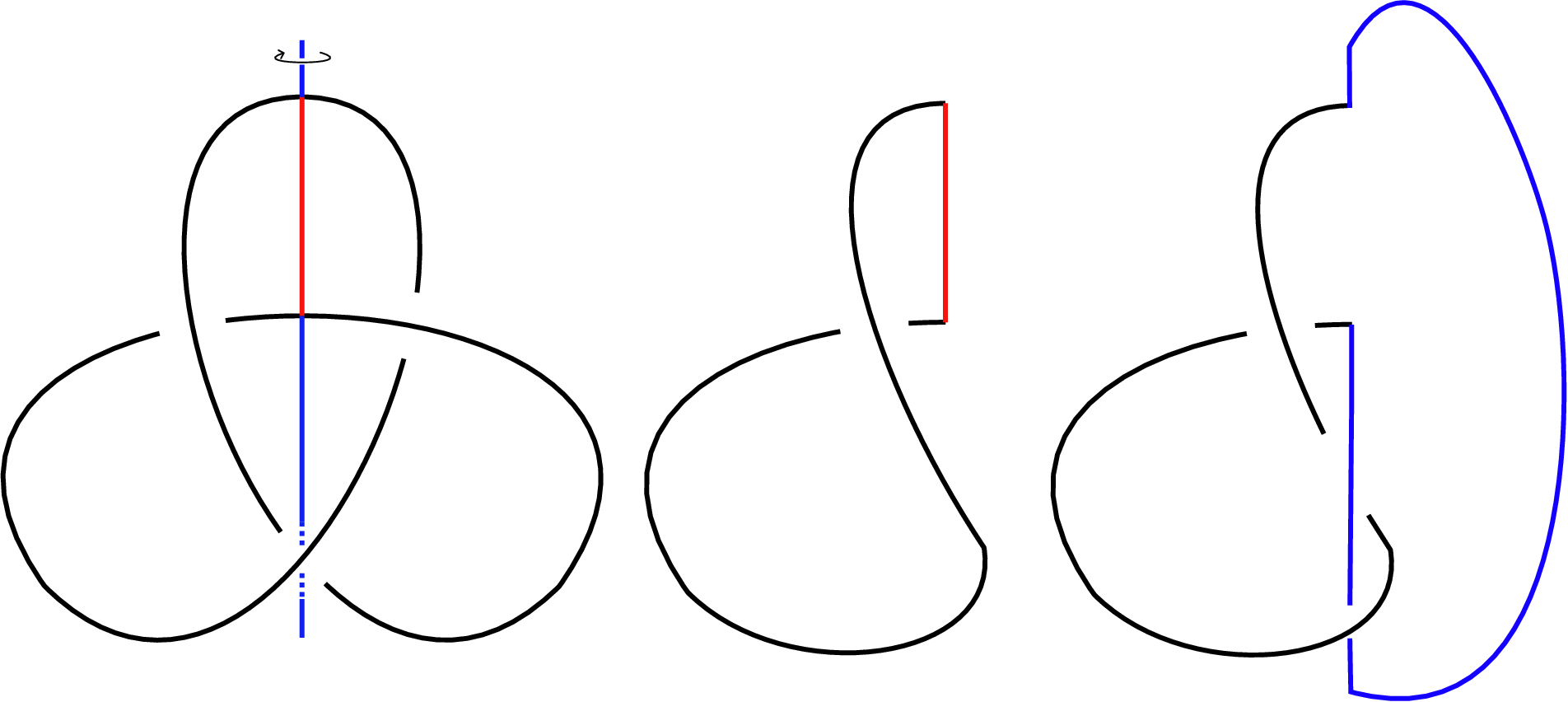}\end{center}
 \caption{Left: The left handed trefoil $3_1$ with strong inversion given by a $\pi$-radian rotation about the vertical symmetry axis. Middle: The quotient unknot of the directed strong inversion on $3_1$ denoted $3_1^+$ with the short vertical red segment chosen as the distinguished half axis. Right: The quotient trefoil of the directed strong inversion on $3_1$ denoted $3_1^-$ with the opposite choice of distinguished half axis.} \label{trefoilwquotients}
 \end{figure}

Let $K$ be a \emph{directed} strongly invertible knot (DSI) with quotient knot $\overline{K}$ as in Definition \ref{defstroquot} and Figure \ref{trefoilwquotients}.

\begin{figure}[htb!]\begin{center}
 \includegraphics[scale=.3]{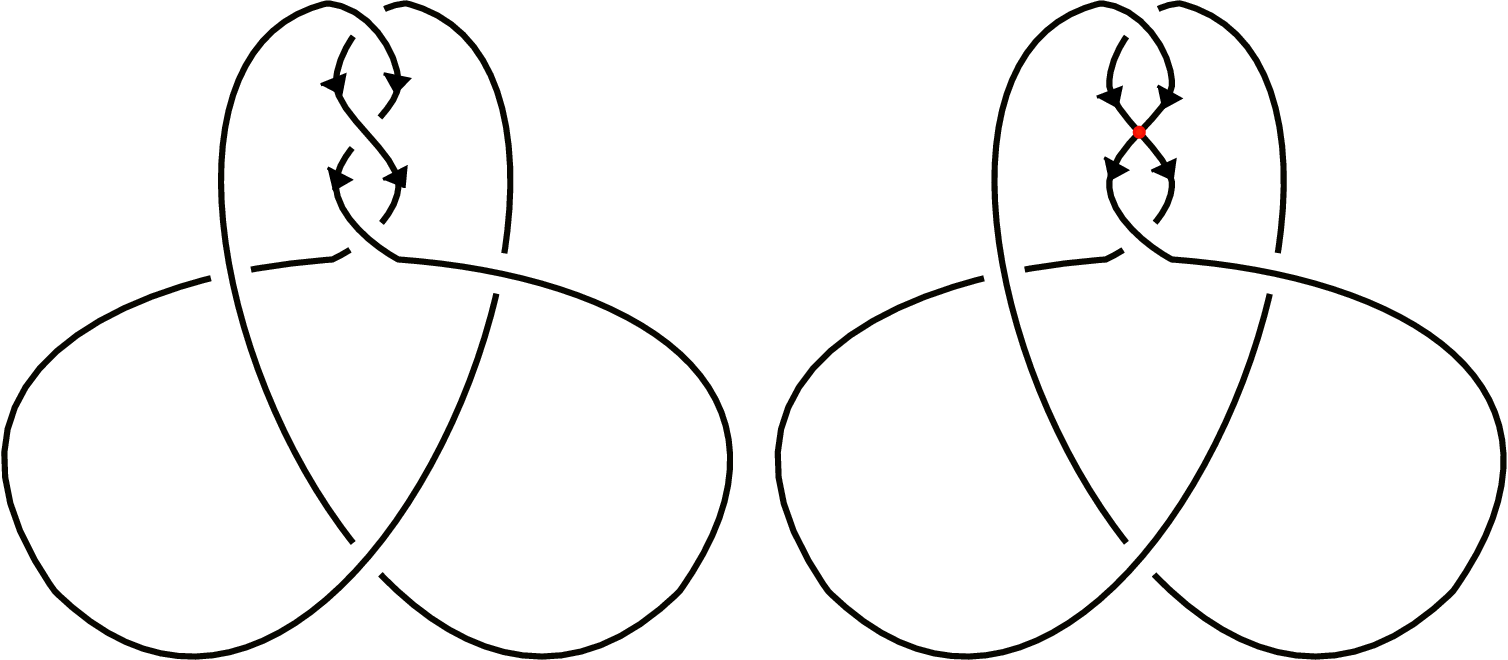}
  \end{center}
 \caption{Left: The butterfly link $L_b(3_1^+)$. Right: The singular butterfly link $S_b(3_1^+)$.} \label{LbSb}
 \end{figure}

The next spectral sequence makes use of
\begin{itemize}
    \item \emph{singular} knot Floer homology, a Floer homology group associated to a singular knot due to Ozsv\'{a}th, Szab\'{o} and Stipsicz \cite{OSsing1}, \cite{OSsing2}, and
    \item a singular knot $S_b(K)$ constructed from Boyle and Issa's \emph{butterfly link} $L_b(K)$, and also their \emph{axis linking number} $\widetilde{lk}(K)$ \cite{BoyAhE4G}.
\end{itemize}Section \ref{stro} contains the definitions of $L_b(K),$ $S_b(K)$ and $\widetilde{lk}(K)$ (see also Figure \ref{LbSb}), and Section \ref{singularlinkFloer} is a brief review of singular knot Floer homology.

\begin{theorem}\label{Mainthm2}
    Let $K$ be a DSI with quotient knot $\overline{K}$. There is a spectral sequence with $E_1$-page equal to
    \[\HFKhat(S_b(K))\otimes \F[\theta,\theta^{-1}]\] and $E_{\infty}$-page isomorphic to  \[\HFKhat(\overline{K})\otimes \F[\theta,\theta^{-1}].\] Furthermore, this spectral sequence splits along Alexander gradings:
    \begin{itemize}
        \item For each $a\in \Z$ there is a spectral sequence with $E_1$ page equal to \[\HFKhat(S_b(K), 2a-\widetilde{lk}(K)-\frac{1}{2})\otimes \F[\theta,\theta^{-1}]\] and $E_{\infty}$ page isomorphic to \[\HFKhat(\overline{K}, a)\otimes \F[\theta,\theta^{-1}].\]
        \item For any Alexander grading $A$ which cannot be written as $A=2a-\widetilde{lk}(K)-\frac{1}{2}$, the spectral sequence has $E_1$ page equal to $\HFKhat(S_b(K), A)\otimes \F[\theta,\theta^{-1}]$ and $E_\infty$ page equal to $0$.
    \end{itemize} 
    Every page of this spectral sequence is an invariant of $K.$
\end{theorem}

\subsection{Background}\label{background}
 There has been a recent burst of activity in the study of strongly invertible knots and their relation to 3- and 4-dimensional topology. For instance, Boyle--Issa studied equivariant versions of 3- and 4-genera \cite{BoyAhE4G}, and Hirasawa--Hiura--Sakuma computed the equivariant 3-genus for all 2-bridge strongly invertible knots \cite{HirasawaHiuraSakuma2023}; Alfieri--Boyle introduced an equivariant knot signature \cite{alfieri2021strongly} and used it to give a lower bound on the \emph{butterfly} 4-genus;  Di Prisa showed that the equivariant concordance group is non-abelian \cite{diprisa2022equivariant}; Dai--Mallick--Stoffregen introduced equivariant concordance invariants derived from knot Floer homology and used them to give lower bounds on the equivariant 4-genus \cite{DMSE}; Dai--Kang--Mallick--Park--Stoffregen studied the $(2,1)$-cable of $4_1$ and proved that it is not slice by showing that the branched double cover $\Sigma((4_1)_{2,1})\cong S_{+1}(4_1\#4_1^r)$ does not bound an equivariant $\Z/2\Z$ homology ball, in part by studying the swapping strong inversion on $4_1\#4_1^r$; Lobb--Watson constructed a spectral sequence involving a refinement of Khovanov homology in the presence of an involution for a strongly invertible knot \cite{Lobb_2021}; Lipshitz-Sarkar constructed another spectral sequence involving Khovanov homology of a DSI and the annular Khovanov homologies of its two annular quotients; and, Hendricks--Mak--Raghunath constructed the analog of the Lipshitz-Sarkar spectral sequence using symplectic Khovanov homology and symplectic annular Khovanov homology \cite{hendricks2024symplecticannularkhovanovhomology}.

Theorems \ref{Mainthm1} and \ref{Mainthm2} are inspired by Lipshitz and Sarkar's spectral sequence for the Khovanov homology strongly invertible knots, established using the Khovanov homotopy type and the localization spectral sequence for Borel equivariant cohomology \cite{KhSI}. An annular link is a link in a thickened annulus $D^2\times S^1$, or equivalently a link in $S^3$ along with a choice of unknotted axis. Given an intravergent diagram (as in Definition \ref{trans&intra}) $D_n$ for a DSI $K$, there are two naturally associated annular quotient knots, $\overline{K}_0$ and $\overline{K}_1$ -- the quotients of the 2-periodic $0$ and $1$ resolutions of the central crossing of $D_n$. The subscript $n$ in $D_n$ is the winding number of the annular knot $\overline{K}_0$ or equivalently the linking number of $\overline{K}_0$ with the symmetry axis $A$. Lipshitz and Sarkar define an \emph{axis moving map} \[f^+\colon ACKh(\overline{K}_1)\to ACKh(\overline{K}_0)\]  on the annular Khovanov chain complexes of $\overline{K}_1$ and $\overline{K}_0$. Let $AKh(\overline{K}_0,\overline{K}_1)=H_*(Cone(f^+)).$ 
\begin{theorem}\cite[~Theorem 1.1]{KhSI}
There is a spectral sequence \[Kh(K, \F)\otimes \F[\theta, \theta^{-1}]\rightrightarrows AKh(\overline{K}_0, \overline{K}_1)\otimes \F[\theta,\theta^{-1}].\]

\end{theorem}
 
Our techniques are inspired by Hendricks' work concerning localization spectral sequences for link Floer homology of doubly periodic knots \cite{HendricksDP}.
A \emph{doubly periodic link} is a link $P\subset S^3$ equipped with an order 2 orientation preserving diffeomorphism $\tau\colon  (S^3,P)\to (S^3,P)$  that preserves the orientation of $P$ and has non-empty fixed set. The fixed point set of $\tau$ is an unknot $A$ disjoint from $P$ called the \emph{axis of symmetry}. The \emph{quotient link} $\overline{P}$ is the image of $P$ under the quotient map $\ell\colon S^3\to S^3/\tau\cong S^3$, and $\overline{A}\coloneq \ell(A)$ is the \emph{quotient axis}.

\begin{theorem}\cite{HendricksDP, Hendricks_2016}\label{HendricksSS}
There is a spectral sequence 
with $E_1$ page equal to $\HFKhat(P)\otimes V\otimes W\otimes \F[\theta,\theta^{-1}]$ and $E_\infty$ page isomorphic to $\HFKhat(\overline{P})\otimes W\otimes  \F[\theta,\theta^{-1}]$. Furthermore this spectral sequence splits along Alexander gradings:
\begin{itemize}
\item For any $a\in \Z$ there is a spectral sequence with $E_1$ page equal to \[\HFKhat(P, 2a+\frac{1-\lambda}{2})\otimes V\otimes W\otimes \F[\theta,\theta^{-1}]\] and $E_\infty$ page isomorphic to \[\HFKhat(\overline{P}, a+\frac{1-\lambda}{2})\otimes W\otimes \F[\theta,\theta^{-1}].\]
\item For any Alexander grading $A$ which cannot be written as $A=2a+\frac{1-\lambda}{2}$, there is a spectral sequence with $E_1$ page equal to $\HFKhat(P)\otimes V\otimes W\otimes \F[\theta,\theta^{-1}]$ and $E_\infty$ page equal to $0$.
\end{itemize}
where $V$ and $W$ are bigraded two dimensional vector spaces over $\F$ and $\lambda=lk(P, A)=lk(\overline{P}, \overline{A}).$ Every page of this spectral sequence is an invariant of $(P,\tau)$.\end{theorem} 
The spectral sequences of Theorems \ref{Mainthm1}, \ref{Mainthm2} and \ref{HendricksSS} are constructed via the following procedure:
\begin{enumerate}
    \item Construct a  $\tau$-equivariant Heegaard diagram $\mathcal{H}$ for the symmetric knot.
    \item Identify the fixed point sets of the symmetric product and tori associated to $\mathcal{H}$ with the symplectic manifold associated to the Heegaard diagram $\overline{\mathcal{H}}=\mathcal{H}/\tau$ for the quotient.
    \item In this setting verify the symplectic hypotheses for a \emph{localization isomorphism} in $\Z/2\Z$-equivariant Floer homology are met.
\end{enumerate}

\begin{figure}
    \center
    \includegraphics[scale=.3]{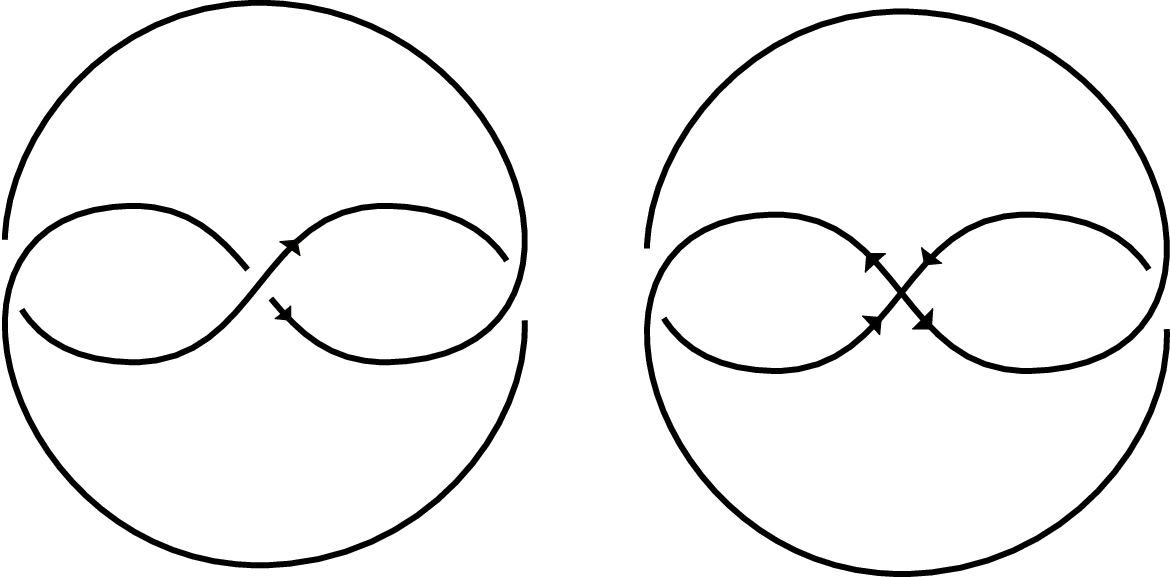}
    \caption{Left: An intravergent diagram for $3_1^+$. Right: An intravergent diagram for $S_b^1(3_1^+)$.}
    \label{fig:31s1b31}
\end{figure}

\begin{figure}
    \centering
    \includegraphics[width=0.5\linewidth]{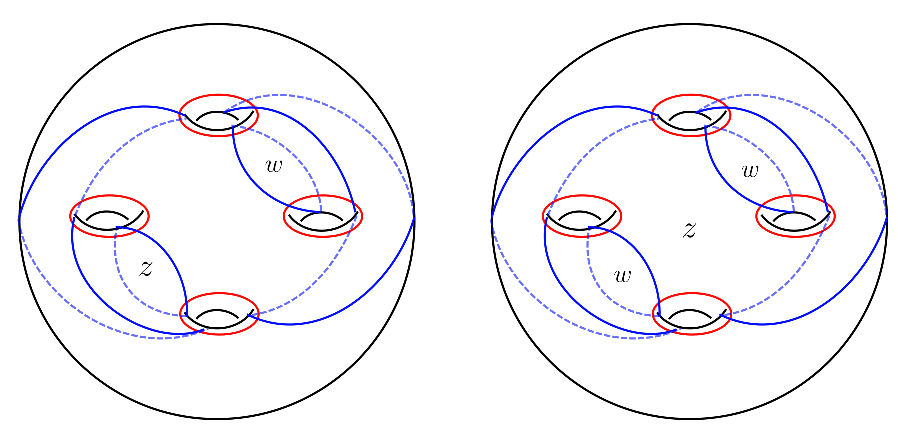}
    \caption{Left: A symmetric Heegaard diagram for $3_1^+$. Right: A symmetric Heegaard diagram for $S_b^1(3_1^+)$.}
    \label{TrefoilHeegaardDiagram}
\end{figure}

In Hendricks' work the result used for (3) is Seidel--Smith's localization isomorphism \cite{SS}, and in this paper we make use of Large's more general localization isomorphism \cite{Large}.
The equivariant Heegaard diagrams constructed in \cite{HendricksDP} have underlying surface $S^2$ and multiple basepoints because this is the setting in which one of the restrictive hypotheses required by Seidel-Smith's localization theorem -- the existence of a \emph{stable normal trivialization} --  can be met. The flexibility afforded by the analogous weaker hypothesis in Large's localization theorem -- a \emph{stable tangent normal isomorphism} -- allows us to make use of higher genus equivariant Heegaard diagrams with minimal basepoints. This is of particular relevance in Theorem \ref{Mainthm1}, because an equivariant Heegaard diagram for a DSI with more than 2 basepoints induces a spectral sequence that abuts to 0. For the benefit of the expert reader, we include an example of the symmetric Heegaard diagrams associated to an intravergent diagram for $3_1^+$ and $S_b^1(3_1^+)$ (cf. Definition \ref{Snb}) in Figures \ref{fig:31s1b31} and \ref{TrefoilHeegaardDiagram}. 

\subsection*{Organization} This paper is organized as follows. In Section \ref{stro} we recall background on strongly invertible knots and singular links, Heegaard Floer homology, link Floer homology and singular link Floer homology. Equivariant Heegaard diagrams for strongly invertible knots are constructed in Section \ref{diagrams}. In Section \ref{equivariant} we introduce Large's localization theorem and use it to prove Theorems \ref{Mainthm1} and \ref{Mainthm2}, deferring the verification of Large's hypotheses to later sections. Section \ref{examples} contains examples of Theorems \ref{Mainthm1} and \ref{Mainthm2}. We study the homotopy type and cohomology of the symmetric product of punctured surfaces in Section \ref{geosym}. Section \ref{stableTNiso} contains a proof of the existence of a stable tangent normal isomorphism. Appendix \ref{Appendix} is an exposition of a grid diagram based proof of the skein exact triangle for singular knot Floer homology, and an explanation of how this proof pins down absolute Alexander gradings for singular knot Floer homology. We use this to determine the Alexander grading shift in Theorem \ref{Mainthm2}.
\subsection*{Acknowledgments}We would like to thank Kristen Hendricks for suggesting this project
and providing guidance and support throughout. We are also grateful to Akram Alishahi,
Anna Antal, Keegan Boyle, Jen Hom, Tye Lidman, Robert Lipshitz, Abhishek Mallick, Danielle O'Donnol and Sriram Raghunath for helpful conversations. In addition, we thank the anonymous referee whose comments greatly enhanced the clarity of exposition in many places. Portions of this work were carried out at the conference ``Gauge Theory and Topology: in Celebration of Peter Kronheimer's 60th Birthday'' and the 2024 Geometry Topology Workshop at UCLA -- we thank the organizers for their hospitality.

\section{Background }\label{stro}

\subsection{Topological preliminaries}
In this section we recall basics about strongly invertible knots and singular links. 
\subsubsection{Strongly invertible knots}
\begin{definition}\label{defstroquot} A knot $K\subset S^3$ is \emph{strongly invertible} if there is an order 2 orientation preserving diffeomorphism $\tau\colon (S^3,K)\to (S^3, K)$ with non-empty fixed point set that reverses the orientation on $K$. The fixed point set of $\tau$ is called the \emph{axis (of symmetry)}, denoted by \[A\coloneq \{x\in S^3|\tau(x)=x\}.\] The axis $A$ intersects $K$ in two points. By the resolution of the Smith conjecture, the axis $A$ is an unknot. The two points of $\textrm{Fix}(\tau)\cap K$ separate $A$ into two pieces $A=A_1\cup A_2$ called half-axes. The pair $(K,\tau)$ along with the choice of an oriented half axis, without loss of generality say $A_1$, is called a \emph{directed strongly invertible knot} (DSI). Say that $\ell\colon S^3\to S^3/ \tau\cong S^3$ is the quotient map induced by $\tau$. Defining $\overline{A}\coloneq \ell(A)$, and $\overline{A}_1=\ell(A_1)$, we see that $\ell(K)$ is an arc with endpoints on the unknotted axis $\overline{A}\subset S^3$. The \emph{quotient knot} $\overline{K}$ of a DSI $(K,\tau,A_1)$ is the union of the arc and chosen half axis: $\overline{K}\coloneq \ell(K)\cup \overline{A}_1$.
\end{definition}
It is useful to have diagrams for strongly invertible knots that display their symmetry.
\begin{definition}\cite[~Definition 3.3]{BoyAhE4G}\label{trans&intra}
Let $(K,\tau)$ be a strongly invertible knot. A knot diagram for $K$ is 
\begin{enumerate}
    \item \emph{transvergent} if $\tau$ acts as rotation around an axis contained in the diagram, and
    \item \emph{intravergent} if $\tau$ acts as rotation around an axis perpendicular to the plane of the diagram. 
\end{enumerate} 
\end{definition}

\begin{figure}[htb!]
 \center
 \includegraphics[scale=.4]{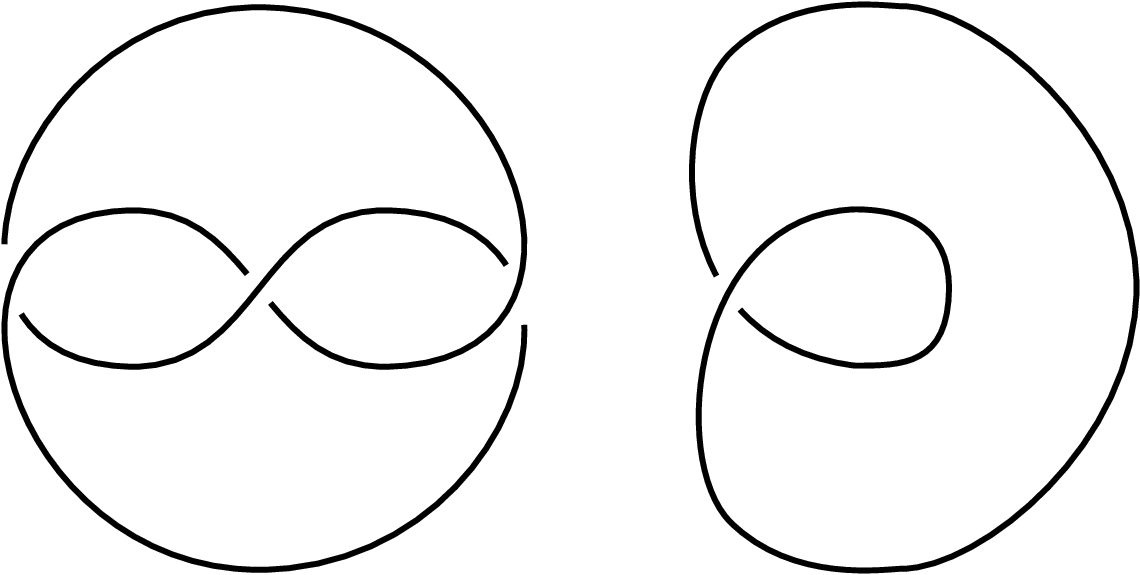}
 \caption{Left: An intravergent diagram for $3_1^+$. Right: The quotient unknot of $3_1^+$.} \label{intratrefquot}
 \end{figure}

Every strongly invertible knot admits both transvergent and intravergent diagrams.
An intravergent diagram for a strongly invertible knot distinguishes a half-axis, namely the one that lies between the over and under strands of the central crossing, and we choose the convention that this half axis is oriented from the over strand to the under strand. See the left hand side of Figure \ref{trefoilwquotients} for a transvergent diagram displaying the unique strong inversion on the trefoil, and the left hand side of Figure \ref{intratrefquot} for an intravergent diagram for the directed strong inversion on the trefoil which we denote by $3_1^+.$ Notice that this DSI has an unknotted quotient, $\overline{3_1^+}=0_1$.

\begin{definition} \label{equivalent}
Strongly invertible knots $(K,\tau)$ and $(C,\rho)$ are \emph{equivariantly isotopic} or \emph{equivalent} if there is an orientation preserving diffeomorphism $\phi\colon S^3\to S^3$ such that $\phi(K)=C$ and $\phi\circ \tau=\rho\circ \phi$. If $K$ and $C$ are directed with oriented half axes $A_1$ and $B_1$ respectively they are \emph{equivalent as DSIs} if they are equivariantly isotopic, $\phi(A_1)=B_1$, and $\phi$ preserves the chosen orientations on the half axes.
\end{definition}
The data of a DSI consists of the knot $K$ along with the involution $\tau$ and the choice of oriented half axis $A_0$. Indeed, the quotient knot of a directed strong inversion depends on the choice of half axis; for example the strong inversion on the trefoil yields quotient knots equal to the unknot or the trefoil depending on the choice of half axis as seen in Figure \ref{trefoilwquotients}.  However, we often abuse notation and refer to $K$ as the DSI.

\begin{figure}[htb!]
 \includegraphics[scale=.25]{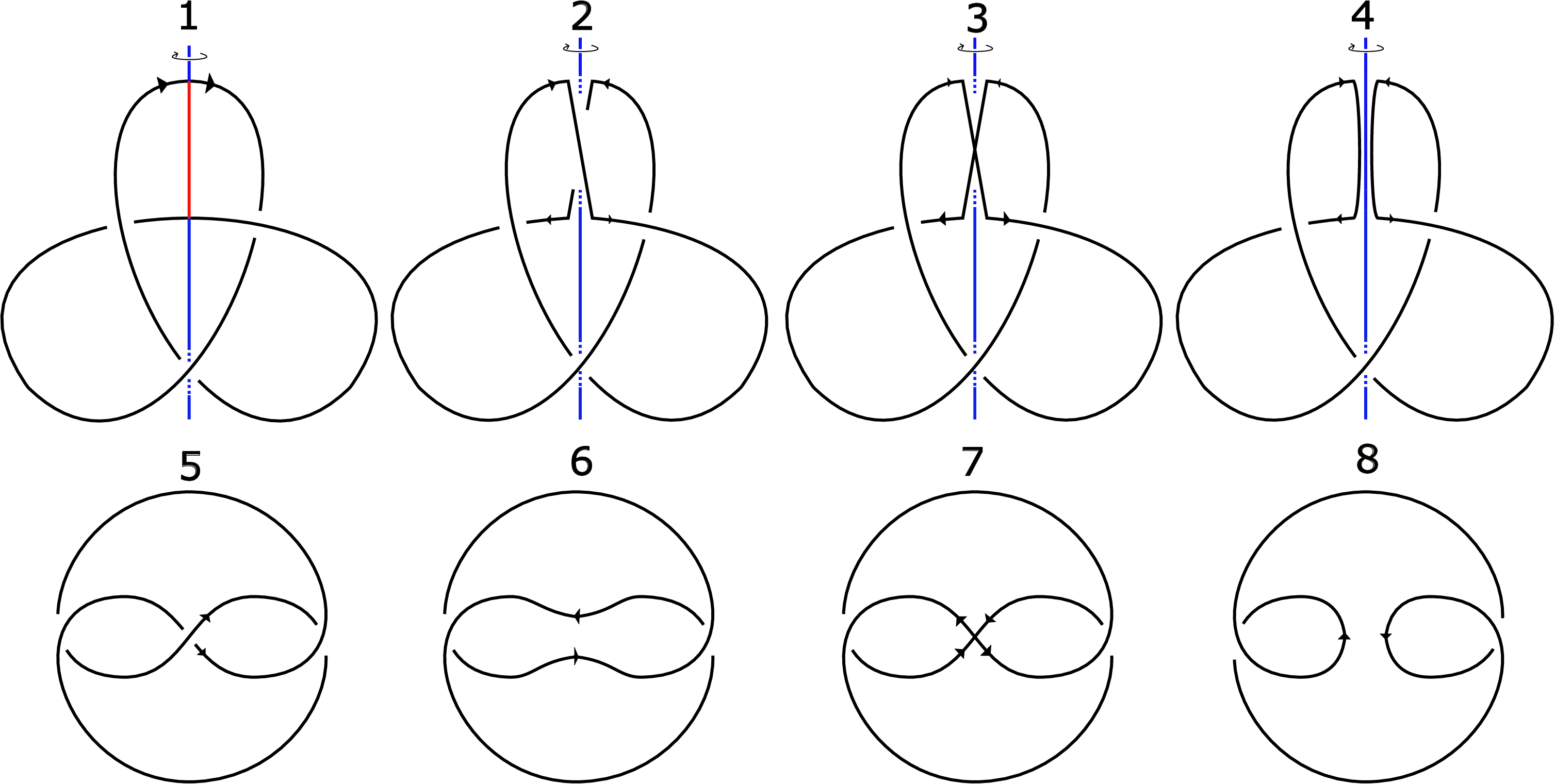}
 \caption{(1) A transvergent diagram showing the strong inversion on left handed trefoil $3_1$. The symmetry axis is split into two halves corresponding to the DSIs $3_1^+$ in red and $3_1^-$ in blue. (2) A transvergent diagram for $L_b^1(3_1^+)$, the positive Hopf link. (3) A transvergent diagram for the singular link $S_b^1(3_1^+)$, as in Definition \ref{Snb}, obtained from singularizing $L_b^1(3_1^+)$. (4) A transvergent diagram for $L_b^{\frac{3}{2}}(3_1^+)$, the unknot. Diagrams (5)-(8) display intravergent perspectives of (1)-(4) respectively.} \label{4trefoil}
 \end{figure}

There is a family of $2$-periodic links naturally associated to a DSI. 
 \begin{definition}\cite[~Definition 4.1]{BoyAhE4G}\cite[~Definition 2.8]{diprisa2023equivariant}\label{n-butterfly link}
The \emph{n-butterfly link} $L_b^n(K)$  of a DSI $K$ is the two component 2-periodic link with linking number $n$ constructed by performing a band surgery on $K$ along the chosen half-axis. The $0$-butterfly link is \emph{the butterfly link} of $K$  and is denoted $L_b(K)\coloneq L_b^0(K)$. 
\end{definition}
See the left hand side of Figure \ref{LbSb} for an illustration of the butterfly link $L_b(3_1^+)$ and $(2)$ of Figure \ref{4trefoil} for an illustration of the $1$-butterfly link $L_b^1(3_1^+)$. We also make use of the following equivariant concordance invariant in the statement of Theorem \ref{Mainthm2}.
\begin{definition}\cite[~Definition 4.6]{BoyAhE4G} The \emph{axis linking number} $\widetilde{lk}(K)$ is the linking number of either component of $L_b(K)$ with the symmetry axis $A$.
\end{definition}
It is natural to consider as well the family of doubly periodic \emph{knots} obtained from band surgery along the axis of $K$.
\begin{definition}
    Given $n\in \Z$, let $L^{n+\frac{1}{2}}_b(K)$ be the doubly periodic knot obtained by performing a band surgery on $K$ along the chosen half axis so that \begin{equation}\label{linkingnumberequation}lk(L_b^{n+\frac{1}{2}}(K), A) = 2\widetilde{lk}(K)+2n+1. \end{equation}
\end{definition}
See part $(4)$ of Figure \ref{4trefoil} for an example of this construction with $K=3_1^+$ and $n=\frac{3}{2}.$ The quotient knot of $L_b^n(K)$ (for integer or half integer value of $n$) -- is the same quotient knot $\overline{K}$ of the DSI $K$. 

\subsection{Singular links and spatial graphs}
Because of the natural relationship of DSIs and trivalent graphs or equivalently singular knots via appending the distinguished half axis, we will also be interested in a family of singular links $S^n_b(K)$.

\begin{figure}
    \centering
    \includegraphics[width=0.5\linewidth]{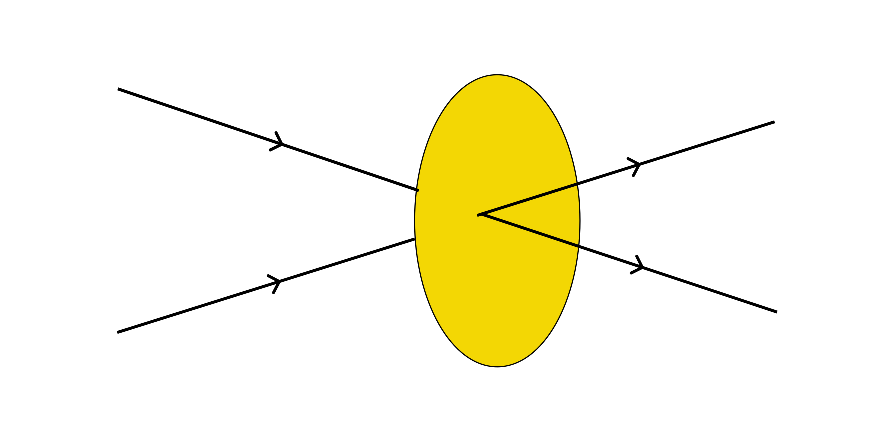}
    \caption{A disk separating the incoming and outgoing edges at a vertex.}
    \label{tdiskv_fig}
\end{figure}

The following definitions are from \cite[~Section 2]{Harvey_2017}.
We will work in the PL-category. A \textit{graph}, $Y$, is a 1-dimensional complex, consisting of a finite set of vertices ($0$-simplices) and edges ($1$-simplices) between them.
A \textit{spatial graph} is an embedding $f$ of a graph $Y$ in $S^3$, $f\colon Y\rightarrow S^3$.  
An \textit{oriented graph} is a graph together with orientations given on each of the edges.  
Let $\mathcal{D}$ be the 2-complex obtained by gluing $3$ copies of the 2-simplex $[e_0,e_1,e_2]$ together so that their union is a disk and with all three of the $e_0$s identified to a single point in $\mathcal{D}$. Note that $e_0$ is the unique vertex in the interior of $\mathcal{D}$. We say that $\mathcal{D}$ is a \textit{standard disk} and $e_0$ is the \textit{vertex associated to $\mathcal{D}$}.
 An \textit{oriented disk graph, $G$,} is a $2$-complex constructed as follows.  Start with an oriented graph $Y$.  Then, for each vertex $v$ of $Y$, glue a standard disk $\mathcal{D}$ to $Y$ by identifying the vertex associated to $\mathcal{D}$ with $v$.  
A \textit{transverse spatial graph} is an embedding $f\colon  G \rightarrow S^3$ of an oriented disk graph $G$ into $S^3$ where each vertex of the graph locally looks like Figure~\ref{tdiskv_fig} and where each standard disk of $G$ lies in a plane.  

A \emph{diagram of a spatial graph} is a projection of $f(Y)$ to $S^2$ with only transverse double points away from vertices, where the over and under crossings are indicated. 
A \textit{regular projection of a transverse spatial graph $f\colon  G \rightarrow S^3$} is a projection that satisfies the following two conditions.  (1) For each point, $pt$, in the image of the underlying graph of $G$, $f^{-1}(pt)$ contains no more than two points and if $f^{-1}(pt)$ contains two points then neither is a graph vertex. (2) All of the standard disks of $f$ are perpendicular to the plane of projection. A \emph{diagram for a transverse spatial graph} $f \colon  G \rightarrow S^3$
is a regular projection of $f$
where all the over and under crossings are indicated.

\begin{figure}
    \centering
    \includegraphics[width=0.3\linewidth]{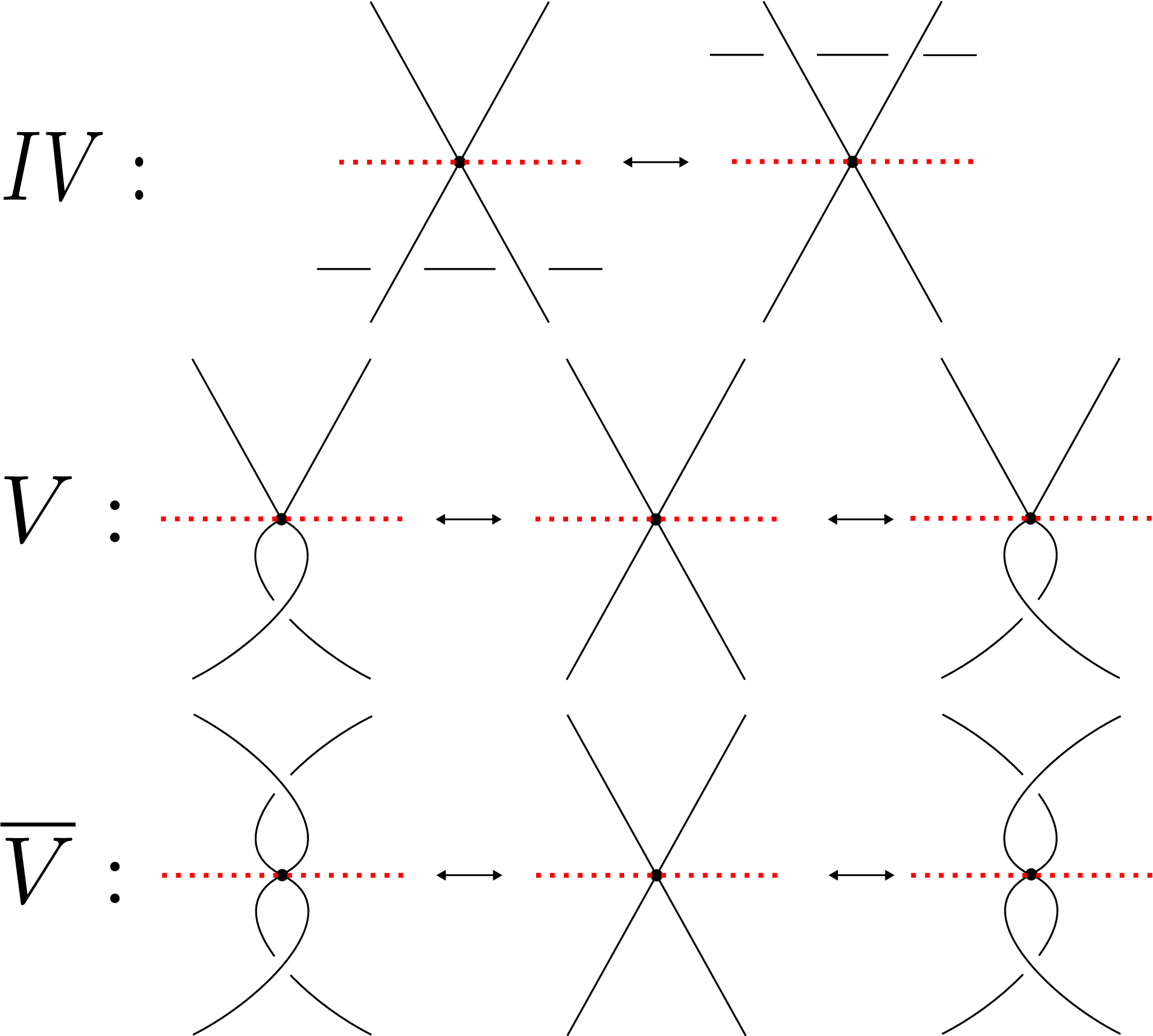}
    \caption{The new Reidemeister moves for transverse spatial graphs. The
dashed line is a disk, transverse to the page, separating incoming and outgoing
edges.}
    \label{Reidemeister moves}
\end{figure}

We will consider two equivalence relations on the set of transverse spatial graphs. 
Transverse spatial graphs $f_1(Y)$ and $f_2(Y)$ are \emph{equivalent as transverse spatial graphs} if diagrams for $f_1(Y)$ and $f_2(Y)$ are related by the standard Reidemeister moves $I-III$, Reidemeister moves $IV$ and $V$ of Figure \ref{Reidemeister moves}, and planar isotopy. Transverse spatial graphs $f_1(Y)$ and $f_2(Y)$ are \emph{equivalent as rigid vertex transverse spatial graphs} if diagrams for $f_1(Y)$ and $f_2(Y)$ are related by the Reidemeister moves I-III, Reidemeister moves $IV$ and $\overline{V}$ of Figure \ref{Reidemeister moves}, and planar isotopy.

\begin{figure}
    \centering
    \includegraphics[width=0.4\linewidth]{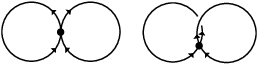}
    \caption{Two transverse spatial graphs which differ by a Reidemeister $V$ move but have distinct Alexander polynomials.}
    \label{fig:RV}
\end{figure}

\begin{definition}
A $4$-valent rigid vertex transverse spatial graph such that each vertex has two incoming and two outgoing edges is a \emph{singular link}. If the singular link has one connected component, then it is a \emph{singular knot}.
\end{definition}
\begin{remark} Ozsv\'{a}th--Stipsicz--Szab\'{o} model singular knots as embedded
oriented trivalent graphs with distinguished thick edges. Contracting each
thick edge gives a four-valent vertex with two incoming and two outgoing
edges, together with the local vertex data encoded here by the separating
disk/rigid-vertex structure. Thus our convention is equivalent to their
singular-knot convention, with the rigid vertex datum made explicit. This
should not be confused with the special class of planar singular knots in
\cite{OSsing1}, namely singular knots admitting an injective projection to the
plane. The singular links considered here need not be planar in that sense. The reason that we use the more restrictive notion of rigid vertex transverse spatial graphs is because the Alexander polynomial -- and hence the absolute Alexander grading in singular knot Floer homology -- of a 4-valent transverse spatial graph is not well defined. For example, consider the two transverse spatial graphs in Figure \ref{fig:RV}. They differ by a $V$ move, but a skein relation calculation using Equation \ref{posskein} shows that the left transverse spatial graph has Alexander polynomial equal to $1$, and the right transverse spatial graph has Alexander polynomial equal to $-t^{-\frac{1}{2}}$.
\end{remark}

\begin{definition}
    Let $D$ be a knot diagram and $c$ be any crossing in $D$.
    The \emph{singularization of} $D$ \emph{at} $c$ is the singular link diagram $S_c(D)$ obtained by performing the local modification seen in Figure \ref{singularize} at $c$.
\end{definition}
    
\begin{figure}[htb!]
\center
\begin{tikzpicture}[scale=0.8]

% Singular Crossing
\draw[->] (3,0) -- (5,2);
\draw[->] (5,0) -- (3,2);
\draw[->] (2.25,1) -- (2.75,1);
\draw[->] (5.75,1) -- (5.25,1);
% Negative Crossing
\draw (6,0) -- (6.75, .75);
\draw[->](7.25, 1.25) -- (8,2);
\draw[->] (8,0) -- (6,2);

%positive crossing
\draw[->] (0,0) -- (2,2);
\draw[dashed][red] (3.5,1) -- (4.5,1);
\draw (2,0) -- (1.25, .75);
\draw[->] (.75, 1.25) -- (0,2);

\end{tikzpicture}
\caption{Singularizing a positive or negative crossing.}\label{singularize}\end{figure}
We are now ready to give the definition of the singularized $n$-butterfly link $S^n_b(K)$.
\begin{definition}\label{Snb}
    Let $D$ be a transvergent diagram for a DSI $K$,  $D_b^n$ be the diagram for $L_b^n(K)$ obtained by surgering in an equivariant band in $D$, and $c$ be a crossing in that band. Then the \emph{singularlized} $n$-\emph{butterfly link of} $K$ is the singular link $S^n_b(K)$ represented by the diagram $S_c(L_b^n(K))$. The \emph{singular butterfly link of} $K$, denoted $S_b(K)$, is the singularized $0$-butterfly link $S_b^0(K)$.
\end{definition}
In Figure \ref{LbSb} the butterfly link $L_b(3_1^+)$ and the singular butterfly link $S_b(3_1^+)$ are depicted. 
\begin{remark}\label{thetaequivclass}
An application of $\overline{V}$ from Figure \ref{Reidemeister moves} shows that the singular knot $S^n_b(K)$ is independent of the choice of the crossing in a diagram for the equivariant band of $L^n_b(K)$. The involutive Reidemeister moves of \cite[~Theorem 2.3]{Lobb_2021} along with the Reidemeister moves for rigid vertex transverse spatial graphs guarantee that $S^n_b(K)$ is independent of chosen transvergent diagram.
The Reidemeister move $V$ in Figure \ref{Reidemeister moves} implies that the transverse spatial graph equivalence class of $S^n_b(K)$ is independent of $n$.
\end{remark}
While it is most convenient to define the $n$-butterfly link and singularized $n$-butterfly link of a DSI from the transvergent perspective, in this paper we primarily examine strongly invertible knots from the intravergent perspective. It is therefore necessary to understand how to take an intravergent diagram $D$ for a DSI $K$, and produce diagrams for $S^n_b(K)$, $L_b^{n}(K)$ and $L_b^{n\pm\frac{1}{2}}(K)$ (where $n$ and the sign $\pm$ depends on $D$) by performing local modifications to the center crossing of $D$.
\begin{definition}\label{intravergentsingularization} Let $c$ be the central crossing of an intravergent diagram $D$ for a strongly invertible knot. The \emph{intravergent singularization of} $D$ at $c$ is the the diagram $IS(D)$ obtained by performing the local modification seen in Figure \ref{intrasing} to $c$, followed by orientations changes on the rest of $D$ that are forced by the new orientations at $c$.

\begin{figure}
    \centering
    \includegraphics[width=0.5\linewidth]{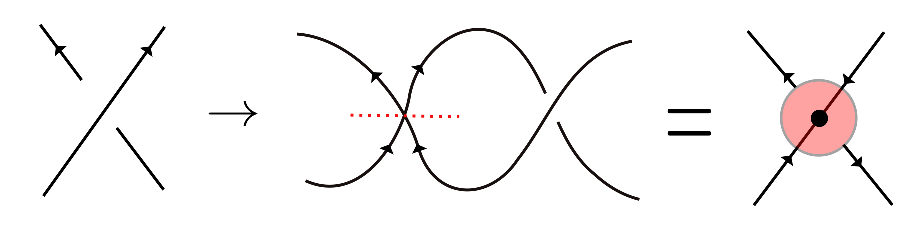}
    \caption{An intravergent singularization of a positive central crossing.}
    \label{intrasing}
\end{figure}
\end{definition}

\begin{remark}Singularization can be done at any crossing on any knot diagram, while we are only defining intravergent singularization for the central crossing in an intravergent diagram of a strongly invertible knot. In general, it is not  possible to orient the rest of the diagram in a way that is compatible with the orientations at the intravergent singularity, but the symmetry of an intravergent diagram for a DSI guarantees that there are global orientations on the diagram compatible with the local ones of Figure \ref{intrasing}.
\end{remark}

It is worth emphasizing the diagrammatic difference between the singularization  and intravergent singularization of a crossing. A transverse double point singularity has two incoming and two outgoing strands. In a singularization, the singularity is presented so that the plane of the diagram is transverse to a disk separating the incoming and outgoing edges. In an intravergent singularization, the singularity is presented with the separating disk \emph{in} the plane of the diagram. From here on, we stop drawing the separating disk of a singular point as its placement is forced by the orientation of the four strands at a singularity.

We may apply an intravergent singularization to an intravergent diagram $D$ to obtain a diagram $IS(D)$ for $S_b^n(K)$. The value of $n$ will depend on the intravergent diagram $D$ that we choose. From now on intravergent diagrams will be labeled with a subscript $n$: if $D_n$ is an intravergent diagram for $K$, then $IS(D_n)$ is a diagram for $S^n_b(K)$. If one applies the equivariant Reidemeister move illustrated in \cite[~Figure 4.1]{KhSI} to the central crossing of a  diagram $D_n$, this has the effect of shifting $n$ by a half integer.
Associated to an intravergent diagram $D_n$ of $K$ is a singular skein triple consisting of $S^n_b(K)$ and two 2-periodic resolutions as shown in Figure \ref{fig:singskeintripintra}. 
\begin{figure}[htb!]
\center
\begin{tikzpicture}[scale=0.8]

% Singular Crossing
\draw[->] (3,0) -- (3.5,.5);
\draw(4,1) -- (3.5,.5);
\draw(4,1) -- (4.5, 1.5);
\draw[->] (5,2) -- (4.5,1.5);
\draw(5,0) -- (4.5, .5);
\draw[<->](4.55,.45) -- (3.45,1.55);
\draw(3,2) -- (3.5, 1.5);

\draw[<-] (2.25,1) -- (2.75,1);
\draw[->] (5.25, 1) -- (5.75, 1);
% Oriented Resolution
\draw[->] (6,0) to [bend right=60] (6,2);
\draw[<-] (8,0) to [bend left=60] (8,2);
\draw[->] (2,2) to [bend left=60] (0,2);
\draw[->] (0,0) to [bend left=60] (2,0);
\end{tikzpicture}
\caption{A singular skein triple from the intravergent perspective.}\label{fig:singskeintripintra}\end{figure}
If the central crossing of $D_n$ is positive this is called a positive singular skein triple, and we denote these 2-periodic resolutions as $S_b^n(K)_+=L_b^n(K)$ and $S_b^n(K)_0=L_b^{n-\frac{1}{2}}(K)$.
If the central crossing of $D_n$ is negative this is called a negative singular skein triple, and we denote these 2-periodic resolutions as $S_b^n(K)_0=L_b^{n+\frac{1}{2}}(K)$ and $S_b^n(K)_-=L_b^n(K)$. These labeling conventions come from looking at Figure \ref{fig:singskeintripintra} from the transvergent perspective; this is illustrated in case of $3_1^+$ in Figure \ref{4trefoil}.\nopagebreak\section{Heegaard Floer homology}\label{HFsection}
In Sections \ref{heegaard} and \ref{link} we give brief reviews of Heegaard Floer homology and link Floer homology, mainly to set up notation. The original references for this standard material are \cite{OS3manifolds1, OSknots, OSlink, RasmussenThesis}. In Section \ref{singularlinkFloer} we review singular link Floer homology \cite{OSsing1,OSsing2}.
\subsection{Heegaard Floer homology}\label{heegaard}
Let $X$ be a integer homology $3$-sphere.
 A \emph{Heegaard splitting} of $X$ is a decomposition $X=U_{\bm{\alpha}}\cup_{\Sigma_g}U_{\bm{\beta}}$ into two genus $g$ handlebodies glued along their common boundary $\Sigma_g$. The handlebodies $U_{\bm{\alpha}}$ and $U_{\bm{\beta}}$ are specified by systems of attaching curves $\bm{\alpha}=\{\alpha_1,\hdots,\alpha_{g+n-1}\}$ and $\bm{\beta}=\{\beta_1,\hdots,\beta_{g+n-1}\}$ on $\Sigma_g$. Given $n$ basepoints $\bm{w}=\{w_1,\hdots,w_n\}\subset \Sigma_g$ such that each connected component of $\Sigma_g-\bm{\alpha}$ and $\Sigma_g-\bm{\beta}$ contains exactly one $w_i$, the data $\mathcal{H}=(\Sigma_g,\bm{\alpha},\bm{\beta},\bm{w})$ defines an $n-$\emph{pointed Heegaard diagram} for $X$.
 
 The connected components of $\Sigma-\cup \alpha_i-\cup \beta_i$ are called \emph{elementary regions}. A linear combination $P$ of elementary regions is called a \emph{periodic domain} if $\partial P$ can be expressed as a linear combination of the $\alpha$ and $\beta$ curves. The Heegaard diagram $\mathcal{H}$ is \emph{weakly admissible} if every periodic domain has both positive and negative local multiplicities. From now on we assume that $\mathcal{H}$ is weakly admissible, as this is a necessary assumption to compute Heegaard Floer homology from $\mathcal{H}.$
 The $\alpha$ and $\beta$ curves define two Lagrangians $\talpha\coloneq \alpha_1\times\hdots\times \alpha_{g+n-1}$ and $\tbeta\coloneq \beta_1\times\hdots\times \beta_{g+n-1}$ in the symplectic manifold $(\Sym^{g+n-1}(\Sigma_g\backslash \bm{w}),\omega_{\bm{\alpha},\bm{\beta}})$ where $\omega_{\bm{\alpha},\bm{\beta}}$ is a symplectic form which agrees with the product symplectic form away from the fat diagonal of the symmetric product as constructed in \cite{perutz}. Then the Heegaard Floer homology of the diagram $\mathcal{H}$ is the Lagrangian Floer homology \begin{equation}\widetilde{HF}(\mathcal{H})\coloneq HF(\talpha,\tbeta)\end{equation} of $(\talpha,\tbeta)$ computed in $\Sym^{g+n-1}(\Sigma_{g}\backslash \bm{w})$. Heegaard Floer homology admits a $\Z$-valued homological (or \emph{Maslov}) grading which we denote by $M$:
 
\begin{equation}
\widetilde{HF}(\mathcal{H})=\bigoplus_{M\in\Z}\widetilde{HF}_M(\mathcal{H})
\end{equation}
The version of Heegaard Floer homology defined above is not quite a three manifold invariant as it depends on the number of basepoints; if $\mathcal{H}'$ is a $k$-pointed Heegaard diagram for $X$ and $\mathcal{H}$ is a singly pointed Heegaard diagram for $X$ then 

\begin{equation}
    \widetilde{HF}(\mathcal{H}')\cong\widetilde{HF}(\mathcal{H})\otimes (\F\oplus\F_{(-1)})^{\otimes (k-1)}.
\end{equation}
where $\F\oplus\F_{(-1)}$ is the two dimensional vector space with generators in Maslov gradings $0$ and $-1.$ The isomorphism class of the homology $\widetilde{HF}(\mathcal{H})$ is independent of the choice of diagram $\mathcal{H}$ for $X$, and hence it is a three manifold invariant denoted $\HFhat(X)$.
\subsection{Link Floer homology}\label{link}
Next we discuss a link invariant due to Ozsv\'{a}th, Szab\'{o} called link Floer homology, whose definition closely mirrors that of Heegaard Floer homology \cite{OSSlink}; the case of knots was also studied independently by Rasmussen \cite{RasmussenThesis}. Given an $\ell$-component link $L\subset S^3$, a $2n$-pointed \emph{Heegaard diagram} for $L$ is a quintuple \[\mathcal{H}=(\Sigma_g, \bm{\alpha}=\{\alpha_1,\hdots,\alpha_{g+n-1}\},\bm{\beta}=\{\beta_1,\hdots,\beta_{g+n-1}\},\bm{w}=\{w_1,\hdots,w_n\},\bm{z}=\{z_1,\hdots,z_n\})\]
such that 
\begin{itemize}
\item The quadruples $(\Sigma_g,\bm{\alpha},\bm{\beta},\bm{w})$ and $(\Sigma_g,\bm{\alpha},\bm{\beta},\bm{z})$ are $n$-pointed Heegaard diagrams for $S^3$.
\item The union of arcs $\xi_i$ connecting each $w_i$ to a $z_j$ on $\Sigma_g-\cup_{i=1}^{g+n-1}\alpha_i$ that are then slightly pushed into the handlebody specified by $\bm{\alpha}$ and arcs $\zeta_i$ connecting each $z_j$ to a $w_i$ on $\Sigma_g-\cup_{i=1}^{g+n-1}\beta_i$ that are then slightly pushed into the handlebody specified by $\bm{\beta}$ yields the link $L$: that is, $\bigcup_i \xi_i\cup \zeta_i=L$.
\end{itemize}
Then the link Floer homology of a weakly admissible $2n$-pointed Heegaard diagram $\mathcal{H}$ is
\begin{equation}\label{linkfloereq}
\widetilde{HFL}(\mathcal{H})\coloneq HF(\talpha,\tbeta)
\end{equation}
where $\talpha\coloneq \alpha_1\times\hdots\times\alpha_{g+n-1}$ and $\tbeta\coloneq \beta_1\times\hdots\times \beta_{g+n-1}$, and the right hand side of the Equation (\ref{linkfloereq}) is the Lagrangian Floer homology $HF(\talpha,\tbeta)$ computed inside of $\Sym^{g+n-1}(\Sigma_g\backslash (\bm{w}\cup\bm{z}))$, again with respect a suitable symplectic form \cite{perutz}. Link Floer homology decomposes as a direct sum along two gradings -- the $\Z$ valued homological/Maslov grading, and the \emph{Alexander grading}.  Write $L=K_1\cup\hdots\cup K_\ell$ for a decomposition of $L$ into its $\ell$ connected components. Define \[a_i\coloneq \begin{cases}
    0&\textrm{ if }\sum_{j\neq i}lk(K_i,K_j)\equiv 0\textrm{ (mod }2\textrm{)}\\
    \frac{1}{2}&\textrm{ if }\sum_{j\neq i}lk(K_i,K_j)\equiv 1\textrm{ (mod }2\textrm{)}
\end{cases}\] and $a\coloneq \sum_{i=1}^\ell a_i.$ Then the Alexander grading takes values in $\Z+a$.
\begin{remark}
    In some variations the Alexander grading of a multi-component link is a multi-grading, but we use the collapsed version of this grading defined by adding up all of the multi-gradings.
\end{remark}Link Floer homology in Maslov and Alexander gradings $(d,s)$ is denoted $\widetilde{HFL}_M(\mathcal{H},A)$:
\begin{equation}
    \widetilde{HFL}(\mathcal{H})=\bigoplus_{(d,s)\in \Z\oplus (\Z+a)} \widetilde{HFL}_M(\mathcal{H},A).
\end{equation}This version of link Floer homology depends on the number basepoints; if we let $V\coloneq \F_{(0,0)}\oplus \F_{(-1,-1)}$ be the two dimensional vector space with generators in gradings $(0,0)$ and $(-1,-1)$, $\mathcal{H}'$ be a $2k$-pointed Heegaard diagram for $L$ for some $k\ge \ell$, and $\mathcal{H}$ be a $2\ell$ pointed Heegaard diagram for $L$, then 
\begin{equation}\label{HFKdependencyonbp}
    \widetilde{HFL}(\mathcal{H}')\cong\widetilde{HFL}(\mathcal{H})\otimes V^{\otimes(k-\ell)}.
\end{equation}
The isomorphism class of the homology $\widetilde{HFL}(\mathcal{H})$ is independent of the choice of $2\ell$ pointed diagram $\mathcal{H}$ for $L$, and hence it is a link invariant denoted $\HFLhat(L)$. An important property of link Floer homology is that it recovers the Alexander polynomial as its bigraded Euler characteristic \cite[~Theorem 1.3]{OSSlink}:
\begin{equation}\label{HFLeulerchar}
    \chi(\HFLhat(L))\coloneq \sum_{M,A}t^A(-1)^M\text{dim}(\HFLhat_M(L,A))=\Delta_L(t)\cdot(t^{\frac{1}{2}}-t^{-\frac{1}{2}})^{\ell-1}.
\end{equation}

\subsection{Singular link Floer homology}\label{singularlinkFloer}
In \cite{OSsing1} and \cite{OSsing2}, 
a variant of Heegaard Floer homology is established for singular links. Here we review this construction. In this section a singular link will refer to a singular link with exactly one transverse double point singularity and one component. The conventions and equations for singular link Floer homology described below reflect this specialization, but can be easily generalized at the cost of additional notation.
\begin{remark}
    We follow the basepoint conventions of \cite{OSsing2} for singular link Floer homology with the understanding that our $\bm{w}$ basepoints are their $\mathbb{O}$ basepoints and our $\bm{z}$ basepoints are their $\mathbb{X}$ basepoints. 
\end{remark}

A \emph{multipointed Heegaard diagram for a singular link} $S$ is a quintuple
\begin{align*}
    \mathcal{H}=&\\&(\Sigma_g,\bm{\alpha}=\{\alpha_1,\hdots,\alpha_{g+n-1}\},\bm{\beta}=\{\beta_1,\hdots,\beta_{g+n-1}\},\bm{w}=\{w_1,\hdots,w_{n+1}\},\bm{z}=\{z_1^*,z_2,\hdots,z_{n}\})
\end{align*}
such that 
\begin{itemize}
\item The sets of attaching curves $\bm{\alpha}$ and $\bm{\beta}$ specify a handlebody decomposition of $S^3.$
\item In each region $R$ of $\Sigma_g-\cup_{i=1}^{g+n-1}\alpha_i$ or $\Sigma_g-\cup_{i=1}^{g+n-1}\beta_i$, $|
R\cap \bm{z}|=1$ and either $z_1^*\not\in R\cap\bm{z}$ and $|R\cap \bm{w}|=1$, or $z_1^*\in R\cap \bm{z}$ and $|R\cap \bm{w}|=2$.
\item If one connects $w$'s to $z$'s by arcs $\xi_i$ on $\Sigma_g-\cup_{i=1}^{g+n-1}\alpha_i$ that are then slightly pushed into the $\bm{\alpha}$ handlebody, and connects $z$'s to $w$'s by arcs $\zeta_i$ on $\Sigma_g-\cup_{i=1}^{g+n-1}\beta_i$ that are then slightly pushed into the $\bm{\beta}$ handlebody, then the union of the $\xi$ and $\zeta_i$ recovers the singular link $S$: \[(\bigcup \xi_i)\cup(\bigcup \zeta_i)=S.\] The point $z_1^*$ will connect to two distinct $w_i$.
\end{itemize}
In the final condition, the equality is meant in the category of
rigid-vertex transverse spatial graphs. More explicitly, the special
basepoint \(z_1^*\) is the singular point. The two arcs in the
\(\alpha\)-handlebody incident to \(z_1^*\) are oriented from \(w\)-basepoints
to \(z_1^*\), while the two arcs in the \(\beta\)-handlebody incident to
\(z_1^*\) are oriented from \(z_1^*\) to \(w\)-basepoints. Thus a sufficiently
small disk \(D_*\subset \Sigma_g\) centered at \(z_1^*\), disjoint from the
attaching curves and from the other basepoints, gives the standard disk at
the vertex: it separates the two incoming \(\alpha\)-handlebody half-edges
from the two outgoing \(\beta\)-handlebody half-edges. In this way the
Heegaard diagram determines the rigid-vertex structure, not merely the
underlying four-valent graph.

The singular link Floer homology of a weakly admissible Heegaard diagram $\mathcal{H}$ is the Lagrangian intersection Floer homology
\begin{equation}
    \widetilde{HFK}(\mathcal{H})\coloneq HF(\talpha,\tbeta)
\end{equation}
where $\talpha\coloneq \alpha_1\times\hdots\times \alpha_{g+n-1},$ $\tbeta\coloneq \beta_1\times\hdots\times \beta_{g+n-1}\subset \Sym^{g+n-1}(\Sigma_g\backslash(\bm{w}\cup\bm{z}))$ and the ambient symplectic manifold is $\Sym^{g+n-1}(\Sigma_g\backslash(\bm{w}\cup\bm{z}))$ equipped with an appropriate symplectic form \cite{perutz}.  Singular link Floer homology also decomposes as a bigraded sum on a Maslov (homological) and Alexander grading:
\begin{equation}
    \widetilde{HFK}(\mathcal{H})=\bigoplus_{(M,A)\in\Z\oplus (\Z+a)}\widetilde{HFK}_M(\mathcal{H},A)
\end{equation}
where $a=\frac{1}{2}$ if resolving the singularity of $S$ yields a $2$ component link and $a=0$ if the zero resolution (cf. Figure \ref{fig:negskeintriple}) of the singularity of $S$ yields a knot. Singular link Floer homology depends on the number of basepoints in the Heegaard diagram for $S$. If $\mathcal{H}'$ is a Heegaard diagram for $S$ with $|\bm{z}|=|\bm{w}|-1=k$ and $\mathcal{H}$ is a Heegaard diagram for $S$ with $|\bm{z}|=|\bm{w}|-1=1$ then
\begin{equation}
    \widetilde{HFK}(\mathcal{H}')\cong\widetilde{HFK}(\mathcal{H})\otimes V^{\otimes(k-1)}
\end{equation}
where $V=\F_{(0,0)}\oplus\F_{(-1,-1)}$ is as in the last section.
The homology $\widetilde{HFK}(\mathcal{H})$ is independent of the choice of triply pointed Heegaard diagram for $S$, and hence it is a singular link invariant denoted $\widehat{HFK}(S)$. Singular links also have Alexander polynomials, characterized by either of the skein relations
\begin{equation}\label{posskein}
    \Delta_{S_+}(t)=\Delta_{S}(t)+t^{\frac{1}{2}}\Delta_{S_0}(t)
\end{equation}
or
\begin{equation}\label{negskein}
    \Delta_{S_-}(t)=\Delta_{S}(t)+t^{-\frac{1}{2}}\Delta_{S_0}(t)
\end{equation}
where $S$, $S_-$, $S_+$ and $S_0$ are as in Figure \ref{fig:negskeintriple}. The Alexander polynomials of $S_+$, $S_-$ and $S_0$ are inductively defined, as each of these potentially singular links will have one fewer singularities than $S.$

\begin{figure}[htb!]
\center
\begin{tikzpicture}[scale=0.8]

% Singular Crossing
\draw[->] (0,0) -- (2,2);
\draw[->] (2,0) -- (0,2);

% Negative Crossing
\draw (3,0) -- (3.75, .75);
\draw[->](4.25, 1.25) -- (5,2);
\draw[->] (5,0) -- (3,2);
\draw[->] (6,0) -- (8,2);
\draw (8,0) -- (7.25,.75);
\draw[->] (6.75,1.25) -- (6, 2);

% Oriented Resolution
\draw[->] (9,0) to [bend right=60] (9,2);
\draw[->] (11,0) to [bend left=60] (11,2);

\end{tikzpicture}
\caption{From left to right: The singular crossing $S$, its negative resolution $S_-$, its positive resolution $S_+$ and its zero resolution $S_0$.}\label{fig:negskeintriple}\end{figure}
Singular link Floer homology recovers the Alexander polynomial of a singular link as its bi-graded Euler characteristic \cite[~Theorem 1.1]{OSsing1}:
\begin{equation}
\chi(\widehat{HFK}(S))=\sum_{M,A}t^A\textrm{dim}(\widehat{HFK}_M(S,A))=\Delta_{S}(t).
\end{equation}

 \section{Equivariant Heegaard diagrams}\label{diagrams}
\subsection{Equivariant Heegaard diagrams for strongly invertible knots}
Let $D_n$ be a $2g-1$ crossing intravergent diagram for a DSI $K$. Recall that this notation means that the intravergent singularization of $D_n$, $IS(D_n)$, is a diagram for the singularized n-butterfly link $S^n_b(K)$. In this section we construct a doubly pointed genus $2g$ $\tau$-equivariant diagram $\mathcal{H}(D_n)$ for $K$, and a triply pointed genus $2g$ $\tau$-equivariant diagram $\mathcal{H}(IS(D_n))$ for $S^n_b(K)$. We also describe the quotients by the $\tau$ action of these Heegaard diagrams.
 
 \begin{figure}
     \centering
     \includegraphics[width=0.5\linewidth]{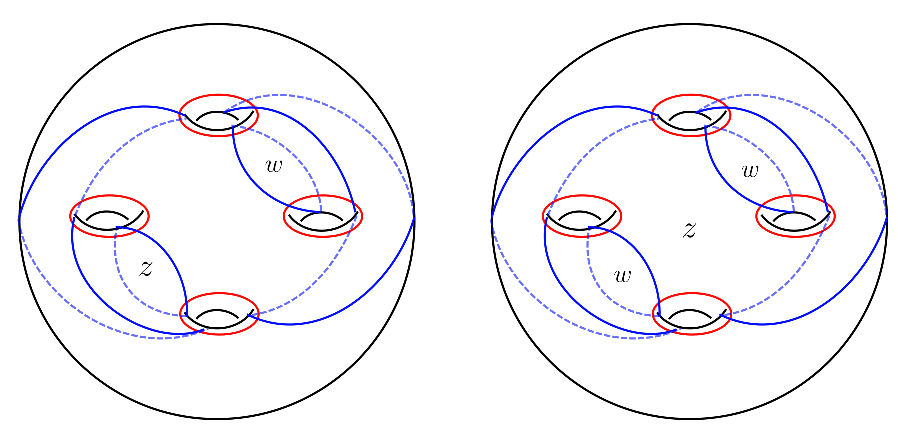}
     \caption{The diagrams $\mathcal{H}(D_n)$ (left) and $\mathcal{H}_S(D_n)$ (right).}
     \label{TrefoilHeegaardDiagram2}
 \end{figure}
\subsubsection{A doubly pointed Heegaard diagram for a strongly invertible knot}\label{twopointdiagram} 
\begin{figure}[htb!]
 \center
 \includegraphics[scale=1.4]{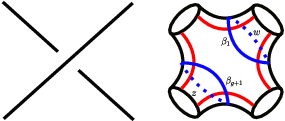}

 \includegraphics[scale=1.4]{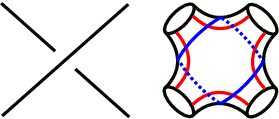}
 \caption{Top: The placement of $\beta$ curves and basepoints for the Heegaard diagram $\mathcal{H}(D_n)$ near the central crossing. Bottom: The placement of $\beta$ curves for the Heegaard diagrams $\mathcal{H}(D_n)$ and $\mathcal{H}_S(D_n)$ near a non-central crossing.} \label{non-sing}
 \end{figure}

Briefly, the doubly pointed Heegaard diagram $\mathcal{H}(D_n)$ for $K$ is the same as the ``pringle chip'' Heegaard diagram in \cite[~Proposition 12.1]{OSsurvey}, except at the central crossing where we instead use the placement of basepoints and $\beta$ curves seen in Figure \ref{non-sing}. We now explain this in more detail. Ignore the crossing data on $D_n$ to obtain an immersed curve in the plane $C$ with $2g-1$ double points. Let $\Sigma$ be the boundary of a  $\tau$-equivariant regular neighborhood of $C$ in $S^3.$ For each of the $2g$ bounded regions $R_1,\hdots,R_{2g}$ in the complement of $C$, put attaching curves $\bm{\alpha}=\{\alpha_1=R_1\cap \Sigma,\hdots,\alpha_{2g}=R_{2g}\cap \Sigma\}$ on $\Sigma.$ Pick the labels on the $\alpha$ curves so that $\tau(\alpha_i)=\alpha_{g+i}$ for $i=1,\hdots,g.$ 
 Corresponding to the central crossing in $D_n$, add two meridians for $K$ as attaching curves $\beta_1$ and $\beta_{g+1}$, and place two basepoints $w$ and $z$ in such a way that $\tau(w)=z$ as indicated in the top part of Figure \ref{non-sing}.  Notice that $\tau(\beta_1)=\beta_{g+1}.$ For each of the non-central $2g-2$ crossings in $D_n$, insert a $\beta$ curve on $\Sigma$ as shown in the bottom part of Figure \ref{non-sing} labeled as $\beta_2,\hdots, \beta_g, \beta_{g+2}, \hdots, \beta_{2g}$ in such a way that $\tau(\beta_i)=\beta_{g+i}$ for $i=2,\hdots,g$. Then we define \[\mathcal{H}(D_n)\coloneq (\Sigma, \bm{\alpha}=\{\alpha_1, \hdots, \alpha_{2g}\}, \bm{\beta}=\{\beta_1,\hdots,\beta_{2g}\}, \bm{w}=\{w\}, \bm{z}=\{z\}).\]

An example of $\mathcal{H}(D_n)$ for $D_n$ equal to intravergent diagram (5) of Figure \ref{4trefoil} is shown on the left side of Figure \ref{TrefoilHeegaardDiagram2}.

\subsubsection{A triply pointed Heegaard diagram for the singularized n-butterfly link}\label{threepointdiagram}
\begin{figure}[htb!]
 \center
 \includegraphics[scale=.5]{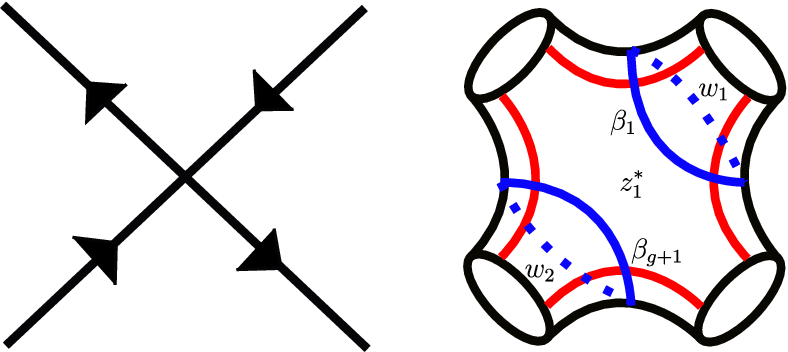}
 \caption{The placement of $\beta$ curves and basepoints for the Heegaard diagram $\mathcal{H}_S(D_n)$ near the central singular crossing.}\label{heegsingularcenter}
 \end{figure}

The triply pointed Heegaard diagram $\mathcal{H}_S(D_n)$ for $S^n_b(K)$ has the same underlying surface, $\Sigma$, and the same set of attaching curves, $\alpha$ and $\beta$, as the Heegaard diagram $\mathcal{H}(D_n)$. The difference is the number, placement and labeling of basepoints: near the central singular crossing of $IS(D_n)$, we have the corresponding piece of the Heegaard diagram shown in Figure \ref{heegsingularcenter}.

In summary, we define
\[\mathcal{H}_S(D_n)\coloneq (\Sigma, \bm{\alpha}=\{\alpha_1, \hdots, \alpha_{2g}\}, \bm{\beta}=\{\beta_1,\hdots,\beta_{2g}\}, \bm{w}=\{w_1, w_2\}, \bm{z}=\{z^*\}).\]

An example of a Heegaard diagram thus constructed from $D_n$ equal to the intravergent diagram (t) of Figure \ref{4trefoil} is shown on the right of Figure \ref{TrefoilHeegaardDiagram2}.

\subsubsection{Quotient diagrams}
Recall the quotient map $\ell\colon S^3\to S^3/\tau$. Let $\overline{\Sigma}=\ell(\Sigma)$ be the genus $g$ quotient of $\Sigma$ by $\tau$. Define $\overline{w}\coloneq \ell(w)=\ell(z)$, $\overline{\bm{\alpha}}=\ell(\bm{\alpha})=\{\overline{\alpha}_1=\ell(\alpha_1),\hdots,\overline{\alpha}_g=\ell(\alpha_g)\}$ and $\overline{\bm{\beta}}=\ell(\bm{\beta})=\{\overline{\beta}_1=\ell(\beta_1),\hdots,\overline{\beta}_g=\ell(\beta_g)\}$. Then the quotient of $\mathcal{H}(D_n)$ by $\tau$ is the genus $g$ singly pointed Heegaard diagram \[\mathcal{H}(D_n)/\tau\coloneq (\overline{\Sigma}, \overline{\bm{\alpha}}, \overline{\bm{\beta}}, \overline{w}).\] 
The quotient diagram $\mathcal{H}_S(D_n)/\tau=(\overline{\Sigma}, \overline{\bm{\alpha}}, \overline{\bm{\beta}}, \overline{w}, \overline{z}=\ell(z^*))$ is the doubly pointed Heegaard diagram of Proposition 12.1 in \cite{OSsurvey} for the quotient diagram $\overline{D}_n$ of $\overline{K}$. The Heegaard diagram $\mathcal{H}_S(D_n)/\tau$ is obtained from $\mathcal{H}(D_n)/\tau$ by adding in a $z$ basepoint, and hence $\mathcal{H}(D_n)/\tau$ is a Heegaard diagram for $S^3.$

\subsubsection{Equivariant Heegaard diagrams of intravergent type}
Let $(K,\tau,A_1)$ be a DSI with symmetry axis $A=A_0\cup A_1$. The following definitions and lemma will be useful when proving invariance of both spectral sequences.

\begin{definition}Let  \[\mathcal{H}=(\Sigma,\bm{\alpha},\bm{\beta}, \bm{w}, \bm{z})\] be a Heegaard diagram for $K$  with 
\[\tau(\Sigma)=\Sigma\textrm{, }
    \tau( \bm{\alpha})=\bm{\alpha}\textrm{, }
   \tau(\bm{\beta})=\bm{\beta}\textrm{, }
    \tau(\bm{w})=\bm{w}\textrm{, }
    \tau(\bm{z})=\bm{z}\textrm{, and }  
    |A\cap \Sigma|=2\]
such that $\tau$ swaps attaching curves in pairs (i.e there are no attaching curves fixed by $\tau$ set-wise). Let $f$ be a self-indexing $\tau$-equivariant Morse function compatible with $\mathcal{H}$. Then the gradient flow-lines through the two points in $A\cap \Sigma$ are $A_0$ and $A_1$ up to orientation. See the proof of Lemma \ref{eqhamiso} for a justification of the last two sentences. Assume that the flow-line through the fixed point corresponding to $A_1$ is oriented the same way as the direction on $A_1$ coming from the DSI. We call this the ``top fixed point'', and the other fixed point the ``bottom fixed point.'' Then $\mathcal{H}$ is a $\tau$\emph{-equivariant Heegaard diagram for} $(K,\tau,A_1)$ \emph{of intravergent type}. 
\end{definition}

The diagram $\mathcal{H}(D_n)$ introduced in a previous subsection is a $\tau$-equivariant Heegaard diagram of intravergent type. For the remainder of this paper we write ``equivariant Heegaard diagram'' in lieu of ``$\tau$-equivariant Heegaard diagram of intravergent type."
\begin{lemma}\label{eqhamiso}
Let $\mathcal{H}$ and $\mathcal{H}'$ be equivariant Heegaard diagrams for the same DSI $K$ with the same number of basepoints. Then there is a sequence of handleslides, isotopies and type 1-2 (de)stabilizations (i.e connect sums with the standard genus 1 Heegaard diagram for $S^3$) that interpolate $\mathcal{H}/\tau$ and $\mathcal{H}'/\tau$ in the complement of (the image under $q\colon S^3\to S^3/\tau$ of) the top and bottom fixed points. 
\end{lemma}
\begin{proof}
Say that $\mathcal{H}=(\Sigma,\bm{\alpha},\bm{\beta},\bm{w},\bm{z})$. Consider the Heegaard diagram $s\mathcal{H}$ constructed by adding a basepoint $z^*$ at the top fixed point and a basepoint $w^*$ at the bottom fixed point. This is a Heegaard diagram for the singular knot $K\cup A$ which has two singularities at $K\cap A$; the flow-lines through the basepoints $\bm{w}\cup\bm{z}$ union to the knot $K$ by assumption, and we can see that the flow-lines through the $w^*$ and $z^*$ union to $A$ as follows. Let $f\colon S^3\to[0,3]$ be a self-indexing $\tau$-invariant Morse function for which $f^{-1}(\frac{3}{2})=\Sigma,$ and such that there exists a $\tau$-invariant metric $g$ so that the ascending manifolds of the index $0$ and $1$ critical points intersect $\Sigma$ in $\bm{w}\cup\{w^*\}$ and $\bm{\alpha}$ respectively, and the descending manifolds of the index $3$ and $2$ critical points intersect $\Sigma$ in $\bm{z}\cup\{z^*\}$ and $\bm{\beta}$ respectively.  A flow-line of the $\tau$ invariant Morse function $f$ under the $\tau$ invariant metric $g$ either is a subset of $A=\textrm{Fix}(\tau)$ or is disjoint from it: if $\gamma $ is a flow-line, then so is $\tau(\gamma)$, and if they intersect at some time on $A$, then by uniqueness of solutions to ODEs it must be the case that $\gamma=\tau(\gamma)$ at all times. It is clear then that the union of flow-lines through $w^*$ and $z^*$ is a subset of $A$, but this subset must actually be all of $A$ since there is only one $\tau$ invariant critical point of index $0$ and one $\tau$ invariant critical point of index $3$. 

We must then have that $s\mathcal{H}/\tau$ is a Heegaard diagram for $\overline{K}\cup \overline{A}$. The same argument for $s\mathcal{H}'/\tau$ shows that it is a Heegaard diagram for $\overline{K}\cup \overline{A}$ as well. Therefore, by \cite[~Proof of Theorem 2.4]{OSsing1}, $s\mathcal{H}/\tau$ and $s\mathcal{H}'/\tau$ are connected by a sequence of pointed handleslides, isotopies and type 1-2 (de)stabilizations. This is of course equivalent to a sequence of pointed handleslides, isotopies and type 1-2 (de)stabilizations for $\mathcal{H}/\tau$ and $\mathcal{H}'/\tau$ which avoids the image under $q$ of the top and bottom fixed points.
\end{proof}

 \subsubsection{The action on gradings} 
 The strong inversion $\tau$ on $K$ induces $\Z/2\Z$ actions on the equivariant Heegaard diagrams $\mathcal{H}(D_n)$ and $\mathcal{H}_S(D_n)$ which in turn induce $\Z/2\Z$ actions \[\tau_\#\colon \widehat{CFK}(\mathcal{H}(D_n))\to \widehat{CFK}(\mathcal{H}(D_n))\]and\[\tau_\#\colon \widehat{CFK}(\mathcal{H}_S(D_n))\to \widehat{CFK}(\mathcal{H}_S(D_n)).\]
 
 \begin{proposition}\label{tauonknot}
     The $\tau_\#$ action on $\widehat{CFK}(\mathcal{H}(D_n))$ has the following properties:
     \begin{enumerate}
         \item For a homogeneously graded element $x\in \widehat{CFK}(\mathcal{H}(D_n))$, $A(x)=-A(\tau(x))$
         
         \item If $y$ is a $\tau_\#$-equivariant element of $\widehat{CFK}(\mathcal{H}(D_n))$, $\tau_\#(y)=y$, then $A(y)=0.$
         \item If $A(x)=0$, then $M(x)=M(\tau(x))$.
        
     \end{enumerate}
 \end{proposition}
 \begin{proof}
 Let $x$ and $y$ be any two homogeneous elements of $\widehat{CFK}(\mathcal{H}(D_n))$. Then if $\phi$ is a Whitney disk connecting $x$ and $y,$ $\tau(\phi)$ is a Whitney disk connecting $\tau_\#(x)$ and $\tau_\#(y)$. Furthermore, since $\tau$ switches the $z$ and $w$ basepoints of $\mathcal{H}(D_n)$, we have that \[A(x)-A(y)=n_z(\phi)-n_w(\phi)=n_w(\tau(\phi))-n_z(\tau(\phi))=A(\tau_\#(y))-A(\tau_\#(x))\]and therefore \[A(x)+A(\tau_\#(x))=A(y)+A(\tau_\#(y)).\]Therefore $A(x)+A(\tau_\#(x))$ is a constant independent of the choice of $x$. The symmetry $\HFKhat(K,A)\cong \HFKhat(K,-A)$ forces this constant to be equal to $0$. The second claim follows immediately from the first. Let $y$ be a $\tau$-equivariant generator, such as the lift of any generator in $\widehat{CFK}(\mathcal{H}(D_n)/\tau)$ to $\widehat{CFK}(\mathcal{H}(D_n))$, and let $\phi\in\pi_2(x,y)$. Then since $A(x)=A(y)=0$, \[n_z(\phi)=n_w(\phi)=n_z(\tau(\phi))=n_w(\tau(\phi)).\]Therefore, \[M(x)-M(y)=\mu(\phi)-2n_w(\phi)=\mu(\tau(\phi))-2n_w(\tau(\phi))=M(\tau_\#(x))-M(\tau_\#(y))=M(\tau_\#(x))-M(y)\]which proves the third claim.
 \end{proof}
We have the analogous proposition for the $\tau_\#$ action on $\widehat{CFK}(\mathcal{H}_S(D_n))$. 
\begin{proposition}\label{tauonsing}
For a homogeneous element $x\in \widehat{CFK}(\mathcal{H}_S(D_n))$, \[A(x)=A(\tau_\#(x))\]and \[M(x)=M(\tau_\#(x)).\]
\end{proposition}
\begin{proof}
 On the Heegaard diagram $\mathcal{H}_S(D_n)$ of $S^n_b(K)$, $\tau$ preserves the $\mathbf{w}$ and $\mathbf{z}$ sets of $\mathcal{H}_S(D_n)$. Therefore, if $y\in \widehat{CFK}(\mathcal{H}_S(D_n))$ is a $\tau$-equivariant generator and $\phi\in\pi_2(x,y)$, we see that 
 \[A(x)-A(y)=n_z(\phi)-n_w(\phi)=n_z(\tau(\phi))-n_w(\tau(\phi))=A(\tau_\#(x))-A(\tau_\#(y))=A(\tau_\#(x))-A(y)\]
 and
 \[M(x)-M(y)=\mu(\phi)-2n_w(\phi)=\mu(\tau(\phi))-2n_w(\tau(\phi))=M(\tau_\#(x))-M(\tau_\#(y))=M(\tau_\#(x))-M(y).\]
\end{proof}
\section{Spectral sequences for strongly invertible knots}\label{equivariant}

Subection \ref{localization} lays out hypotheses which, when met, allow for the application of Large's localization theorem for Lagrangian Floer homology, recapped in Theorem \ref{Largelocal}. In Section \ref{mainthmsection}, we apply Theorem \ref{Largelocal} \emph{under the assumption that its technical hypotheses are met} to derive Theorems \ref{Mainthm1} and \ref{Mainthm2}. Sections \ref{geosym} and \ref{stableTNiso} are devoted to demonstrating that the hypotheses of Theorem \ref{Largelocal} are met in our setup.
\subsection{Polarization data and Large's localization theorem}\label{localization}
The following definitions and theorem originate from Section 3.2 of Large's paper \cite{Large}.
\begin{definition}
Let $(M,L_0,L_1$) be a symplectic manifold equipped with two Lagrangian submanifolds.
A set of \emph{polarization data for} $(M, L_0, L_1)$ is a triple $\mathfrak{p} = (E, F_0, F_1)$ where $E\twoheadrightarrow M$ is a symplectic vector bundle and $F_i$ is a Lagrangian subbundle of $E|_{L_i}$ for $i = 0, 1.$
\end{definition}
Letting $\underline{\C}$ and $\underline{\R}$ denote the trivial dimension 1 complex and real vector bundles over $M$, we may \emph{stabilize} polarization data $\mathfrak{p}=(E,F_0,F_1)$ by direct summing with $(\underline{\C},\underline{\R},i\underline{\R})$ to obtain \[\mathfrak{p}\oplus (\underline{\C},\underline{\R},i\underline{\R})\coloneq (E\oplus \underline{\C},F_0\oplus \underline{\R},F_1\oplus i\underline{\R}).\]
\begin{definition}\label{isopol}
    Let $\mathfrak{p}=(E,F_0,F_1)$ and $\mathfrak{p}'=(E',F_0',F_1')$ be two sets of polarization data for $(M,L_0,L_1)$. An \emph{isomorphism of polarization data} is a symplectic vector bundle isomorphism $\alpha\colon E\to E'$ such that $\alpha(F_i)$ is homotopic through Lagrangian subbundles to $F_i'$ for $i=0,1.$ We say that polarization data $\mathfrak{p}$ and $\mathfrak{p}'$ are \emph{stably isomorphic} if $\mathfrak{p}\oplus (\underline{\C},\underline{\R},i\underline{\R})^{\oplus k}$ and $\mathfrak{p}'\oplus (\underline{\C},\underline{\R},i\underline{\R})^{\oplus k'}$ are isomorphic polarization data for some $k,k'\in \N.$
\end{definition}

\begin{definition}
   Let $M$ be a symplectic manifold equipped with two Lagrangians $L_0$ and $L_1$. Suppose $\tau\colon (M,L_0,L_1)\to(M,L_0,L_1)$ is a symplectic involution with $\tau-$invariant sets \[(\M,\Lo,\Li)\subset(M,L_0,L_1).\] Define the \emph{tangent polarization} \begin{equation}
    \mathfrak{p}_T\coloneq (T\M,T\Lo,T\Li)
\end{equation} and the \emph{normal polarization} \begin{equation}\mathfrak{p}_N\coloneq (N\M,N\Lo,N\Li).\end{equation} Then a \emph{stable tangent normal isomorphism} is a stable isomorphism of the polarization data $\mathfrak{p}_T$ and $\mathfrak{p}_N$.
\end{definition}
The following theorem is the main technical tool that we use to establish the spectral sequences of Theorems \ref{Mainthm1} and \ref{Mainthm2}.
\begin{theorem}\label{Largelocal}
    Suppose that 
    \begin{enumerate}
        \item $M$ is an exact symplectic manifold and convex at infinity, and $L_0$, $L_1$ are compact exact Lagrangians.
        \item There is a symplectic involution $\tau\colon (M,L_0,L_1)\to(M,L_0,L_1)$ and an associated stable tangent normal isomorphism from $\mathfrak{p}_T$ to $\mathfrak{p}_N.$

    \end{enumerate}
Then there is a spectral sequence with $E_1$ page equal to $HF(L_0,L_1)\otimes\F[\theta,\theta^{-1}]$ and $E_\infty$ page isomorphic to $HF(\Lo, \Li)\otimes\F[\theta,\theta^{-1}]$. 
\end{theorem}

\subsection{Proofs of main theorems}\label{mainthmsection}

Now let $D_n$ be a $2g-1$ crossing intravergent diagram for the DSI $K$, and recall the genus $2g$ $\tau$-equivariant Heegaard diagrams \[\mathcal{H}(D_n)=(\Sigma,\bm{\alpha},\bm{\beta},z,w) \textrm{ and } \mathcal{H}_S(D_n)=(\Sigma,\bm{\alpha},\bm{\beta},\{w_1,w_2\},z^*),\] along with their genus $g$ quotients by $\tau$,
\[\mathcal{H}(D_n)/\tau=(\overline{\Sigma},\overline{\bm{\alpha}},\overline{\bm{\beta}},\overline{w}) \textrm{ and } \mathcal{H}_S(D_n)/\tau=(\overline{\Sigma},\overline{\bm{\alpha}},\overline{\bm{\beta}},\overline{w},\overline{z})\]constructed in Section \ref{diagrams}.

\subsubsection{Proof of Theorem \ref{Mainthm1}}\label{pfmainthm1}
We work in the symmetric product
\begin{equation}
M\coloneq \Sym^{2g}(\Sigma\backslash\{z,w\})
\end{equation}
equipped with submanifolds
\begin{equation}
L_0\coloneq \mathbb{T}_{\bm{\alpha}}=\alpha_1\times\hdots\times\alpha_{2g}\textrm{ and } L_1\coloneq \mathbb{T}_{\bm{\beta}}=\beta_1\times\hdots\times\beta_{2g}.
\end{equation}

Let $\tau\colon (M,L_0,L_1)\to (M,L_0,L_1)$  be the involution \begin{equation}
    \tau((x_1\hdots x_{2g}))\coloneq (\tau(x_1)\hdots\tau(x_{2g})).
\end{equation}
As explained in Section \ref{HFsection}, the work of Perutz shows that $M$ carries a symplectic form $\omega'\in \Omega^2(M)$ that agrees with the product symplectic form induced by an area form on $\Sigma\backslash\{z\}$ away from the fat diagonal of $\Sym^g(\Sigma\backslash \{z\})$, and with respect to which $L_0$ and $L_1$ are Lagrangians \cite{perutz}. The techniques of \cite[~Section 4]{Hendricks2012rank} -- there applied to the case of symmetric products of multi-punctured spheres corresponding to weakly admissible genus 0 Heegaard diagrams for branched double covers of three manifolds -- generalize immediately to show that $\omega'$ can be modified to a $\tau$-equivariant symplectic form $\omega\in\Omega^2(M)$ with respect to which $M$ is convex at infinity and $L_0$ and $L_1$ are still Lagrangians. The proof given in \cite[~Proposition 4.2]{hendricks2022rank} then shows that there exists a primitive $\lambda$ of $\omega$ with $\lambda|_{L_0}=df_0$ and $\lambda|_{L_1}=df_1$ for some $f_i\in C^\infty(L_i)$ so that $L_0$ and $L_1$ are exact Lagrangians. These remarks tell us that $(1)$ of Theorem \ref{Largelocal} is satisfied for $(M, L_0, L_1, \omega, \tau).$ The proof that $(2)$ of Theorem \ref{Largelocal} is also satisfied in this setup is the subject of Sections \ref{geosym} and \ref{stableTNiso}. From here, we proceed assuming that this assumption is satisfied. 

In order to better understand the fixed point sets $(\M, \Lo, \Li)$, we define 
\begin{equation}\label{M'}
M'=\Sym^g(\overline{\Sigma}\backslash \{\overline{z}\}),
\end{equation}
\begin{equation}\label{L0'L1'}
L_0'\coloneq \mathbb{T}_{\overline{\bm{\alpha}}}=\overline{\alpha}_1\times\hdots\times\overline{\alpha}_g\textrm{ and }
L_1'\coloneq \mathbb{T}_{\overline{\bm{\beta}}}=\overline{\beta}_1\times\hdots\times\overline{\beta}_{g}.
\end{equation}

Recall the quotient map $\ell\colon S^3\to S^3/\tau$ . The map $\phi\colon (M', L_0', L_1')\to (\M, \Lo, \Li)$ defined by 
\begin{equation}
    \phi((\ell(x_1)\hdots \ell(x_g)))=(x_1\tau(x_1)\hdots x_g\tau(x_g))
\end{equation}
identifies $(M', L_0', L_1')$ with $(\M, \Lo, \Li)$; specifically, $\phi$ is a biholomorphic map -- see \cite[~Appendix 1]{Hendricks2012rank} for a proof in charts of this fact. Letting $\phi^*(\omega|_{\M})$ be the symplectic form on $M',$ the Lagrangian intersection Floer homologies $HF(\Lo, \Li)$ and $HF(L_0', L_1')$ are isomorphic. 

Applying Theorem \ref{Largelocal} along with the isomorphisms $HF(\Lo, \Li)\cong HF(L_0', L_1')$ gives us a spectral sequence with $E_1$ page equal to \[HF(L_0,L_1)\otimes \F[\theta,\theta^{-1}]=\HFKhat(K)\otimes\F[\theta,\theta^{-1}]\] and $E_\infty$ page isomorphic to \[HF(L_0', L_1')\otimes \F[\theta,\theta^{-1}]=\HFhat(S^3)\otimes \F[\theta,\theta^{-1}].\]

Now we prove that this spectral sequence is independent of the $\tau$-equivariant Heegaard diagram inducing it, and in particular is independent of the intravergent diagram $D_n$ for $K$ that we start with. The argument is essentially the proof of Corollary 1.10 in \cite{Hendricks_2016}.  By \cite[~Proposition 4.3]{Hendricks_2016} along with the identification of Large's $\Z/2\Z$-equivariant Floer homology with Seidel-Smith's $\Z/2\Z$-equivariant Floer homology, the spectral sequence above is the same as the spectral sequence induced by $CF_{\Z/2}(L_0, L_1)$, the equivariant complex from non-invariant complex structures defined in Section 3 of \cite{Hendricks_2016}. The aforementioned identification is nicely summarized in the proof of \cite[~Theorem 2.4]{hendricks2022rank}. By Propositions $3.23$ and $3.24$ of \cite{Hendricks_2016}, this spectral sequence is independent of the choices in its construction and is invariant under equivariant Hamiltonian isotopies of the Lagrangians. Lemma \ref{eqhamiso} shows that the Heegaard diagrams $\mathcal{H}(D_n)/\tau$ and $\mathcal{H}(D_n')/\tau$ are related by a sequence of handleslides, isotopies of the attaching circles, and type 1-2 (de)stabilizations which avoid the fixed points of $\tau$. By \cite[~Lemma 5.1]{Hendricks_2016}, such handleslides of the $\overline{\alpha}$ or $\overline{\beta}$ curves induce Hamiltonian isotopies of $L_0'$ and $L_1'$ which avoid the divisors $x\times \Sym^{g-1}(\overline{\Sigma})$ for $x$ equal to $\overline{w}$ or the top or bottom fixed point. By \cite[~Section 7.3]{OS3manifolds1}, isotopies of the $\overline{\alpha}$ and $\overline{\beta}$ attaching circles avoiding $\overline{w}$ and the fixed points can be realized as deformations of the almost complex structures and Hamiltonian deformations of $L_0'$ and $L_1'$ avoiding the same divisors. These Hamiltonian isotopies in $\Sym^g(\overline{\Sigma}\backslash\{{\overline{w}\}})$ lift to equivariant Hamiltonian isotopies in $\Sym^{2g}(\Sigma/\{w,z\})$ since they avoid the branch points of $\Sym^{2g}(\Sigma)\to \Sym^g(\overline{\Sigma})$. Finally, we may assume that index 1-2 stabilizations occur near the basepoints and so, if we choose our complex structures to be split near the basepoints, these stabilizations have no effect at all on the equivariant Floer complex (cf. \cite[~Proof of Theorem 10.1]{OS3manifolds1}).

The key observation to reduce the $E_1$ page of this spectral sequence to $\HFKhat(K,0)\otimes \F[\theta,\theta^{-1}]$ is that any generator in Alexander grading not equal to $0$ must be killed by the $d_2$ differential. We explain this as follows. First let us recap how Large's Theorem \ref{Largelocal} works in our setup. We have the chain complex $CF(L_0,L_1)$ equipped with the (knot Floer) differential $\partial$, and we may also consider the map $\tau_\#\colon CF(L_0,L_1)\to CF(L_0,L_1)$ which takes $x\in L_0\cap L_1$ to $\tau(x)\in L_0\cap L_1$. The map $\tau_\#$ is not in general a chain map. Assuming that the hypotheses of Theorem \ref{Largelocal} are satisfied, then Large shows that there is an equivariant exact Lagrangian isotopy which replaces $(L_0,L_1)$ with new Lagrangians $(\widetilde{L}_0,\widetilde{L}_1)$ and fixes the invariant sets, and also a family of complex structures $\bf{J}$ so that the map $\widetilde{\tau}_\#$ on $CF(\widetilde{L}_0,\widetilde{L}_1)$ (defined by taking $x\in \widetilde{L}_0\cap\widetilde{L}_1$ to $\tau(x)\in \widetilde{L}_0\cap\widetilde{L}_1$) is a chain map. It should be noted that even in the case that $\tau_\#$ was a chain map, $\widetilde{\tau}_\#$ is not necessarily chain homotopy equivalent to $\tau_\#$. However, the manipulation of the Lagrangians described above to obtain a chain map $\widetilde{\tau}_\#$ crucially does not affect the Alexander grading. Proposition \ref{tauonknot} then tells us that $\widetilde{\tau}_\#$ negates Alexander gradings.
 The spectral sequence of Theorem \ref{Largelocal} can be computed from the Leray spectral sequence induced from the vertical filtration of the double complex $CF(\widetilde{L}_0,\widetilde{L}_1)\otimes \F[\theta,\theta^{-1}]$ equipped  differential $\widetilde{\partial}$ in the vertical direction and $(1+\widetilde{\tau}_\#)\theta$ in the horizontal direction. 
 \begin{center}
\begin{tikzcd}
&\vdots\arrow{d}&\vdots\arrow{d}&\\
\dots\arrow{r}&CF(\widetilde{L}_0,\widetilde{L}_1) \arrow{r}{\theta(1+\widetilde{\tau}_{\#})} \arrow{d}{\partial} & \theta\cdot CF(\widetilde{L}_0,\widetilde{L}_1) \arrow{r}{\theta(1+\widetilde{\tau}_\#)}\arrow{d}{\partial}&\dots \\
\dots\arrow{r}&CF(\widetilde{L}_0,\widetilde{L}_1) \arrow{r}{\theta(1+\widetilde{\tau}_{\#})} \arrow{d} & \theta\cdot CF(\widetilde{L}_0,\widetilde{L}_1) \arrow{r}{\theta(1+\widetilde{\tau}_{\#})}\arrow{d}&\dots\\
&\vdots&\vdots&
\end{tikzcd}.
\end{center}
In particular, the $d_2$ differential is given by $\theta(1+\widetilde{\tau}_*)$ acting on $HF(\widetilde{L}_0,\widetilde{L}_1)\otimes \F[\theta,\theta^{-1}]$. If $x$ is an element of $HF(\widetilde{L}_0,\widetilde{L}_1)$ with $A(x)\neq 0$, then since $A(\widetilde{\tau}_*(x))=-A(x)$, it must be the case that $x\neq \widetilde{\tau}_*(x)$ and hence if a sum of homology classes $x+\hdots $ is in the kernel of $1+\tau_*$, then $\widetilde{\tau}_*(x)$ must also appear with a coefficient of $1$ in that linear combination.$\qed$\\

\subsubsection{Proof of Theorem \ref{Mainthm2}}\label{pfmainthm2}
Now we work in the symmetric product
\begin{equation}
\check{M}\coloneq \Sym^{2g}(\Sigma\backslash\{w_1,w_2,z^*\})
\end{equation}
equipped with the submanifolds
\begin{equation}
L_0\coloneq \mathbb{T}_{\bm{\alpha}}=\alpha_1\times\hdots\times\alpha_{2g}\textrm{ and } L_1\coloneq \mathbb{T}_{\bm{\beta}}=\beta_1\times\hdots\times\beta_{2g}.
\end{equation}

Recall that $z=w_1$ and $w=w_2$, so that we can think of $\check{M}$ as a subset of $M$, where $M$ is defined in the proof of Theorem \ref{Mainthm1}. Pulling back the form $\omega$ constructed on $M=\Sym^{2g}(\Sigma\backslash\{w,z\})$ to $\check{M}$ tells us that $(1)$ of Theorem \ref{Largelocal} is satisfied for $(\check{M}, L_0, L_1, \omega, \tau).$ As mentioned before, the proof that $(2)$ of Theorem \ref{Largelocal} is also satisfied in this setup is the subject of Sections \ref{geosym} and \ref{stableTNiso}, and so we assume for this section that this assumption is satisfied.

In order to better understand the fixed point sets $(\check{M}^{\textrm{inv}}, \Lo, \Li)$, we define 
\begin{equation}\label{barM'}
\check{M}'=\Sym^g(\overline{\Sigma}\backslash \{\overline{z},\overline{w}\}),
\end{equation}
and define $L_0'$ and $L_1'$ as in Equation \ref{L0'L1'}.

We again have that the map $\phi\colon (\check{M}', L_0', L_1')\to (\check{M}^{\textrm{inv}}, \Lo, \Li)$ defined by 
\begin{equation}
    \phi((\ell(x_1)\hdots \ell(x_g)))=(x_1\tau(x_1)\hdots x_g\tau(x_g))
\end{equation}
is biholomorphic. And, again letting $\phi^*(\omega|_{\check{M}^{\textrm{inv}}})$ be the symplectic form on $\check{M}',$ the Lagrangian intersection Floer homologies $HF(\Lo, \Li)$ and $HF(L_0', L_1')$ are isomorphic. 

Applying Theorem \ref{Largelocal} along with the isomorphism $HF(\Lo, \Li)\cong HF(L_0', L_1')$ gives us a spectral sequence with $E_1$ page equal to $HF(L_0,L_1)\otimes \F[\theta,\theta^{-1}]=\HFKhat(S_b^n(K))\otimes \F[\theta,\theta^{-1}]$ and $E_\infty$ page isomorphic to $HF(L_0', L_1')\otimes \F[\theta,\theta^{-1}]=\HFKhat(K)\otimes \F[\theta,\theta^{-1}].$ 

This spectral sequence is invariant under choice of $\tau$-equivariant Heegaard diagram for $S^n_b(K)$, and in particular under the choice of the intravergently singularized diagram $IS(D_n)$ for $S^n_b(K)$; this is an identical argument to the invariance proof of the spectral sequence of Theorem \ref{Mainthm1}. 

We now argue that this spectral sequence splits along Alexander gradings. As explained in the proof of Theorem \ref{Mainthm1}, to compute the spectral sequence of Theorem \ref{Largelocal}, one first replaces $CF(L_0,L_1)$ with the chain homotopy equivalent complex $CF(\widetilde{L}_0,\widetilde{L}_1)$, and $\tau_\#$ with the map $\widetilde{\tau}_\#$ which has the same behavior on Alexander gradings as $\tau_\#$. In particular, Proposition \ref{tauonsing} states that $\tau_\#$ preserves Alexander gradings, and hence $\widetilde{\tau}_\#$ does as well. The spectral sequence of Theorem \ref{Largelocal} is then

\begin{enumerate}
    \item the vertically filtered Leray spectral sequence of the double complex \[(CF(\widetilde{L}_0,\widetilde{L}_1)\otimes \F[\theta],\widetilde{\partial},\theta(1+\widetilde{\tau}_\#))\] -- which has $E^\infty$ page equal to $CF_{\Z/2\Z}(\widetilde{L}_0,\widetilde{L}_1)\otimes \F[\theta]$ -- followed by
    \item identification of $CF_{\Z/2\Z}(\widetilde{L}_0,\widetilde{L}_1)\otimes \F[\theta,\theta^{-1}]$ with $CF(\Lo,\Li)\otimes \F[\theta,\theta^{-1}]$ via the \emph{localization isomorphism} of \cite[~Theorem 1.1]{Large}.
\end{enumerate}
We have just argued above that the vertical and horizontal differentials in the first step preserve Alexander gradings, and hence the spectral sequence of this double complex splits along Alexander gradings. The localization isomorphism of the second step also preserves Alexander gradings as it is constructed by counting pseudoholomorphic disks of varying index and multiplication and division by powers of $\theta$ and no flowlines pass over the basepoint divisors $V_{\overline{z}}=\overline{z}\times \Sym^{g-1}(\overline{\Sigma})$ or $V_{\overline{w}}=\overline{w}\times \Sym^{g-1}(\overline{\Sigma})$. So we have that the spectral sequence splits over Alexander gradings. For the argument that Alexander grading $2a-\widetilde{lk}(K)-\frac{1}{2}$ is carried to grading $a$, see the appendix.
\qed

The spectral sequence of Theorem \ref{Mainthm1} collapses on the $E_2$ page for every strongly invertible knot with less than $12$ crossings; see Remark \ref{thm1collapse<12}. However, as stated in Conjecture \ref{pagecollapse} we do not expect that this is the general behavior of this spectral sequence for all strongly invertible knots.
As explained in the proof above, the spectral sequence of Theorem \ref{Mainthm1} is the same as the spectral sequence \begin{equation}\label{ss3}\HFKhat(K)\otimes \F[\theta]\rightrightarrows CF_{\Z/2\Z}(L_0, L_0)\end{equation}induced from the equivariant Floer complex from non-invariant complex structures defined in Section 3 of \cite{Hendricks_2016}, followed by the localization isomorphism 
\begin{equation}\label{localiso}\theta^{-1}CF_{\Z/2\Z}(L_0, L_0)\cong \HFhat(S^3)\otimes \F[\theta,\theta^{-1}]\end{equation} of Large. With this in mind, we can make the following definition.

\begin{definition}
    The localization isomorphism Equation \ref{localiso} establishes that, thought of as an $\F_2[\theta]-$module, $CF_{\Z/2\Z}(L_0, L_0)$ has one $\theta$ tower. Let $s_\tau$ be the minimum grading of an element of the tower in $CF_{\Z/2\Z}(L_0, L_0)$. More precisely, set
    \[s_\tau(K)\coloneq \min\{d\;|\;\exists [x]\in CF_{\Z/2\Z}(L_0, L_0), M([x])=d, \theta^n[x]\neq 0\;\forall n\}.\]
\end{definition}

 \section{Examples}\label{examples}
 Before analyzing specific examples, we first record some general facts about the behavior of the spectral sequences of Theorem \ref{Mainthm1} and \ref{Mainthm2}. A knot $K$ is said to be Floer $\delta$-thin if the difference $\delta\coloneq M-A$ between the Maslov and Alexander gradings of any homogeneous element of $\widehat{HFK}(K)$ is a constant $\delta.$ We remind the reader that for a Floer $\delta$-thin knot, $\delta$ is equal to equivalently (suitably normalized versions of) the classical knot signature, Ozsvath-Szabo's concordance homomorphism $\tau$ \cite{OS4ball}, or Manolescu-Owens' concordance homomorphism $\delta$ \cite{manolescu2005concordanceinvariantfloerhomology}: \[\delta=-\frac{\sigma(K)}{2}=-\frac{\tau(K)}{2}=\delta(K).\] By a similar token, we will say that a singular knot $S$ is Floer $\delta$-thin if the same condition holds on $\widehat{HFK}(S)$. 
 \begin{proposition}\label{thinknotss1}
     If $(K,\tau)$ is a Floer $\delta$-thin strongly invertible knot, then $E_2=E_\infty$ for the spectral sequence of Theorem \ref{Mainthm1}. Furthermore, in this case the invariant $s_{\tau}$ is equal to $\delta$. 
\end{proposition}
\begin{proof}
    The assumptions that $K$ is Floer $\delta$-thin or L-space both imply that $\HFKhat(K)$ has exactly one occupied Maslov grading in Alexander grading $0$. Every differential $d_k$ for $k\ge 2$ does not preserve the Maslov grading because $\partial$, the vertical differential, lowers the Maslov grading by one, and $\theta\cdot(1+\widetilde{\tau}_{\#})$, the horizontal differential, preserves the Maslov grading by Proposition \ref{tauonknot}. Hence by the assumption, these higher differentials must be identically $0$. The unique occupied Maslov grading in Alexander grading $0$ is clearly $\delta$, and so the generator that survives the spectral sequence must have Maslov grading $\delta$.
\end{proof}

There is an analogous statement for the spectral sequence of Theorem \ref{Mainthm2} which follows by the same proof using Proposition \ref{tauonsing} in place of Proposition \ref{tauonknot}.
\begin{proposition}
    If $(K,\tau)$ is a strongly invertible knot such that $S^n_b(K)$ is Floer $\delta$-thin, then $E_2=E_\infty$ for the spectral sequence of Theorem \ref{Mainthm2}.
\end{proposition}
\begin{remark}\label{thm1collapse<12}Using Proposition \ref{thinknotss1}, computations of $\HFKhat$ from \cite{baldwin2007computationsheegaardfloerknothomology}, computations of the $CFK^\infty$ type of non-$\delta$-thin strongly invertible knots from \cite{antal2023note}, and properties of the map $\tau_{K}$ on $CFK^\infty(K)$ from \cite[~Theorem 1.7]{DMSE}, the author was able to conclude that the spectral sequence of Theorem \ref{Mainthm1} collapses on the $E_2$ page for all strongly invertible knots with at most $11$ crossings.
\end{remark}Despite this, we state the following conjecture:
\begin{Conjecture}\label{pagecollapse}
    There exists a strongly invertible knot $K$ such that the spectral sequence of Theorem \ref{Mainthm1} collapses on the $E_k$ page for some $k\ge 3.$
\end{Conjecture}

\subsection{The trefoil with unknotted quotient}\label{31plus}
    The spectral sequence of Theorem \ref{Mainthm1} collapses on the $E_1$ page as $dim(\HFKhat(3_1,0))=1=dim(\HFhat(S^3))$. The singular skein triple associated to $S^1_b(3_1^+)$ includes the positive Hopf link and the unknot, and completely determines the singular knot Floer homology of $S^1_b(3_1^+)$:\[\HFKhat(S^1_b(3_1^+))\cong \F_{(0,\frac{1}{2})}\oplus\F^2_{(-1,-\frac{1}{2})}.\]The linking number of $(S^1_b(3_1^+))_0=L_b^{\frac{3}{2}}(3_1^+)=0_1$ with the reversed symmetry axis is \[\lambda=lk(L_b^{\frac{3}{2}}(3_1^+),\widetilde{A})=1.\] Therefore, the proof in the appendix of the grading shift of Theorem \ref{Mainthm2} implies that Alexander grading $2a+\frac{2-1}{2}=2a+\frac{1}{2}$ of $\HFKhat(S^1_b(3_1^+))$ gets sent under the spectral sequence to Alexander grading $a+\frac{1-1}{2}=a$ of $\HFKhat(0_1)\cong \F_{(0,0)}$. The only occupied Alexander grading in the knot Floer homology of the quotient unknot is $a=0$ which corresponds to Alexander grading $2\cdot 0+\frac{1}{2}=\frac{1}{2}$ in the singular knot Floer homology of $S_b^1(3_1^+)$. This tells us that the two generators in grading $(-1,-\frac{1}{2})$ must be swapped by $\tau_*$ so that the $1+\tau_*$ differential on the $E_1$ page kills both of them, and the generator in grading $(0,\frac{1}{2})$ survives the spectral sequence and gets sent to the non-zero element of the $E_\infty=E_2$ page. This can all be seen explicitly from the Heegaard diagrams of Figure \ref{TrefoilHeegaardDiagram} as well.
\subsection{The trefoil with trefoil quotient}\label{31minus}
\begin{figure}
    \center
    \includegraphics[scale=.5]{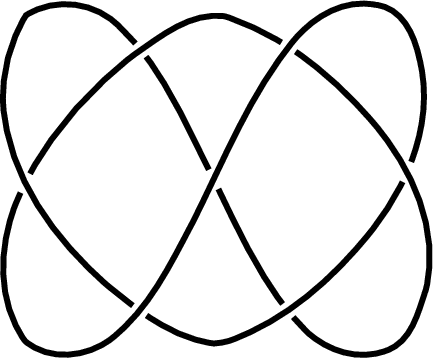}
    \caption{An intravergent diagram $D_3$ for $3_1^-$}
    \label{fig:31minus}
\end{figure}

If we consider the same strong inversion on the trefoil but with the opposite choice of half-axis, we obtain a DSI that we denote $3_1^-$; it is pictured in Figure \ref{fig:31minus}. Notice that $3_1^-$ is a right handed trefoil with quotient equal to a left handed trefoil, while $3_1^+$ was a left handed trefoil with unknotted quotient. We choose to use the right hand trefoil for $3_1^-$ so as to have a convenient intravergent diagram with negative central crossing.
The spectral sequence of Theorem \ref{Mainthm1} does not depend on the choice of half axis and will therefore be the same as described in the last example up to mirroring. The singular skein triple associated to $D_3$ consists of $S^3_b(3_1^-),$ $L^3_b(3_1^-)=T(2,6)$ and $L^{\frac{7}{2}}_b(3_1^-)=(3_1)_L\#(3_1)_L,$ where the subscript $L$'s denote that we are taking a connect sum of \emph{left hand} trefoils. The skein exact triangle associated to this triple along with the Alexander polynomial $\Delta_{S^3_b(3_1^-)}$ completely determine the singular knot Floer homology of $S^3_b(3_1^-):$
\[\HFKhat(S^3_b(3_1^-))\cong \F_{(5,\frac{5}{2})}\oplus \F_{(2,\frac{1}{2})}\oplus \F^2_{(1,-\frac{1}{2})}\oplus \F_{(0,-\frac{3}{2})}.\]This singular knot is not Floer $\delta$-thin, but nonetheless we will see that the spectral sequence of Theorem \ref{Mainthm2} dies on the $E^2$ page. From the diagram obtained by doing a resolution to $D_3$ it can be see that $\lambda=lk(L^{\frac{7}{2}}_b(3_1^-),\widetilde{A})=-1$ which means that
\begin{itemize}
    \item Alexander grading $2\cdot 2+\frac{-1-2}{2}=\frac{5}{2}$ of $\HFKhat(S^3_b(3_1^-))$ gets sent to Alexander grading $2+\frac{-1-1}{2}=1$ of $\HFKhat((3_1)_L)\cong \F_{(2,1)}\oplus \F_{(1,0)}\oplus \F_{(0,-1)}$
    \item Alexander grading $2\cdot 1+\frac{-1-2}{2}=\frac{3}{2}$ of $\HFKhat(S^3_b(3_1^-))$ gets sent to Alexander grading $1+\frac{-1-1}{2}=0$ of $\HFKhat((3_1)_L)$
    \item Alexander grading $2\cdot 0+\frac{-1-2}{2}=\frac{-1}{2}$ of $\HFKhat(S^3_b(3_1^-))$ gets sent to Alexander grading $0+\frac{-1-1}{2}=-1$ of $\HFKhat((3_1)_L)$
    \item Alexander grading $\frac{3}{2}$ of $\HFKhat(S^3_b(3_1^-))$ is not of the form $2a+\frac{\lambda-2}{2}$ and hence dies in the spectral sequence. 
\end{itemize}
These observations together tell us that three generators in Alexander grading not equal to $-\frac{1}{2}$ of of $\HFKhat(S^3_b(3_1^-))$ abut to the generators of the quotient left hand trefoil homology in the spectral sequence, while the two generators in Alexander grading $-\frac{1}{2}$ must be swapped by the $\tau_*$ action so that they cancel each other out in homology on the $E_2$ page.

\subsection{The figure 8 with cinquefoil quotient}\label{41minus}
\begin{figure}
    \center
    \includegraphics[width=.3\linewidth]{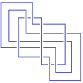}
    \caption{An intravergent diagram $D_n$ for $4_1^-$}
    \label{fig:41cinquefoil}
\end{figure}

The unique strong inversion on $4_1$ yields quotient knot equal to the unknot or the cinquefoil $5_1$ depending on which half axis is chosen. Here we analyze the choice of half axis that yields a $5_1$ quotient, denoted by $4_1^-$. This is displayed as a transvergent symmetry in Figure \ref{fig:41cinquefoil}. We have no need to determine the number $n$ associated to this diagram, and so we will leave it as an abstract label.

\begin{figure}
    \center
    \includegraphics[width=0.3\linewidth]{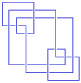}
    \caption{The two component two periodic butterfly link $S^n_b(4_1^-)_0$}
    \label{fig:41res}
\end{figure}

While we do not attempt to discern the precise behavior of the spectral sequence of Theorem \ref{Mainthm2} for this DSI, we can argue that $E_2\neq E_\infty$, thus demonstrating an answer in the affirmative to the analog of Conjecture \ref{pagecollapse} for Theorem \ref{Mainthm2}. The two other members of the negative skein triple associated to the diagram $D$ are the knot $12_{n725}$, found using Miller's Knotfol.io website \cite{knotfol}, and the link pictured in Figure \ref{fig:41res} which we shall just denote $S^n_b(4_1^-)_0$. We used the knot Floer homology calculations in \cite{baldwin2007computationsheegaardfloerknothomology} along the formula for the knot Floer homology of a mirror knot to find $\HFKhat(12_{n725})$, and a modification of the grid homology program from that paper to find $\HFKhat(S^n_b(4_1^-)_0)$.  These computations lead us to the following skein exact triangle:
\vspace{1cm}
\begin{center}
\begin{tikzpicture}[transform canvas={scale=.9}]

% Define the matrix
\matrix (m) [matrix of math nodes, row sep=2em, column sep=2em]
{
    \substack{\F_{(10,5)}\oplus \F_{(9,4)}^2\oplus \F_{(9,3)}\oplus \F_{(8,3)}\oplus\\  \F^4_{(8,2)}\oplus \F^6_{(7,1)}\oplus \F^6_{(6,0)}\oplus \F^6_{(5,-1)}\oplus\\ \F^4_{(4,-2)}\oplus \F_{(3,-3)}\oplus \F_{(2,-3)}\oplus \F^2_{(1,-4)}\oplus \F_{(0,-5)}} && \substack{\F_{(10,5)}\oplus\F^2_{(9,4)}\oplus\F_{(8,3)}\oplus\F_{(6,2)}\oplus\\ \F^3_{(5,1)}\oplus\F^5_{(4,0)}\oplus\F^5_{(3,-1)}\oplus\F^3_{(2,-2)}\oplus\\ \F_{(1,-4)}\oplus\F_{(1,-3)}\oplus\F^2_{(0,-5)}\oplus\F_{(-1,-6)}} \\
    & (\HFKhat(S^n_b(4_1^-))\otimes V)[[0,-\frac{1}{2}]] &\\
};
% Draw the arrows
\draw[->] (m-1-1) -- (m-1-3) node[midway, above] {$(0,0)$};
\draw[->] (m-2-2) -- (m-1-1) node[midway, below left] {$(0,0)$};
\draw[->] (m-1-3) -- (m-2-2) node[midway, below right] {$(-1,0)$};
\end{tikzpicture}
\end{center}
\vspace{1cm}
Observe that $\F_{(9,3)}\oplus \F^4_{(8,2)}$ necessarily maps to $0$ under the top horizontal map since that map has bidegree $(0,0)$ and there are not generators in gradings $(9,3)$ or $(8,2)$ on the top right. Therefore $\F_{(9,3)}\oplus \F^4_{(8,2)}$ is in the image of the right diagonal arrow of bidegree $(0,0)$ which implies that $\HFKhat(S^n_b(4_1^-))\otimes V)[[0,-\frac{1}{2}]]$ contains a $\F_{(9,3)}\oplus \F^4_{(8,2)}$ direct summand. Notice also that there is no generator in grading $(10,4)$ in the top left or in grading $(11,4)$ in the top right, and hence there is no $(10,4)$ graded generator in $\HFKhat(S^n_b(4_1^-))\otimes V)[[0,-\frac{1}{2}]]$. This implies that  there is a generator in grading $(9,3)$ in $\HFKhat(S^n_b(4_1^-))[[0,-\frac{1}{2}]]$. Similarly, we see that one of the four $\F_{(8,2)}$ summands in $\HFKhat(S^n_b(4_1^-))\otimes V)[[0,-\frac{1}{2}]]$ comes from the $\F_{(9,3)}\otimes \F_{(-1,-1)}$ tensor product occurring in $\HFKhat(S^n_b(4_1^-))\otimes V)[[0,-\frac{1}{2}]]$, while the other three must belong to $\HFKhat(S^n_b(4_1^-))[[0,-\frac{1}{2}]$. We have shown there is exactly one generator in grading $(9,\frac{7}{2})$ and three generators in grading $(8,\frac{5}{2})$ in $\HFKhat(S^n_b(4_1^-))$. Taking homology under the $1+\tau_*$ differential, which preserves both Alexander and Maslov gradings by Proposition \ref{tauonsing}, must leave an odd (and in particular positive) number of generators in gradings $(9,\frac{7}{2})$ and $(8,\frac{5}{2})$ on the $E_2$ page. By Theorem \ref{Mainthm2}, generators that survive the spectral sequence must have an even difference in Alexander gradings, and hence not both of these generators persist to the $E_\infty$ page. 

\section{Geometry of the symmetric product}\label{geosym}
Recalling the definitions of $M'$, $L_0'$ and $L_1'$ from the proof of Theorem \ref{Mainthm1}, define
\begin{equation}
    X'\coloneq (L_0'\times\{0\})\cup (L_1'\times\{1\})\subset M'\times[0,1].
\end{equation}
In this section we analyze the (co)homology of the spaces $M'\times[0,1],$ $X'$ and $M'\times[0,1]/X'$. The goal is to prove the following proposition.

\begin{proposition}\label{coM'/X'freeabelianfiniterank}
    The reduced cohomology $\widetilde{H}^*(M'\times[0,1]/X')$ is a free abelian group of finite rank.
\end{proposition}
Let $\nu_1,\hdots,\nu_{2g}:[0,1]\to \Sigma\backslash\{z\}$ be parameterizations of simple closed curves with the following properties. By a small abuse of notation we shall let $\nu_i$ denote the parameterization, its image, and its homology class where convenient.
\begin{enumerate}
    \item The $\nu_i$ all start at the same point, which we think of as their wedge point: \[\nu_1(0)=\hdots =\nu_{2g}(0).\]
    \item The $\nu_i$ are disjoint besides at the wedge point: for  $0<t_1,t_2<1$ and $i\neq j$ \[\nu_i(t_1)\neq\nu_j(t_2).\]
    \item The $\nu_i$ generate the first homology of $\Sigma\backslash\{z\}$: \[H_1(\Sigma\backslash\{z\})=\Z\langle \nu_1,\hdots,\nu_{2g}\rangle.\] 
    \item The punctured surface $\Sigma\backslash\{z\}$ deformation retracts onto the wedge product $\nu_1\vee\hdots \vee \nu_{2g}$.
\end{enumerate}
Taking symmetric products preserves homotopy equivalence, and so $M'=\Sym^g(\Sigma\backslash\{z\})$ deformation retracts onto $\Sym^g(\bigvee_{j=1}^{2g}\nu_j)$. The following lemma of Ong, proved in detail in \cite{Hendricks2012rank}, describes the homotopy type of this space. 
\begin{lemma}\label{symprodwedgecirc}\cite[~Lemma 5.1]{Hendricks2012rank} If $r\le k$, then the $r^{\textrm{th}}$ symmetric product of a wedge of $k$ circles $\Sym^r(\bigvee_{j=1}^kS_j^1)$ 
is homotopy equivalent to the $r$-skeleton of the $k$ torus 
$\prod_{j=1}^kS_j^1$, where each circle is given a CW structure consisting of the wedge point and a single one-cell, and the torus has the natural product CW structure.
\end{lemma}

Composing the homotopy equivalence of Lemma \ref{symprodwedgecirc} with the deformation retraction $M'\approx\Sym^g(\bigvee_{i=1}^{2g}\nu_i)$ gives a homotopy equivalence $G$ of $M'$ with the $g$-skeleton of $\prod_{i=1}^{2g} \nu_i$. Define
\[\bar{\nu}_i=G_*^{-1}([t\to (1,\hdots, \nu_i(t),\hdots,1)])\in H_1(M').\]
Then the first homology of $M'$ is 
\begin{equation}
H_1(M')=\Z\langle \bar{\nu}_1,\hdots ,\bar{\nu}_{2g}\rangle.
\end{equation}
Applying the Kunneth theorem and the fact that $M'$ is homotopy equivalent to the $g$-skeleton of $\prod_{i=1}^{2g} \nu_i$, for $m=1,\hdots ,g$ the homology group $H_m(M')$ is spanned by $m$-tensors of the $\bar{\nu}_i$ with no repeated indices; we will write this as \begin{equation}\label{hlgyM'}
    H_m(M')\cong \Lambda^m H_1(M)=\Lambda^m \Z\langle\bar{\nu}_1,\hdots ,\bar{\nu}_{2g}\rangle.
\end{equation}
The first homology of $M'$ has an alternative description \cite[~Lemma 2.6]{OS3manifolds1}
\begin{equation}\label{ozsziso}H_1(\Sigma\backslash\{z\})\cong H_1(\Sym^g(\Sigma\backslash\{z\})) =H_1(M').\end{equation}
The homology classes of the attaching curves \[\bm{\alpha}\cup\bm{\beta}=\{\alpha_1,\hdots ,\alpha_g,\beta_1,\hdots ,\beta_g\}\] form a basis for $H_1(\Sigma\backslash\{z\})$ since they specify a Heegaard decomposition for $S^3$. Let $\bar{\alpha_i}$ and $\bar{\beta}_i$ denote the images in $H_1(M')$ of $[\alpha_i]$ and $[\beta_i]$ under the isomorphism of Equation \ref{ozsziso}. Explicitly, letting $\alpha_i(t)$ be a parameterization of $\alpha_i$, then $\bar{\alpha}_i=[t\to (\alpha_i(t)x\hdots x)]$ where $x$ is any point in $\Sigma\backslash\{z\}$. We recast Equation (\ref{hlgyM'}) as 
\begin{equation}
H_m(M'\times[0,1])\cong H_m(M')\cong\Lambda^m \Z\langle \bar{\alpha}_1,\hdots ,\bar{\alpha}_g, \bar{\beta}_1,\hdots ,\bar{\beta}_g\rangle\textrm{ for } m=1,\hdots, g,
\end{equation}
and
\begin{equation}H_m(M'\times [0,1])=0\textrm{ for } m>g.\end{equation}
Next we consider the homology of the tori $L_0'=\mathbb{T}_{\bm{\alpha}}=\alpha_1\times\hdots\times\alpha_g$ and $L_1'=\mathbb{T}_{\bm{\beta}}=\beta_1\times\hdots\times\beta_g$. Choosing basepoints $x_i$ on each $\alpha_i$ and $y_i$ on each $\beta_i$ for $i=1,\hdots ,g$, and defining
\begin{align}
&\hat{\alpha}_i=[t\to (x_1\hdots x_{i-1}\alpha_i(t)x_{i+1}\hdots x_g)]\in H_1(L_0')\textrm{ and }\\
&\hat{\beta}_i=[t\to (y_1\hdots y_{i-1}\beta_i(t)y_{i+1}\hdots y_g)]\in H_1(L_1'),
\end{align}
the homology groups of $L_0'$ and $L_1'$ for $m=1,\hdots,g$ are 
\begin{align}
&H_m(L_0')=\Lambda^m\Z\langle \hat{\alpha}_1,\hdots, \hat{\alpha}_g\rangle,\\
&H_m(L_1')=\Lambda^m\Z\langle \hat{\beta}_1,\hdots, \hat{\beta}_g\rangle,
\end{align}
and are zero for $m>g.$ Abusing notation, we then write the homology of $X'=(L_0'\times\{0\})\cup(L_1'\times\{1\})$  as 
\begin{align}
&H_m(X')=\Lambda^m\Z\langle \hat{\alpha}_1,\hdots, \hat{\alpha}_g\rangle\oplus \Lambda^m\Z\langle \hat{\beta}_1,\hdots, \hat{\beta}_g\rangle\textrm{ for }m=1,\hdots,g\textrm{ and }\\
&H_m(X')=0\textrm{ for }m>g.
\end{align}
Letting $I\colon  X'\hookrightarrow M'\times[0,1]$ denote the inclusion map (lowercase $i$ is reserved for a different inclusion map in Section \ref{stableTNiso}), then the homology pushforward \[I_*\colon H_m(X')\to H_m(M'\times[0,1])\] admits the following description:
\begin{equation}
    I_*(\Lambda^m_{t=1}\hat{\alpha}_{j_t}\oplus \Lambda^m_{t=1}\hat{\beta}_{j'_t})= \Lambda^m_{t=1}\bar{\alpha}_{j_t}+\Lambda^m_{t=1}\bar{\beta}_{j'_t}
\end{equation}
where $(j_1,\hdots ,j_k)$ and $(j_1',\hdots ,j_k')$ are any collection of distinct integers between $1$ and $g$ inclusive. Notice that $\Lambda^m_{t=1}\bar{\alpha}_{j_t}= \Lambda^m_{t=1}\bar{\alpha}_{p_t}$ if and only if $(j_1,\hdots ,j_t)=(p_1,\hdots ,p_t)$, and similarly for wedges of the $\hat{\beta}_i$. Since $\Lambda^m_{t=1}\hat{\alpha}_{j_t}\oplus \Lambda^m_{t=1}\hat{\beta}_{j'_t}$ is a generic basis element of $H_m(X')$, and $\Lambda^m_{t=1}\bar{\alpha}_{j_t}$, and $\Lambda^m_{t=1}\bar{\beta}_{j'_t}$ are basis elements of $H_1(M'\times [0,1])$ we have showed the following.

\begin{lemma}The homology push-forward map $I_*\colon H_m(X')\to H_m(M'\times[0,1])$
is an injection for $m\ge1$.\end{lemma} 

We have also showed that $H_*(X')$ and $H_*(M'\times[0,1])$ are free abelian of finite rank because $X'$ is a disjoint union of tori and $M'\times[0,1]$ has the homotopy type of a skeleton of a torus. Every homology and cohomology group in sight will have finite rank from these observations combined with universal coefficients and the (co)homology long exact sequence, so from now on we do not mention rank. From the universal coefficient theorem the cohomology groups $H^*(X')$ and $H^*(M'\times [0,1])$ are free abelian as well. These observations tell us that in the commutative square

\[
\begin{tikzcd}
H^m(M'\times[0,1]) \arrow{r}{} \arrow[swap]{d}{I^*} & \textrm{Hom}(H_m(M'\times[0,1]),\Z) \arrow{d}{(I_*)^T} \\
H^m(X')  \arrow{r}{} & \textrm{Hom}(H_m(X'), \Z)
\end{tikzcd}
\]
the horizontal arrows are isomorphisms. The fact that the homology groups are free abelian further tells us that injectivity of $p_*$ implies surjectivity of $(I_*)^T$. Surjectivity of $(I_*)^T$ along with the fact that the above commutative square has isomorphisms for horizontal arrows yields the following lemma.

\begin{lemma}\label{M'Xsurj}
    The cohomology pullback map induced by inclusion \[I^*\colon H^m(M'\times[0,1])\to H^m(X')\] is a surjection for $m\ge 1$.
\end{lemma}
Lemma \ref{M'Xsurj} is of course equivalent to the same statement for reduced cohomology: the map $I^*\colon \widetilde{H}^m(M'\times[0,1])\to \widetilde{H}^m(X')$ for $m\ge 1$ is a surjection. 

We are now ready to prove Proposition \ref{coM'/X'freeabelianfiniterank}.

\begin{proof}[Proof of Proposition \ref{coM'/X'freeabelianfiniterank}.]
In the following we use $Q$ instead of $q$ because the lowercase letter already denotes the quotient map $S^3\to S^3/\tau$. Surjectivity of $I^*$ for $m\ge 1$ tells us that the cohomology long exact sequence of the pair $(M'\times[0,1],X')$
\[\hdots \to\widetilde{H}^m(M'\times[0,1]/X')\xrightarrow{Q^*}\widetilde{H}^m(M'\times[0,1])\xrightarrow{I^*}\widetilde{H}^m(X')\xrightarrow{\partial}\widetilde{H}^{m+1}((M'\times[0,1]/X')\to\hdots\]
breaks up into short exact sequences
\[
0\to\widetilde{H}^m(M'\times[0,1]/X')\xhookrightarrow{}\widetilde{H}^m(M'\times[0,1])\twoheadrightarrow\widetilde{H}^m(X')\to 0
\]
for $m\ge 2$. So, for $m\ge 2$, $\widetilde{H}^m(M'\times[0,1]/X')$ injects into the free abelian group $\widetilde{H}^m(M'\times[0,1])$, and hence it is free abelian itself. All that is left is to analyze $\widetilde{H}^1(M'\times[0,1]/X')$, but the first cohomology of a space with finitely generated homology is always free.
\end{proof}

\section{Stable tangent normal isomorphism} \label{stableTNiso}
 If $\pi\colon \M\times[0,1]\rightarrow\M$ is the projection map, define \begin{equation}E_N\coloneq \pi^*(N\M)=N(\M)\times[0,1]\end{equation}\begin{equation} E_T\coloneq \pi^*(T\M)=T(\M)\times[0,1]\end{equation}\begin{equation}X\coloneq (\Lo\times\{0\})\cup(\Li\times\{1\})\subset \M\times[0,1]\end{equation} and let $J$ be an almost complex structure tamed by $\omega$; in Heegaard Floer theory our usual choice is that $J$ is of the form $Sym^g(j)$ for $j$ a generic choice of almost complex structure on the Heegaard surface $\Sigma$, which is tamed by Perutz's symplectic form \cite{perutz}. We prove the following proposition which is an analog of a result alluded to in \cite[~Section 3d]{SS} and proven in detail in  \cite[~Proposition 7.1]{Hendricks2012rank}.
\begin{proposition}\label{symptocomp}
    Let \begin{equation*}
        f_N,f_T\colon (\M\times [0,1],X)\to (BU,BO)
    \end{equation*}
be the classifying maps of the (complex bundle $\xi$ over $\M\times [0,1]$, real subbundle of $\xi|_X$) pairs\[(E_N,(N(\Lo)\times\{0\})\cup (JN(\Li)\times\{1\}))\twoheadrightarrow(\M\times[0,1],X)\] and \[(E_T,(T(\Lo)\times\{0\})\cup (JT(\Li)\times\{1\}))\twoheadrightarrow(\M\times[0,1],X)\] 
respectively. Then a homotopy between $f_N$ and $f_T$ implies the existence of a stable tangent normal isomorphism. That is, to demonstrate a stable tangent normal isomorphism it suffices to produce a stable complex isomorphism between $E_N$ and $E_T$ that restricts to a stable real isomorphism of $(N(\Lo)\times\{0\})\cup (JN(\Li)\times\{1\})$ and $(T(\Lo)\times\{0\})\cup (JT(\Li)\times\{1\})$ over $X.$
\end{proposition}

The following proof is identical to the one given in Proposition 7.1 of \cite{Hendricks2012rank} except modifications to talk about a stable tangent normal isomorphism instead of a stable normal trivialization. This proposition should be thought of mainly as the statement that we can argue purely in terms of complex bundle theory instead of symplectic bundle theory when proving the existence of a stable tangent normal isomorphism. Indeed, right after this proposition, we will argue further that Proposition \ref{symptocomp}, and hence the stable tangent normal isomorphism, can be established by studying the reduced relative complex K theory classes of $E_N$ and $E_T$. 
\begin{proof}
Recall that any symplectic vector bundle
can be equipped with a compatible complex structure which is unique up to
homotopy, and any complex vector bundle similarly admits a
symplectification. Moreover, two symplectic vector bundles are isomorphic
if and only if their underlying complex vector bundles are isomorphic, and
isomorphisms of symplectic vector bundles map Lagrangian subbundles to
Lagrangian subbundles (see McDuff and
Salamon~\cite[Theorem~2.62]{McDuffSal}). Let $\omega_N$ and $\omega_T$ be the natural
symplectic structures on $N(\M)$ and $T(\M)$ coming from the symplectic
structure on $TM$.

Let the dimension of $N(\M)$ and $T(\M)$ be $k$. Suppose that $g_N \co N(\M) \rightarrow BU$ and $g_T\co T(\M) \rightarrow BU$ are classifying
maps of $N(\M)$ and $T(\M)$ thought of as a unitary vector bundles, so
that the images of $g_T$ and $g_N$ lie inside $BU_{k}$. Let
$$\zeta_{k}\co EU_{k} \rightarrow
BU_{k}$$
be the complex $k$--dimensional
universal bundle, and similarly let
$$\eta_{k}\co
EO_{k} \rightarrow BO_{k}$$
be the real
$k$--dimensional universal bundle. Equip
$EU_{k}$ with a symplectic structure $\omega_{\zeta}$ such
that $\eta_{k} \subset \zeta_{k}$ is a
Lagrangian subbundle. Then the bundles $(N(\M), \omega_N)$
and $(N(\M), g_N^*(\omega_{\zeta}))$ are isomorphic (indeed,
equal) as complex vector bundles, so there is a symplectic vector bundle
isomorphism
$$\chi_N\co(N(\M), \omega_M) \rightarrow
(N(\M), g^*(\omega_{\zeta})).$$
Similarly, the bundles $(T(\M), \omega_T)$
and $(T(\M), g_T^*(\omega_{\zeta}))$ are equal as complex vector bundles, so there is a symplectic vector bundle
isomorphism
$$\chi_T\co(T(\M), \omega_T) \rightarrow
(T(\M), g_T^*(\omega_{\zeta})).$$

These extend to symplectic
vector bundle isomorphisms
\begin{align*}&\widetilde{\chi}_N\co (E_N,
\pi^*(\omega_N)) \rightarrow (E_N,
(g_N\circ \pi)^*(\omega_{\zeta})= f_N^*(\omega_\zeta) )\textrm{ and }\\&\widetilde{\chi}_T\co (E_T, \pi^*(\omega_T)) \rightarrow (E_T,(g_T\circ \pi)^*(\omega_{\zeta})= f_T^*(\omega_\zeta))\end{align*}
where in all cases the
symplectic forms are constant with respect to the interval $[0,1]$, as are the maps
$\widetilde{\chi}_N$ and $\widetilde{\chi}_T$. From now on, we assume that we have first applied 
isomorphisms of this form to $E_N$ and $E_T$ so that the maps
$g_N \co E_N \rightarrow BU$ and $g_T\co E_T\rightarrow BU$ are in fact symplectic
classifying maps. We can if necessary pre- and post-compose the resulting stable
tangent normal isomorphism with $\widetilde{\chi}_N$ and $\widetilde{\chi}_T$.

Consider a homotopy $H$ of $f_N$ and $f_T$.
\begin{align*}
H \co (\M \times [0,1], \smash{\Lo} \times \{0\} \cup \smash{\Li} \times \{1\}) \times [0,1] &\rightarrow (BU,BO) \\
(x,t,s) &\mapsto h_s(x,t)
\end{align*}
Here the map $h_0$ is equal to $f_N$ and the map $h_1$ is equal to $f_T$.

Since $\M$ is homotopy equivalent to a compact subspace of itself, we may assume there is some $K>0$ such that if $s=k +K$, the image of $H$ lies inside $(BU_s, BO_s)$. Let $\zeta_s \co EU_s \rightarrow BU_s$ be the complex $s$--dimensional universal bundle with subbundle $\eta_s \co EO_s \rightarrow BO_s$ the real $s$--dimensional universal bundle. Then the pullbacks of $\zeta_s$ and $\eta_s$ along $h_1$ and $h_0$ are as follows.
\begin{align*}
h_0^*(\zeta_s) &= (E_N \oplus \mathbb C^K, h_0^*\omega_{\zeta}) \\
h_1^*(\zeta_s) &= (E_T\oplus \C^K, h_1^*\omega_{\zeta}) \\
(h_0|_{\smash{\Lo}\times \{0\}})^*(\eta_s) &= (N(\smash{\Lo}) \times \{0\}) \oplus \mathbb R^K \\
(h_0|_{\smash{\Li} \times \{1\}})^*(\eta_s) &= (J(N(\smash{\Li})) \times \{1\}) \oplus \mathbb R^K \\
(h_1|_{\smash{\Lo}\times \{0\}})^*(\eta_s) &= (T(\smash{\Lo}) \times \{0\}) \oplus \mathbb R^K \\
(h_1|_{\smash{\Li} \times \{1\}})^*(\eta_s) &= (J(T(\smash{\Li})) \times \{1\}) \oplus \mathbb R^K \\
\end{align*}
For clarity's sake it should be borne in mind that $\mathbb R^K$ and $\mathbb R^s$ refer to the canonical real subspaces in $\mathbb C^K$ and $\mathbb C^s$.

Since $H$ is a homotopy, it induces a stable isomorphism $\psi$ of $E_N$ with $E_T$. Write an arbitrary vector in $E_N$ as $(x,t,v)$ where $(x,t) \in \M \times [0,1]$ and $v$ is an element of the fiber over $(x,t)$.
\begin{align*}
\psi \co E_N \oplus \mathbb C^K = h_0^*(\zeta_s)
&\stackrel{\hbox{\footnotesize$\sim$}}{\longrightarrow} h_1^*(\zeta_s) = E_T \oplus \C^K\\
(x,t,v) & \mapsto \psi(x,t,v)
\end{align*}
The restrictions of $\psi$ to $(N(\smash{\Lo}) \times \{0\}) \oplus \mathbb  R^K$ and to $(J(N(\smash{\Li}) \times\{1\})) \oplus \mathbb R^K$ are real stable isomorphisms of these two bundles with $(T(\smash{\Lo}) \times \{0\}) \oplus \mathbb  R^K$ and  $(J(T(\smash{\Li}) \times\{1\})) \oplus \mathbb R^K$ respectively.
\begin{align*}
&\psi|_{(N(\smash{\Lo} \times \{0\}) \oplus \mathbb R^K)} \co (N(\smash{\Lo}) \times \{0\}) \oplus \mathbb R^K \rightarrow (T(\smash{\Lo}) \times \{0\}) \oplus \mathbb R^K \\
&\psi|_{(J(N(\smash{\Li}) \times \{1\}) \oplus \mathbb R^K)} \co (J(N(\smash{\Lo})) \times \{1\}) \oplus \mathbb R^K \rightarrow (J(N(\smash{\Lo})) \times \{1\}) \oplus \mathbb R^K
\end{align*}
Since $E_N$ is the pullback of $N(\M)$ to $\M \times [0,1]$, the map $h_0 = f_N$ is constant with respect to the interval $[0,1]$. That is, $h_0(x,t_1) = h_0(x,t_2)$ for all $x \in \M$ and $t_1,t_2 \in [0,1]$. Then each $\psi_t = \psi|_{\M \times \{t\}}$ is a stable isomorphism of $N(\M) \times \{t\} = N(\M)$ with $T(\M)$.  More concretely, we have symplectic isomorphisms
\begin{align*}
\psi_t \co N(\M) &\rightarrow T(\M)\\
(x,v) &\mapsto \psi(x,t,v).
\end{align*}
We will use the family of isomorphisms $\psi_t$ to produce the stable tangent normal isomorphism. Consider a map $\alpha$ defined by applying $\psi_0$ to each $\M \times \{t\} \subset \M \times [0,1]$.
\begin{align*}
\alpha \co E_N &\rightarrow E_T\\
(x,t,v) &\mapsto \psi_0((x,v)) = \psi(x,0,v).
\end{align*}
This is a stable isomorphism of $E_N$ with $E_T$. Because the symplectic structure on $E_N$ is constant with respect to the interval $[0,1]$, it is in fact a symplectic isomorphism of vector bundles. Consider the following Lagrangian subbundles $\Lambda_0$ and $\Lambda_1$ of $E_N$; we make these definitions with the intention to evaluate $\alpha$ on $\Lambda_0$ to produce a homotopy of Lagrangian subbundles $\alpha((N(\Lo)\times\{0\})\oplus \R^K)$ with $(T(\Lo)\times\{1\})\oplus \R^K$, and to evaluate $\alpha$ on $\Lambda_1$ to produce a homotopy of Lagrangian subbundles $\alpha((N(\smash{\Li}) \times \{0\}) \oplus i\mathbb R^K)$ with $(T(\smash{\Li}) \times \{1\}) \oplus i\mathbb R^K$, as required by Definition \ref{isopol}.
\begin{align*}
\Lambda_0|_{\smash{\Lo} \times \{t\}} &= (N(\smash{\Lo}) \times \{t\}) \oplus \mathbb R^K \\
\Lambda_1|_{\smash{\Li} \times \{t\}} &= \psi_0^{-1} \circ \psi_t(N(\smash{\Li}) \times \{t\} \oplus i \mathbb R^K).
\end{align*}
Since the maps $\psi_t$ form a homotopy, $\Lambda_1$ is a smooth subbundle as desired. Both subbundles are Lagrangian since their restriction to each $L_i^{\textrm{inv}} \times \{t\}$ is Lagrangian.    The next thing to check is the restriction of $\Lambda_i$ to $L^{\textrm{inv}}_i \times \{0\}$ for $i=0,1$.
\begin{align*}
\Lambda_0|_{\smash{\Lo} \times \{0\}} &= (N(\smash{\Lo}) \times \{0\}) \oplus \mathbb R^K \\
\Lambda_1|_{\smash{\Li} \times \{0\}} &= \psi_0^{-1} \circ \psi_0((N(\smash{\Li}) \times \{0\}) \oplus i\mathbb R^K)
= (N(\smash{\Li}) \times \{0\}) \oplus i \mathbb R^K.
\end{align*}
The other condition to check is that $\phi(\Lambda_0|_{L^{\textrm{inv}}_0 \times \{1\}})$ is $(T(\Lo)\times\{0\})\oplus \R^K$ and $\phi(\Lambda_1|_{L^{\textrm{inv}}_1 \times \{1\}})$ is $(T(\smash{\Li}) \times \{1\}) \oplus i\mathbb R^K$.
\begin{align*}
\phi(\Lambda_0|_{\smash{\Lo} \times \{1\}}) &= \psi_0((N(\smash{\Lo}) \times \{0\}) \oplus \mathbb R^K) \\
&= \psi((N(\smash{\Lo}) \times \{0\}) \oplus \mathbb R^K) \\
&=  (T(\Lo)\times\{0\})\oplus \R^K\\
\phi(\Lambda_1|_{\smash{\Li} \times \{1\}}) &= \psi_0( \psi_0^{-1} \circ \psi_1((N(\smash{\Li}) \times \{0\}) \oplus i\mathbb R^K) \\
&= \psi_1(J(J((N(\smash{\Li}) \times \{0\}) \oplus i \mathbb R^K)))	\\
&=J(\psi_1(J(N(\smash{\Li}	) \times \{0\}) \oplus \mathbb R^k))	\\								
&= J((J(T(\smash{\Lo})) \times \{1\}) \oplus \mathbb R^K)\\
&= (T(\smash{\Li}) \times \{1\}) \oplus i\mathbb R^K
\end{align*}
Therefore, $\alpha(\Lambda_0|_{\smash{\Lo} \times \{t\}})$ provides a homotopy through Lagrangians of the Lagrangian subbundles $\alpha((N(\Lo)\times\{0\})\oplus \R^K)$ and $(T(\Lo)\times\{1\})\oplus \R^K$, and $\alpha(\Lambda_1|_{\smash{\Li} \times \{t\}})$ provides a homotopy through Lagrangians of the Lagrangian subbundles $\alpha((N(\Li)\times\{0\})\oplus i\R^K)$ and $(T(\Li)\times\{1\})\oplus i\R^K$.\end{proof}

We show that the stable tangent normal isomorphism hypothesis of Theorem \ref{Largelocal} is met for $(M,L_0,L_1,\tau)$ defined as in Section \ref{equivariant} by using Proposition \ref{symptocomp}; it will suffice to find a stable complex isomorphism of $E_N$ with $E_T$ that restricts to a stable real isomorphism of \[(N(\Lo)\times\{0\})\cup J(N(\Li)\times\{1\})\] with \[(T(\Lo)\times\{0\})\cup J(T(\Li)\times\{1\}).\] These are trivializable totally real subbundles of $E_N|_X$ and $E_T|_X$ respectively, as they are tangent bundles of tori. After picking real trivializations, we may tensor with $\C$ to obtain complex trivializations \[\phi_N\colon E_N|_X\to \C^k\] and \[\phi_T\colon E_T|_X\to \C^k.\]The composition \[f\coloneq \phi_T^{-1}\circ \phi_N\colon E_N|_X\to E_T|_X\] is then an isomorphism of complex bundles over $X.$ We will show that $f$ can be stably extended to all of $E_N$ and $E_T$, which in particular implies that the hypotheses Proposition \ref{symptocomp} is satisfied. The strategy for showing that such an extension exists will be to show equality of relative reduced $K$-theory classes \begin{equation}\label{relKtheoryequality}[E_N]_{rel}=[E_T]_{rel}\in \widetilde{K}^0(\M\times[0,1],X).\end{equation}

\subsection{K-theory}
In this section we recall some results about reduced complex K-theory.  Let $\textrm{Vect}^{\mathbb{C}}(B)$ be the set of all complex finite dimensional vector bundles over a space $B$. For our purposes it suffices to take $B$ to be a compact manifold, although the following remarks and theorems apply to vector bundles over compact Hausdorff topological spaces. We refer the reader to \cite{HatcherVB} for a complete treatment of this subject.

    A stable isomorphism of complex vector bundles $E,E'\in \textrm{Vect}^{\mathbb{C}}(B)$ is a complex vector bundle isomorphism of $E\oplus \underline{\C}^m\cong E'\oplus \underline{\C}^n$ for some $m,n\in \Z$. Stable isomorphism is an equivalence relation on $\textrm{Vect}^{\mathbb{C}}(B)$ which we denote by $\sim$. Consider the set of equivalence classes $\textrm{Vect}^{\mathbb{C}}(B)/\sim$. Define addition and multiplication on these classes by $[E]+[E']=[E\oplus E']$ and $[E]\times [E']=[E\otimes E']$. Then
    \begin{equation}
        \widetilde{K}^0(B)=(\textrm{Vect}^{\mathbb{C}}(B)/\sim,+,\times).
    \end{equation}
    is a ring called the reduced zeroeth complex K-theory of $B$. Reduced complex K-theory can be made into a reduced cohomology theory by defining the higher $K$ theory groups \[\widetilde{K}^i(B)\coloneq \widetilde{K}^0(\Sigma^i(B))\] where $\Sigma^i(B)$ is the reduced suspension functor applied $i$  times to $B$ (to see that this definition actually yields a reduced cohomology theory \emph{Bott periodicity}, $\widetilde{K}^0(B)\cong \widetilde{K}^0(\Sigma^2(B))$, is employed), and the relative K-theory groups of a pair $X\subset B$ \[\widetilde{K}^i(B,X)\coloneq \widetilde{K}^i(B/X)\] where $X$ is a closed subspace of $B$. The Brown representability theorem then implies that each reduced complex K-theory group can be represented by homotopy classes of maps from $B$ into a classifying space; in particular if $BU$ is the classifying space of the infinite unitary group then
\begin{proposition}There is an isomorphism \begin{equation}\widetilde{K}^0(B)\cong[B,BU].\end{equation}This image of the (equivalence class of a) vector bundle under this isomorphism is called the \emph{classifying map} of the vector bundle.\end{proposition}
The following two propositions relate K-theory to ordinary (singular) cohomology.
\begin{proposition}\cite[~Theorem 4.5]{HatcherVB}\label{cherncharacter}
    There is a homomorphism
    \begin{equation}
    \widetilde{ch}\colon \widetilde{K}^0(B) \to \widetilde{H}^{\textrm{even}}(B,\Q)
    \end{equation}
called the \emph{reduced Chern character} between reduced complex K-theory and reduced singular cohomology with rational coefficients that enjoys the 
following properties:
\begin{itemize}
\item tensoring with $\Q$ makes $\widetilde{ch}$ into an isomorphism (which we will also call $\widetilde{ch}$) \[\widetilde{ch}\colon \widetilde{K}^0(B)\otimes_{\Z} \Q\xrightarrow{\sim} \widetilde{H}^{\textrm{even}}(B,\Q),\] and
\item the reduced Chern character $\widetilde{ch}([E])$ is a polynomial in the Chern classes of the (stable equivalence class of) the vector bundle $E$.
\end{itemize}
\end{proposition}
\begin{proposition}\cite[~Section 2.5]{MR139181}\label{HomKtorsionfree}
    If $\widetilde{H}^*(B,\Z)$ is torsion free and of finite rank, then so is $\widetilde{K}^0(B)$.
\end{proposition}
Propositions \ref{cherncharacter} and \ref{HomKtorsionfree} along with the fact that the $\widetilde{H}^*(\M\times[0,1],X)$ is torsion free and finite rank (c.f Proposition \ref{coM'/X'freeabelianfiniterank}) tell us that demonstrating equality of relative Chern classes $c((E_T)_{rel})=c((E_N)_{rel})$ would in particular imply Equation (\ref{relKtheoryequality}), and hence demonstrate the existence of a stable tangent normal isomorphism. This proof is the subject of the next subsection.
 
  \subsection{Algebraic topology of symmetric products}\label{chern}Here we show
\begin{proposition}  
\label{chernTNrel}
    The relative Chern classes $c((E_T)_{rel})$ and $c((E_N)_{rel})\in \widetilde{H}^{*}(\M\times [0,1],X)$ coincide:
    \begin{equation}\label{relequal}
c((E_T)_{rel})=c((E_N)_{rel})\in \widetilde{H}^{*}(\M\times[0,1],X).
\end{equation}
\end{proposition}
The lemmas and corollaries necessary to prove Proposition \ref{chernTNrel} are collected below.

\begin{lemma}\label{chernsquare}
The pullback of the total Chern class of $M$ along $i\colon \M \xhookrightarrow{} M$ is the square of the total Chern class of $M^{inv}:$
\begin{equation}
c(T\M)^2=i^*(c(TM))=c(i^*(TM))\in \widetilde{H}^{*}(\M).
\end{equation}
\end{lemma}
We prove Lemma \ref{chernsquare} later in this section. Assuming it, note the following corollary.
\begin{corollary}\label{chernTN}
The Chern classes of the tangent and normal bundles of $\M\subset M$ coincide:
\begin{equation}
c(T\M)=c(N\M)\in \widetilde{H}^{*}(\M).
\end{equation}
\end{corollary}
\begin{proof}[Proof of Corollary \ref{chernTN}]
The vector bundle isomorphism $T\M\oplus N\M\cong i^*(TM)$ implies \[
c(T\M)c(N\M)=c(i^*(TM))\in \widetilde{H}^{*}(\M).
\]
From Lemma \ref{chernsquare} we then get
\[
c(T\M)c(N\M)=c(T\M)^2\implies c(T\M)=c(N\M).
\]For this last implication, if $c(T\M)=1+a_1+\hdots +a_g$ and $c(N\M)=1+b_1+\hdots +b_g$ then \begin{equation}\label{cherncomp}(1+a_1+\hdots +a_g)(1+b_1+\hdots +b_g)=(1+\hdots +a_g)^2.\end{equation}We inductively see that $a_i=b_i$ for all $i$; the case $i=1$ is evident from analyzing the degree $1$ piece of the above equation, and if $a_i=b_i$ for $i< k\le g$, then from analyzing the degree $k$ piece of Equation \ref{cherncomp} we see that \[a_k+b_k+(a_{k-1}b_1+\hdots +b_{k_1}a_1)=a_k+a_k+(a_{k-1}a_1+\hdots +a_1a_{k-1})\] which reduces to $a_k=b_k$ after using $a_i=b_i$ for $i<k$.
\end{proof}

\begin{lemma}\label{injection}
The even degree cohomology map $q^*\colon \widetilde{H}^{\textrm{even}}(\M\times [0,1],X)\to \widetilde{H}^{\textrm{even}}(\M\times[0,1])$ is an injection.
\end{lemma}
\begin{proof}
Looking at the cohomology long exact sequence for the pair $(\M\times[0,1],X)$
\[\hdots\to \widetilde{H}^{2r-1}(\M\times[0,1])\xrightarrow{\iota^*}\widetilde{H}^{2r-1}(X)\xrightarrow{\delta^*}\widetilde{H}^{2r}(\M\times[0,1]/X)\xrightarrow{q^*}\widetilde{H}^{2r}(\M\times[0,1])\to\hdots\]
we see that \[q^*\colon \widetilde{H}^{2r}(\M\times[0,1]/X)\to \widetilde{H}^{2r}(\M\times[0,1])\] is an injection iff \[\iota^*\colon \widetilde{H}^{2r-1}(\M\times[0,1])\to \widetilde{H}^{2r-1}(X)\] is a surjection. Recall the biholomorphism $\phi\colon (M', L_0', L_1')\to (\M, \Lo, \Li)$, where $M'=\Sym^g(\overline{\Sigma}\backslash \{\overline{w}\})$ was defined in the proof of Theorem \ref{Mainthm1}; clearly it induces a biholomorphism \[\Phi\colon (M'\times[0,1], X')\to (\M\times[0,1], X).\] Naturality of the cohomology long exact sequence and Lemma \ref{M'Xsurj} then tell us that $\iota^*$ is a surjection so long as $r\ge1$. This argument also proves injectivity of $q^*$ on $\widetilde{H}^{2r+1}$ for $r\ge1$ as well, but we don't have any need for this result.  For $r=0$ injectivity of $q^*$ is immediate since the $0^{\textrm{th}}$ reduced cohomology of both $\M\times[0,1]$ and $\M\times[0,1]/X$ is zero.
 \end{proof}

Assuming Lemma \ref{chernsquare} we prove Proposition \ref{chernTNrel}.

\begin{proof}[Proof of Proposition \ref{chernTNrel}]
From functoriality of Chern classes \[q^*(c((E_T)_{rel}))=c(E_T)\textrm{ and } q^*(c((E_T)_{rel}))=c(E_T).\] By definition, $E_T=\pi^*(T\M)$ and $E_N=\pi^*(N\M)$ and so \[(\pi^*)^{-1}\circ q^*(c((E_T)_{rel}))=(\pi^*)^{-1}(c(E_T))=(\pi^*)^{-1}c(\pi^* T\M)=c(T\M)\] and  \[(\pi^*)^{-1}\circ q^*(c((E_N)_{rel}))=(\pi^*)^{-1}(c(E_N))=(\pi^*)^{-1}c(\pi^* N\M)=c(N\M).\]From Corollary \ref{chernTN} we know that \[c(T\M)=c(N\M)\] and therefore by injectivity of $(\pi^*)^{-1}\circ q^*$ (Lemma \ref{injection}) we conclude  \[c((E_T)_{rel})=c((E_N)_{rel}).\]
\end{proof}

Say that the non-zero cohomology groups of $\Sigma$ is given by 
   \[H^0(\Sigma)=\Z,\;H^1(\Sigma)=\bigoplus_{i=1}^{2g}\Z\langle a_i\rangle\oplus \Z\langle b_i\rangle,\; H^2(\Sigma)= \Z\langle c\rangle\]
   where $a_i\cdot a_j=b_i\cdot b_j=0$ and $a_i\cdot b_j=\delta_{ij}c$.  Let $e\colon \Sigma\backslash\{w,z\}\to \Sigma$ be the inclusion map. Then,  the non-zero cohomology of $\Sigma\backslash\{w,z\}$ is given by 
   \[H^0(\Sigma\backslash\{w,z\})=\Z,\;H^1(\Sigma\backslash\{w,z\})=\Z\langle\rho\rangle\oplus\bigoplus_{i=1}^{2g}\Z\langle \sigma_i\rangle\oplus \Z\langle \mu_i\rangle\]
   where $\sigma_i=e^*(a_i)$, $\mu_i=e^*(b_i)$ and $\rho$ is the hom dual of a homology class of a loop about $w$ or $z$. Let the non-zero cohomology of  $\overline{\Sigma}$ be 

\[H^0(\overline{\Sigma})=\Z,\; H^1(\overline{\Sigma})=\bigoplus_{i=1}^{g}\Z\langle \overline{a}_i\rangle\oplus \Z\langle \overline{b}_i\rangle,\; H^2(\overline{\Sigma})=\Z\langle \overline{c}\rangle\]
   where for $1\le i,j\le g$, $\overline{a}_i\cdot \overline{a}_j=\overline{b}_i\cdot \overline{b}_j=0$, $\overline{a}_i\cdot \overline{b}_j=\delta_{ij}\overline{c}$, and (recalling the quotient map $\ell\colon \Sigma\to \overline{\Sigma}$) $\ell^*(\overline{a}_i)=a_i+a_{i+g}$ and $\ell(\overline{b}_i)=b_i+b_{i+g}$. Finally, defining $\bar{e}\colon \overline{\Sigma}\backslash \{\overline{w}\}\to \overline{\Sigma}$ to be the inclusion map, the non-zero cohomology of $\overline{\Sigma}\backslash \{\overline{w}\}$ is given by
   \[H^0(\overline{\Sigma}\backslash \{\overline{w}\})=\Z\;, H^1( \overline{\Sigma}\backslash \{\overline{w}\})=\bigoplus_{i=1}^{g}\Z\langle \overline{\sigma}_i\rangle\oplus \Z\langle \overline{\mu}_i\rangle\]
   where $\overline{\sigma}_i=e^*(\overline{a}_i)$ and $\overline{\mu}_i=e^*(\overline{b}_i)$.
 With this notation in place, we record the following proposition about the cohomology and Chern classes of the $M$ and $M'$. For more about the cohomology ring $H^*(\Sym^{2g}(\Sigma))$ see \cite{Mac}.
\begin{proposition}
    The cohomology of $M=\Sym^{2g}(\Sigma\backslash \{w,z\})$ is given by \[H^*(M)\cong (H^*(\Sigma\backslash \{w,z\})^{\otimes 2g})^{S_{2g}},\] that is, the tensors in $H^*(\Sigma\backslash \{w,z\})^{\otimes 2g}$ which are invariant under the permutation action of the symmetric group $S_{2g}$. If for a cohomology class $u$, $u_k\coloneq 1\otimes\hdots\otimes \underbrace{u}_{\textup{kth slot}}\otimes\hdots\otimes 1$, then this ring is generated by the elements 
    \[\gamma\coloneq \sum_{k=1}^{2g}\rho_k,\;\epsilon_i\coloneq \sum_{k=1}^{2g}(\sigma_i)_k\textrm{ and } \chi_i\coloneq \sum_{k=1}^{2g}(\mu_i)_k.\]
    Similarly, the cohomology of $M'=\Sym^g(\overline{\Sigma}\backslash \{\overline{w}\})$ is given by \[H^*(M')\cong (H^*(\overline{\Sigma}\backslash\{\overline{w}\})^{\otimes g})^{S_g}\]which is generated by the elements 
     \[\overline{\epsilon}_i\coloneq \sum_{k=1}^{g}(\overline{\sigma}_i)_k,\textrm{ and }\overline{\chi}_i\coloneq \sum_{k=1}^{g}(\overline{\mu}_i)_k.\]
    The total Chern classes of $M$ and $M'$ are given by the formulas 
    \begin{equation}\label{chernclassM}
    c(TM)=\prod_{i=1}^{2g}(1-\epsilon_i\chi_i)\end{equation}and \begin{equation}\label{chernclassM'}c(TM')=\prod_{i=1}^g(1-\overline{\epsilon}_i\overline{\chi}_i).\end{equation}
\end{proposition}
\begin{proof}
 For the cohomology of $M$,  \cite[~Equation 1.2]{Mac} gives an isomorphism \[H^*(\Sym^{2g}(\Sigma\backslash \{w,z\}))\cong (H^*(\Sigma\backslash \{w,z\})^{\otimes 2g})^{S_{2g}}.\] This isomorphism is stated with cohomology coefficients in a characteristic $0$ field $K$, but it also applies with $\Z$ coefficients so long as the space in question has torsion free cohomology. It is a simple exercise to determine that $\gamma,\epsilon_i$ and $\chi_i$ for $1\le i\le 2g$ forms a basis for the $S_{2g}$ invariant tensors of $H^*(\Sigma\backslash \{w,z\}$.

In order to establish Equation (\ref{chernclassM}), we appeal to \cite[~(14.5)]{Mac} which states the following formula for the total Chern class of $\Sym^{2g}(\Sigma)$:
\[c(T\text{Sym}^{2g}(\Sigma))=(1+\eta)\prod_{i=1}^{2g}(1+\eta-\xi_i\zeta_{i})
\]
where \[\eta\coloneq \sum_{k=1}^{2g}c_k,\;\xi_i\coloneq \sum_{k=1}^{2g}(a_i)_k,\textrm{ and } \zeta_i\coloneq \sum_{k=1}^{2g}(b_i)_k.\]
 The induced inclusion map $\Sym^{2g}(e)^*\colon H^*(\text{Sym}^{2g}(\Sigma))\to H^*(M)$ acts by \[\Sym^{2g}(e)^*(\eta)=0,\;\Sym^{2g}(e)^*(\xi_i)=\epsilon_i\textrm{ and }\Sym^{2g}(e)^*(\zeta_i)=\chi_i.\]  Naturality of Chern classes then yields \ref{chernclassM}. The claims about $M'$ follow identically.
\end{proof}
We are now ready to prove Lemma \ref{chernsquare}, thus establishing Equation (\ref{relequal}) and verifying the existence of a stable tangent normal isomorphism.
\begin{proof}[Proof of Lemma \ref{chernsquare}]\label{chernsquareproof}Pullback along $i\circ \phi\colon M'\xhookrightarrow{}M$ acts on the generators of $\widetilde{H}^*(M)$ as follows. 
  \begin{itemize}
  \item \[(i\circ \phi)^*(\gamma)=0.\]
  \item For $1\le i\le g$ we have  \[(i\circ \phi)^*(\epsilon_i)=(i\circ \phi)^*(\epsilon_{i+g})=\overline{\epsilon}_i\] and \[(i\circ \phi)^*(\chi_i)=(i\circ \phi)^*(\chi_{i+g})=\overline{\chi}_i.\]
  \end{itemize}
Therefore, from equations \ref{chernclassM} and \ref{chernclassM'} we see that pullback of the total Chern class of $M$ along $i\circ \phi$ is given by \[\phi^*i^*(c(TM))=c((i\circ \phi)^*TM)=c(TM')^2=\phi^*c(T\M)^2.\]The result $i^*(c(TM))=c(T\M)^2$ then follows because $\phi^*$ is invertible.
 \end{proof}
\appendix
\section{Grid homology and Alexander gradings}\label{Appendix}
Let $\mathcal{H}_S$ be a $\tau$-equivariant Heegaard diagram for $S^n_b(K)$ and let $\overline{\mathcal{H}}\coloneq \mathcal{H}_S/\tau$ be the induced quotient Heegaard diagram for $\overline{K}$. The goal of this appendix is to determine the relation between the Alexander grading of a $\tau$-equivariant generator $x\in\widetilde{CFK}(\mathcal{H}_S)$ and the Alexander grading of its image $x\in \widetilde{CFK}(\overline{\mathcal{H}})$, thereby determining the Alexander grading shift of Theorem \ref{Mainthm2}. 

The strategy we employ is as follows. We introduce a class of Heegaard diagrams that we call \emph{spherical grid diagrams}, which are essentially grid diagrams on genus $0$ surfaces. Spherical grid diagrams diagrams are especially suited to studying the knot Floer homology of strongly invertible and doubly periodic knots/links combinatorially since the $\tau$ action on a symmetric spherical grid diagram does not fix any $\alpha$ or $\beta$ curves, unlike the $\tau$ action on a grid diagram. In particular, we introduce $\tau$-equivariant spherical grid diagrams $\mathbb{G}_S^\textrm{sp}$ for $S^n_b(K)$, and $Ax\mathbb{G}_{S_0}^\textrm{sp}$ for $L_b^{n+\frac{1}{2}}(K)\cup A$. The skein exact triangle in singular link Floer homology allows us to determine the precise relation (Proposition \ref{alexandergradingstatement2}) between the Alexander gradings of $\widetilde{CFK}(\mathbb{G}_S^{\textrm{sp}})$ and the Alexander gradings of $\widetilde{CFK}(Ax\mathbb{G}_{S_0}^\textrm{sp})$. In addition, the relation between the Alexander gradings in $\widetilde{CFK}(Ax\mathbb{G}_{S_0}^\textrm{sp})$ and $\widetilde{CFK}(Ax\mathbb{G}_{S_0}^\textrm{sp}/\tau)$ are determined by Hendricks in \cite[~Section 3]{HendricksDP}, since $Ax\mathbb{G}_{S_0}^\textrm{sp}$ is a $\tau$-equivariant Heegaard diagram for the doubly periodic knot $L_b^{n+\frac{1}{2}}(K)$. In this admittedly roundabout way, we are able to pin down the relation between Alexander gradings of $\widetilde{HFK}(\mathbb{G}_S^\textrm{sp})$ and $\widetilde{HFK}(\mathbb{G}_S^\textrm{sp}/\tau)$.

\subsection{Singular skein triangles from grid homology}\label{gridsection} 
In direct analogy to the skein triangle of link Floer homology of \cite[~Section 8]{OSknots}, there are \emph{singular skein exact triangles} that recover the singular skein relations Equations (\ref{posskein}) and (\ref{negskein}) upon taking bi-graded Euler characteristics. 
\begin{proposition}\cite[~Theorem 4.1]{OSsing2}\label{skeintriangles}
If $(S, S_0, S_-)$ is a negative singular skein triple, then for $m$ sufficiently large there is an exact triangle

\begin{center}
\begin{tikzpicture}
% Define the matrix
\matrix (m) [matrix of math nodes, row sep=2em, column sep=2em]
{
     \HFKhat(S_-)\otimes V^{\otimes (m-1)} && \HFKhat(S_0)\otimes V^{\otimes m} \\
    & \widehat{HFK}(S)\otimes V^{\otimes m} [[0,-\frac{1}{2}]]&\\
};
% Draw the arrows
\draw[->] (m-1-1) -- (m-1-3) node[midway, above] {$(0,0)$};
\draw[->] (m-2-2) -- (m-1-1) node[midway, below left] {$(0,0)$};
\draw[->] (m-1-3) -- (m-2-2) node[midway, below right] {$(-1,0)$};
\end{tikzpicture}.
\end{center}
\end{proposition}
For a bigraded module $C$, the notation $C[[a,b]]$ denotes the grading shifted module $C[[a,b]]_{(d,s)}\coloneq C_{d+a,b+s}$; in the context of knot Floer homology, the first grading is the Maslov grading and the second is the Alexander grading.
The grid homology based proof of Proposition \ref{skeintriangles}, originally given in \cite[~Section 6]{OSsing2} modulo gradings, is very similar to the ordinary link Floer homology case. After giving a quick review of grid homology, we spell it out in detail below  mainly to keep track of Alexander gradings, but also as a warm-up for the similar proof given in the following subsection using spherical grid diagrams. We do not attempt show that the Maslov grading shift for the module at the bottom of the triangle is $0$, but content ourselves with a proof that it is even, as our main focus is on Alexander gradings. The proof is based on the following lemma. A singly graded version of this lemma is proved in \cite[~Lemma A.3.10]{OSSgridhomologybook}, but the proof is easily adapted to the bigraded setting.
\begin{lemma}\label{conetriangle}
    Suppose that $C, C'$
, and $C''$ are three bigraded chain complexes, and $f \colon  C' \to C''$ and $g \colon  C \to C'$ are chain maps that are homogeneous of
degrees $(a,p)$ and $(b,q)$ respectively. Then, there is a chain map  $\Phi\colon \textrm{Cone}(f) \to \textrm{Cone}(g)$
which is homogeneous of degree $(-a-1,-p)$ and whose induced map on homology fits into the following exact triangle with bidegrees as indicated\\
\begin{center}
\begin{tikzpicture}
% Define the matrix
\matrix (m) [matrix of math nodes, row sep=2em, column sep=2em]
{
    H(\textrm{Cone}(f)) && H(\textrm{Cone}(g)) \\
    & H(\textrm{Cone}(f\circ g)) &\\
};
% Draw the arrows
\draw[->] (m-1-1) -- (m-1-3) node[midway, above] {$(-a-1,-p)$};
\draw[->] (m-2-2) -- (m-1-1) node[midway, below left] {$(0,0)$};
\draw[->] (m-1-3) -- (m-2-2) node[midway, below right] {$(a,p)$};
\end{tikzpicture}
\end{center}
\end{lemma}

We give a lightning review of a combinatorial reformulation of link Floer homology called grid homology here. For many more details the reader is referred to \cite{OSSgridhomologybook}. 

A \emph{grid diagram} $\mathbb{G}$ for a possibly singular link $L$ is an $m\times m$ satisfying the following properties:
\begin{enumerate}
    \item Every box is empty, has an $X$, or has an $O$.
    \item There is exactly one $X$ in every row and column
    \item In the case that $L$ is not singular, there is exactly one $O$ in every row and column.
    \item In the case that $L$ is singular, there are either one or two $O$'s
 in every row and column. 
    \item Connecting $O$'s to $X$'s in the columns with oriented vertical segments, and connecting $X's$ to $O's$ with oriented horizontal segments that pass below the vertical segments, we obtain a diagram $D(\mathbb{G})$ for the (singular) link $L$.
 \end{enumerate}In the case that $L$ is singular, $X's$ which have two $O's$ in their row and column are instead written as $XX$. Identifying the top and bottom gridlines, and right and left gridlines on $\mathbb{G}$, thinking of the horizontal and vertical grid lines as the $\alpha$ and $\beta$ curves, and calling replacing $X/O\to z/w$, we see that a grid diagram is a genus 1 multipointed Heegaard diagram for $L$ in the sense of Section \ref{singularlinkFloer}. It is typical to use the notation $\widetilde{GC}(\mathbb{G})\coloneq\widetilde{CFK}(\mathbb{G})$ and to call a generator of $\widetilde{CFK}$ a \emph{grid state}. A grid diagram is an example of a \emph{nice Heegaard diagram}, meaning that all of the elementary regions that don't contain a basepoint are rectangles or bigons. One can compute the differential of the Floer chain complex associated to a nice Heegaard diagram purely combinatorially by work of Sarkar and Wang \cite{SWNice}. This incarnation of link Floer homology was elaborated by Manolescu, Ozsv\'{a}th and Thurston so that many features of link Floer homology have explicit formulas in terms of combinatorics of the grid diagram $\mathbb{G}.$ In particular, there are formulas for the Alexander and Maslov grading of grid states  \emph{in the case that}  $L$ \emph{is non-singular} as follows. Let $x^{NWO}$ and $x^{NWX}$ denote the grid states obtained by marking the top left corner of each square containing an $O$ or $X$ respectively. Notice that $x^{NWO}$ is only a valid grid state if $L$ is non-singular. Define functions $M_{\mathbb{O}}, M_{\mathbb{X}}\colon \mathbf{S}(\mathbb{G})\to \Z$ by setting \begin{equation}\label{Maslov normalization}M_{\mathbb{O}}(x^{NWO})=M_{\mathbb{X}}(x^{NWX})=0,\end{equation}
\begin{equation}\label{relativeOMaslov}
M_{\mathbb{O}}(x)-M_{\mathbb{O}}(y)=1-2|r\cap \mathbb{O}|+2|x\cap \textrm{Int}(r)|
\end{equation}
and
\begin{equation}\label{relativeXMaslov}
M_{\mathbb{X}}(x)-M_{\mathbb{X}}(y)=1-2|r\cap \mathbb{X}|+2|x\cap \textrm{Int}(r)|
\end{equation}
where $x$ and $y$ are grid states that differ at exactly two points and $r\in Rect(x,y)$ is any rectangle with bottom left and top right corners in $x$, and bottom right and top left corners in $y.$ Then the Maslov grading $M(x)$ is given by \begin{equation}\label{Maslovgrid}M(x)\coloneq M_{\mathbb{O}}(x)\end{equation} and if $L$ has $\ell$ components, the Alexander grading $A(x)$ is given by \begin{equation}\label{Alexandergrid}A(x)\coloneq \frac{1}{2}(M_{\mathbb{O}}(x)-M_{\mathbb{X}}(x))-\frac{m-\ell}{2}.\end{equation}

\begin{figure}[htb!]
\center

\begin{tikzpicture}[scale=1, every node/.style={draw, minimum size=1cm}]

  % Grid
  \draw (0,0) grid (3,3);

  % Markings
  \foreach \i/\j/\marking in {
    1/1/, 1/2/O, 1/3/, 
    2/1/O, 2/2/XX, 2/3/O, 
    3/1/, 3/2/O, 3/3/
  } {
    \node at (\i-0.5, \j-0.5) {\marking};
  }
  \draw[red] (1.5, .5) -- (1.5, 2.5);
  \draw[red] (.5,1.5)-- (2.5, 1.5);

\end{tikzpicture}

\begin{tikzpicture}[scale=1, every node/.style={draw, minimum size=1cm}]

  % Grid
  \draw (0,0) grid (4,4);

  % Markings
  \foreach \i/\j/\marking in {
    1/1/, 1/2/O, 1/3/,1/4/, 
    2/1/, 2/2/X, 2/3/,2/4/O,
    3/1/O, 3/2/, 3/3/X,3/4/,
    4/1/, 4/2/, 4/3/O, 4/4/
  } {
    \node at (\i-0.5, \j-0.5) {\marking};
  }
  \draw[red] (.5,1.5) -- (1.5,1.5);
  \draw[red] (2.5, 2.5) -- (3.5, 2.5);
  \draw[red] (1.5, 1.5) -- (1.5, 3.5);
  \draw[red] (2.5, .5) -- (2.5, 2.5);
\fill[color=cyan] (2,2) circle (0.1);
\end{tikzpicture}
\begin{tikzpicture}[scale=1, every node/.style={draw, minimum size=1cm}]

  % Grid
  \draw (0,0) grid (4,4);

  % Markings
  \foreach \i/\j/\marking in {
    1/1/, 1/2/O, 1/3/,1/4/, 
    2/1/, 2/2/, 2/3/X,2/4/O,
    3/1/O, 3/2/X, 3/3/,3/4/,
    4/1/, 4/2/, 4/3/O, 4/4/
  } {
    \node at (\i-0.5, \j-0.5) {\marking};
  }
  \draw[red] (.5, 1.5) -- (2.5, 1.5);
  \draw[red] (1.5, 2.5) -- (3.5, 2.5);
  \draw[red] (1.5,2.5) -- (1.5, 3.5);
  \draw[red] (2.5, .5) -- (2.5, 1.5);
\fill[color=cyan] (2,2) circle (0.1);
\end{tikzpicture}

\caption{Top: The center $3\times3$ subgrid of the grid diagram $\mathbb{G}_S$ for a singular knot $S$. Bottom left: Local modifications to move from $\mathbb{G}_S$ to a diagram $\mathbb{G}_{S_-}$ for the negative resolution $S_-$. Bottom right: Local modifications to move from $\mathbb{G}_S$ to a diagram $\mathbb{G}_{S_0}$ for the zero resolution $S_0$}\label{centercross}
\end{figure}

\begin{proof}[Proof of Proposition \ref{skeintriangles}] 
Consider an $m\times m$ grid diagram $\mathbb{G}_S$ for a singular knot $S$ which contains a $3\times 3$ subgrid decorated as in the top of Figure \ref{centercross}. Performing the local modifications seen in the bottom right and left parts of Figure \ref{centercross} to $\mathbb{G}_S$ yields $(m+1)\times (m+1)$ grid diagrams, $\mathbb{G}_{S_-}$ for the negative resolution $S_-$ of $S$, and $\mathbb{G}_{S_0}$ for the zero resolution $S_0$ of $S$.

Let $\mathbf{S}_-$ and $\mathbf{S}_0$ denote the set of generators/grid states of $\widetilde{GC}(\mathbb{G}_{S_-})$ and $\widetilde{GC}(\mathbb{G}_{S_0})$ respectively. Partition $\widetilde{GC}(\mathbb{G}_{S_-})$ and $\widetilde{GC}(\mathbb{G}_{S_0})$ as follows: $\widetilde{GC}(\mathbb{G}_{S_-})=I_-\cup N_-$ and $\widetilde{GC}(\mathbb{G}_{S_0})=I_0\cup N_0$ where $I_-$ and $I_0$ are generated by those states containing the center point $c$ (marked in blue in Figure \ref{centercross}), and $N_-$ and $N_0$ are generated by those states not containing $c.$ The subspaces $I_-$ and $N_0$ are subcomplexes while the subspaces $I_0$ and $N_-$ are quotient complexes, and hence all of these subspaces inherit differentials from the complexes $\widetilde{GC}(\mathbb{G}_{S_-})$ and $\widetilde{GC}(\mathbb{G}_{S_0})$; we write the resulting complexes as $(I_-, \partial_{I_-}^{I_-})$, $(N_-, \partial_{N_-}^{N_-})$, $(I_0, \partial_{I_0}^{I_0})$, and $(N_0, \partial_{N_0}^{N_0})$. Consider also the three chain maps \[\partial_{I_0}^{N_0}\colon (I_0, \partial_{I_0}^{I_0})\to (N_0, \partial_{N_0}^{N_0})\]\[\partial_{N_-}^{I_-}\colon (N_-, \partial_{N_-}^{N_-})\to (I_-, \partial_{I_-}^{I_-})\]\[\psi\colon (N_0, \partial_{N_0}^{N_0})\to (N_-, \partial_{N_-}^{N_-})\]where the first two are induced by the differentials on $\widetilde{GC}(\mathbb{G}_{S_-})$ and $\widetilde{GC}(\mathbb{G}_{S_0})$ respectively while $\psi$ is the map identifying grid states missing $c$ in $\widetilde{GC}(\mathbb{G}_{S_0})$ with grid states missing $c$ in $\widetilde{GC}(\mathbb{G}_{S_-}).$ Now compose these chain maps: 
\begin{equation}
    (I_0, \partial_{I_0}^{I_0})\xrightarrow[(-1,0)]{\partial_{I_0}^{N_0}}(N_0, \partial_{N_0}^{N_0})\xrightarrow[(0,0)]{\psi}(N_-, \partial_{N_-}^{N_-})\xrightarrow[(-1,0)]{\partial_{N_-}^{I_-}}(I_-, \partial_{I_-}^{I_-}).
\end{equation}
The tuples under the arrows indicate the bigraded degree of each chain map. The degree shift associated to $\partial_{I_0}^{N_0}$ and $\partial_{N_-}^{I_-}$ is $(-1,0)$ because both of these differentials are induced from the standard differential on the grid chain complex which has bi-degree $(-1,0)$. We use the grid Maslov and Alexander grading formulas, Equations \ref{Maslovgrid} and \ref{Alexandergrid}, to prove the bi-degree of $\psi$ is $(0,0)$.  Let $x_-^{NWX}$ and $x_0^{NWX}$ be the grid states with markings on the northwest corner of the $X$ markings in $\mathbb{G}_{S_-}$ and $\mathbb{G}_{S_0}$ respectively. Notice that $x_-^{NWX}$ is an element of $N_0$ and $N_-$ because it does not contain the center point $c$. Comparing the Alexander and Maslov gradings of $x_-^{NWX}$ considered as an element of $\mathbf{S}_0$ and $\mathbf{S}_-$ will therefore tell us the bi-degree of $\psi$. Firstly notice that $\mathbb{G}_{S_-}$ and $\mathbb{G}_{S_0}$ have the same $\mathbb{O}$ sets and hence $M_{\mathbb{O},-}$ and $M_{\mathbb{O},0}$ (the extra subscript labels the domain of the function, either $\mathbf{S_-}$ or $\mathbf{S_0}$) take the same value on $x_-^{NWX}$. This verifies the claim that the Maslov grading shift of $\psi$ is $0$. Next, using equations \ref{Maslov normalization}, \ref{relativeXMaslov} and the fact that the $1\times1$ square $NW_c$ to the northwest of $c$ is a rectangle going from $x_-^{NWX}$ to $x_0^{NWX}$, we see that 
\begin{equation}\label{MX-x0}
    -M_{\mathbb{X},-}(x_0^{NWX})=M_{\mathbb{X},-}(x_-^{NWX})-M_{\mathbb{X},-}(x_0^{NWX})=1-0+0=1\implies M_{\mathbb{X},-}(x_0^{NWX})=-1
\end{equation}
since $NW_c$ in $\mathbb{G}_{S_-}$ contains no $X$ marking, and 
\begin{equation}\label{MX0x-}
M_{\mathbb{X},0}(x_-^{NWX})=M_{\mathbb{X},0}(x_-^{NWX})-M_{\mathbb{X},0}(x_0^{NWX})=1-2\cdot 1+0=-1. 
\end{equation}
since $NW_c$ in $\mathbb{G}_{S_0}$ contains one $X$ marking.
We only need Equation (\ref{MX0x-}) to compute the Alexander grading shift of $\psi,$ but Equation (\ref{MX-x0}) will be used later in this proof. Now we may compare the Alexander gradings of $x_-^{NWX}$ in $\mathbf{S_-}$ and $\mathbf{S_0}$ by using Equation (\ref{Alexandergrid}).
\begin{equation}\label{alex0x-}
     A_{0}(x_-^{NWX})=\frac{1}{2}(M_{\mathbb{O}, 0}(x_-^{NWX})-M_{\mathbb{X},0}(x_-^{NWX}))-\frac{(m+1)-1}{2}
\end{equation}
\begin{equation}\label{alex-x-}
    A_{-}(x_-^{NWX})=\frac{1}{2}(M_{\mathbb{O}, -}(x_-^{NWX})-M_{\mathbb{X},-}(x_-^{NWX}))-\frac{(m+1)-2}{2}
\end{equation}
Recalling Equations \ref{Maslov normalization} and \ref{MX0x-} along with the fact that $M_{\mathbb{O},0}(x_-^{NWX})=M_{\mathbb{O},-}(x_-^{NWX})$, we have shown  $A_{\mathbb{X},0}(x_-^{NWX})=A_{\mathbb{X},0}(x_-^{NWX})$. The Alexander grading shift of $\psi$ is 0 as claimed.

In Lemma \ref{conetriangle}, let $f=\partial_{N_-}^{I_-}$ and $g=\psi\circ \partial_{I_0}^{N_0}$. Then we have $deg(f)=(-1,0)$, and hence the following exact triangle is established: 
\begin{center}
\begin{tikzpicture}
% Define the matrix
\matrix (m) [matrix of math nodes, row sep=2em, column sep=2em]
{
    H(\textrm{Cone}(\partial_{N_-}^{I_-})) && H(\textrm{Cone}(\psi\circ\partial_{I_0}^{N_0})) \\
    & H(\textrm{Cone}(\partial_{N_-}^{I_-}\circ \psi\circ\partial_{I_0}^{N_0})) &\\
};
% Draw the arrows
\draw[->] (m-1-1) -- (m-1-3) node[midway, above] {$(0,0)$};
\draw[->] (m-2-2) -- (m-1-1) node[midway, below left] {$(0,0)$};
\draw[->] (m-1-3) -- (m-2-2) node[midway, below right] {$(-1,0)$};
\end{tikzpicture}
\end{center}
There are isomorphisms \begin{equation}H(\textrm{Cone}(\psi\circ\partial_{I_0}^{N_0}))\cong H(\textrm{Cone}(\partial_{I_0}^{N_0}))\cong \widetilde{GH}(\mathbb{G}_{S_0})\end{equation}
where the first follows because $\psi$ is an isomorphism of bi-degree $(0,0)$, and the second follows because the differential on $\textrm{Cone}(\partial_{I_0}^{N_0})$ is by definition $\begin{bmatrix}
    \partial_{I_0}^{I_0}      & 0  \\
    \partial_{N_0}^{I_0}       & \partial_{N_0}^{N_0}
\end{bmatrix}$ which is the full differential on $\widetilde{GH}(\mathbb{G}_{S_0})=I_0\oplus N_0$ since $N_0$ is a subcomplex ($\partial^{N_0}_{I_0}=0$). Similarly there is an isomorphism \begin{equation}H(\textrm{Cone}(\partial_{N_-}^{I_-}))\cong \widetilde{GH}(\mathbb{G}_{S_-}).\end{equation}Also notice that $\partial_{N_-}^{I_-}\circ \psi\circ\partial_{I_0}^{N_0}\colon (I_0, \partial_{I_0}^{I_0})\to(I_-, \partial_{I_-}^{I_-})$ is the zero map and has degree $(-2,0).$ We therefore have the following equalities and isomorphisms of chain complexes: 
\begin{equation}\textrm{Cone}(\partial_{N_-}^{I_-}\circ \psi\circ\partial_{I_0}^{N_0})_{(d,s)}=(I_0)_{(d+1,s)}\oplus (I_-)_{(d,s)}\cong (I_0)_{(d+1,s)}\oplus (I_0)_{(d,s-1)}=(I_0\otimes V)[[0,-1]]_{(d,s)}.\end{equation} 
The bidegree of the isomorphism $(I_-)_{(d,s)}\cong (I_0)_{(d,s-1)}$ follows from our previous observation that the $\mathbb{O}$ sets of $\mathbb{G}_{S_-}$ and $\mathbb{G}_{S_0}$ being identical implies $M_{\mathbb{O},-}=M_{\mathbb{O},0}$, and from Equation (\ref{MX-x0}). 
Putting these identifications into the exact triangle yields
\begin{center}
\begin{tikzpicture}
% Define the matrix
\matrix (m) [matrix of math nodes, row sep=2em, column sep=2em]
{
    \widetilde{GH}(\mathbb{G}_{S_-}) && \widetilde{GH}(\mathbb{G}_{S_0}) \\
    & (H(I_0)\otimes V)[[0,-1]] &\\
};
% Draw the arrows
\draw[->] (m-1-1) -- (m-1-3) node[midway, above] {$(0,0)$};
\draw[->] (m-2-2) -- (m-1-1) node[midway, below left] {$(0,0)$};
\draw[->] (m-1-3) -- (m-2-2) node[midway, below right] {$(-1,0)$};
\end{tikzpicture}
\end{center}

Taking bi-graded Euler characteristics (which we denote by the symbol $\chi$) of this exact triangle we get
\begin{equation}\label{eulercharoftri}\chi(\widetilde{GH}(\mathbb{G}_{S_-}))=\chi(\widetilde{GH}(\mathbb{G}_{S_0}))+\chi((H(I_0)\otimes V)[[0,-1]]).\end{equation}
Equations \ref{HFKdependencyonbp} and \ref{HFLeulerchar} tell us that\[\chi(\widetilde{GH}(\mathbb{G}_{S_-}))= \chi(\HFKhat(S_-)\otimes V^{(m+1)-2})=\Delta_{S_-}(t)(t^{\frac{1}{2}}-t^{-\frac{1}{2}})(1-t^{-1})^{m-1}\]
and similarly 
\[\chi(\widetilde{GH}(\mathbb{G}_{S_0}))=\chi(\HFKhat(S_0)\otimes V^{(m+1)-1})=\Delta_{S_0}(t)(1-t^{-1})^{m}.\]In addition,
\[\chi((H(I_0)\otimes V)[[0,-1]])=\chi(H(I_0))(1-t^{-1})\cdot t.\]
Making these replacements in Equation (\ref{eulercharoftri}) and dividing through by $t^{\frac{1}{2}}(1-t^{-1})^m$ yields
\begin{equation}\label{eqq}
    \Delta_{S_-}(t)=\frac{t^{\frac{1}{2}}\chi(H(I_0))}{(1-t)^{m-1}}+t^{-\frac{1}{2}}\Delta_{S_0}(t).
\end{equation}
Comparing Equation (\ref{eqq}) to Equation (\ref{negskein}) shows that \begin{equation}\label{eqq2}\frac{t^{\frac{1}{2}}\chi(H(I_0))}{(1-t)^{m-1}}=\Delta_{S}.\end{equation}
Substituting $t^{\frac{1}{2}}\chi(H(I_0))=\chi(H(I_0)[[0,-\frac{1}{2}]])$ and $\chi(\widetilde{GH}(\mathbb{G}_S))=\Delta_S(t)(1-t^{-1})^{m-1}$ into Equation (\ref{eqq2}) we conclude

\begin{equation}\label{eulerchars}
    \chi(H(I_0)[[0,-\frac{1}{2}]])=\chi(\widetilde{GH}(\mathbb{G}_S)).
\end{equation}
The complex $\widetilde{GC}(\mathbb{G}_S)$ is, up to a bidegree shift, isomorphic to the complex $I_0$ through an isomorphism \[F\colon \widetilde{GC}(\mathbb{G}_S)\to I_0\] taking a generator in $\widetilde{GC}(\mathbb{G}_S)$ and mapping it to the corresponding generator in $I_0$ that contains the center point $c$. More explicitly, we have

\begin{equation*}
\begin{split}
F(\{(1,x_1), \hdots,& (m,x_m)\})= \\
& \{(1,x_1),\hdots,(\frac{m-1}{2}, x_{\frac{m-1}{2}}), (\frac{m+1}{2}, \frac{m+1}{2}), (\frac{m+3}{2}, x_{\frac{m+1}{2}}), \hdots,(m+1, x_m)\}.
\end{split}
\end{equation*}
By Equation (\ref{eulerchars}) $F$ increases Alexander grading by $\frac{1}{2}$ and changes Maslov grading by some even integer $2j$. Making the replacement $H(I_0)\cong \widetilde{GH}(\mathbb{G}_S)[[2j,\frac{1}{2}]]$ in the exact triangle  gives us
\begin{center}
\begin{tikzpicture}
% Define the matrix
\matrix (m) [matrix of math nodes, row sep=2em, column sep=2em]
{
    \widetilde{GH}(\mathbb{G}_{S_-}) && \widetilde{GH}(\mathbb{G}_{S_0}) \\
    & (\widetilde{GH}(\mathbb{G}_S)\otimes V)[[2j,-\frac{1}{2}]] &\\
};
% Draw the arrows
\draw[->] (m-1-1) -- (m-1-3) node[midway, above] {$(0,0)$};
\draw[->] (m-2-2) -- (m-1-1) node[midway, below left] {$(0,0)$};
\draw[->] (m-1-3) -- (m-2-2) node[midway, below right] {$(-1,0)$};
\end{tikzpicture}
\end{center}
which is the exact triangle of Proposition \ref{skeintriangles} (besides showing $j=0$) after identifying grid homology with (singular) knot Floer homology.

\end{proof}

The main reason that we gave a proof of Proposition \ref{skeintriangles} was the following corollary that was shown along the way:
\begin{corollary}\label{gradingshift}
    Let $\mathbb{G}_S$ and $I_0$ be as in the proof of Proposition \ref{skeintriangles}. Then there is an isomorphism of chain complexes
    \begin{equation}
    I_0\cong \widetilde{GC}(\mathbb{G}_S)[[2j,\frac{1}{2}]]
    \end{equation}
    for some integer $j$.
\end{corollary}

Although the skein exact triangle of Proposition \ref{skeintriangles} is completely general, from now on we have in mind the situation that $\mathbb{G}_S$ is a $(2k+1)\times (2k+1)$ grid diagram for $S^n_b(K)$ which is invariant under rotation by $180$ about the center square, and such that the central $3\times3$ subgrid of $\mathbb{G}_S$ is given by the top of Figure \ref{centercross}. Then $\mathbb{G}_{S_-}$ is a grid diagram for $S^n_b(K)_-=L^b_n(K)$ and $\mathbb{G}_{S_0}$ is a grid diagram for $S^n_b(K)_0=L_b^{n+\frac{1}{2}}(K)$. In Figure \ref{fig:gridsSbLb} such a triple of grid diagrams associated to $S^1_b(3_1^+), L^1_b(3_1^+)$ and $L^{\frac{3}{2}}_b(3_1^+)$ is displayed.

\begin{figure}[htb!]
\center

\begin{tikzpicture}[scale=1, every node/.style={draw, minimum size=1cm}]

  % Grid
  \draw (0,0) grid (5,5);

  % Markings
  \foreach \i/\j/\marking in {
    1/1/, 1/2/X, 1/3/, 1/4/, 1/5/O,
    2/1/X, 2/2/, 2/3/O, 2/4/, 2/5/,
    3/1/, 3/2/O, 3/3/XX, 3/4/O, 3/5/,
    4/1/, 4/2/, 4/3/O, 4/4/, 4/5/X,
    5/1/O, 5/2/, 5/3/, 5/4/X, 5/5/
  } {
    \node at (\i-0.5, \j-0.5) {\marking};
  }
    \draw[green] (.5,1.5) -- (.5, 4.5);
  \draw[green] (1.5,.5) -- (1.5, 2.5);
  \draw[red] (2.5,1.5) -- (2.5, 3.5);
  \draw[green] (3.5,2.5) -- (3.5, 4.5);
  \draw[green] (4.5,.5) -- (4.5, 3.5);
  \draw[green] (.5,4.5) -- (3.5, 4.5);
  \draw[green] (2.5,3.5) -- (3.25, 3.5);
  \draw[green] (3.75, 3.5) -- (4.5, 3.5);
  \draw[red] (1.5, 2.5) -- (3.5, 2.5);
  \draw[green] (.5, 1.5) -- (1.25, 1.5);
  \draw[green] (1.75, 1.5) -- (2.5, 1.5);
  \draw[green] (1.5, 0.5) -- (4.5, 0.5);

\end{tikzpicture}

\begin{tikzpicture}[scale=1, every node/.style={draw, minimum size=1cm}]

  % Grid
  \draw (0,0) grid (6,6);

  % Markings
  \foreach \i/\j/\marking in {
    1/1/, 1/2/X, 1/3/, 1/4/, 1/5/,1/6/O,
    2/1/X, 2/2/, 2/3/O, 2/4/, 2/5/,2/6/,
    3/1/, 3/2/, 3/3/X, 3/4/, 3/5/O,3/6/,
    4/1/, 4/2/O, 4/3/, 4/4/X, 4/5/,4/6/,
    5/1/, 5/2/, 5/3/, 5/4/O, 5/5/,5/6/X,
     6/1/O, 6/2/, 6/3/, 6/4/, 6/5/X,6/6/
  } {
    \node at (\i-0.5, \j-0.5) {\marking};
  }
  \draw[green] (.5, 1.5) -- (.5, 5.5);
  \draw[green] (1.5, .5) -- (1.5, 2.5);
  \draw[red] (2.5, 2.5) -- (2.5, 4.5);
  \draw[red] (3.5, 1.5) -- (3.5, 3.5);
  \draw[green] (4.5, 3.5) -- (4.5, 5.5);
  \draw[green] (5.5, .5) -- (5.5, 4.5);
  \draw[green] (.5,5.5) -- (4.5, 5.5);
  \draw[green] (2.5, 4.5) -- (4.25, 4.5);
  \draw[green] (4.75, 4.5) -- (5.5, 4.5);
  \draw[red] (3.5, 3.5) -- (4.5, 3.5);
  \draw[red] (1.5, 2.5) -- (2.5, 2.5);
  \draw[green] (.5, 1.5) -- (1.25, 1.5);
  \draw[green] (1.75, 1.5) -- (3.5, 1.5);
  \draw[green] (1.5, .5) -- (5.5, .5);
  \fill[color=cyan] (3,3) circle (0.1);

\end{tikzpicture}
\begin{tikzpicture}[scale=1, every node/.style={draw, minimum size=1cm}]

  % Grid
  \draw (0,0) grid (6,6);

  % Markings
  \foreach \i/\j/\marking in {
    1/1/, 1/2/X, 1/3/, 1/4/, 1/5/,1/6/O,
    2/1/X, 2/2/, 2/3/O, 2/4/, 2/5/,2/6/,
    3/1/, 3/2/, 3/3/, 3/4/X, 3/5/O,3/6/,
    4/1/, 4/2/O, 4/3/X, 4/4/, 4/5/,4/6/,
    5/1/, 5/2/, 5/3/, 5/4/O, 5/5/,5/6/X,
     6/1/O, 6/2/, 6/3/, 6/4/, 6/5/X,6/6/
  } {
    \node at (\i-0.5, \j-0.5) {\marking};
  }
  \draw[green] (.5, 1.5) -- (.5, 5.5);
  \draw[green] (1.5, .5) -- (1.5, 2.5);
  \draw[red] (2.5, 3.5) -- (2.5, 4.5);
  \draw[red] (3.5, 1.5) -- (3.5, 2.5);
  \draw[green] (4.5, 3.5) -- (4.5, 5.5);
  \draw[green] (5.5, .5) -- (5.5, 4.5);
  \draw[green] (.5,5.5) -- (4.5, 5.5);
  \draw[green] (2.5, 4.5) -- (4.25, 4.5);
  \draw[green] (4.75, 4.5) -- (5.5, 4.5);
  \draw[red] (2.5, 3.5) -- (4.5, 3.5);
  \draw[red] (1.5, 2.5) -- (3.5, 2.5);
  \draw[green] (.5, 1.5) -- (1.25, 1.5);
  \draw[green] (1.75, 1.5) -- (3.5, 1.5);
  \draw[green] (1.5, .5) -- (5.5, .5);
\fill[color=cyan] (3,3) circle (0.1);
\end{tikzpicture}

\caption{Top: A grid diagram $\mathbb{G}_S$ for the singular knot $S^1_b(3_1^+)$. Bottom left: The grid diagram $\mathbb{G}_{S_-}$ for the Hopf link $L_b^{1}(3_1^+)$. Bottom right: The grid diagram $\mathbb{G}_{S_0}$ for the unknot $L_b^{\frac{3}{2}}(3_1^+)$.} \label{fig:gridsSbLb}
\end{figure}

 \subsection{Skein triangles from spherical grid diagrams}

\begin{figure}[htb!]
 \center
 \includegraphics[scale=.7]{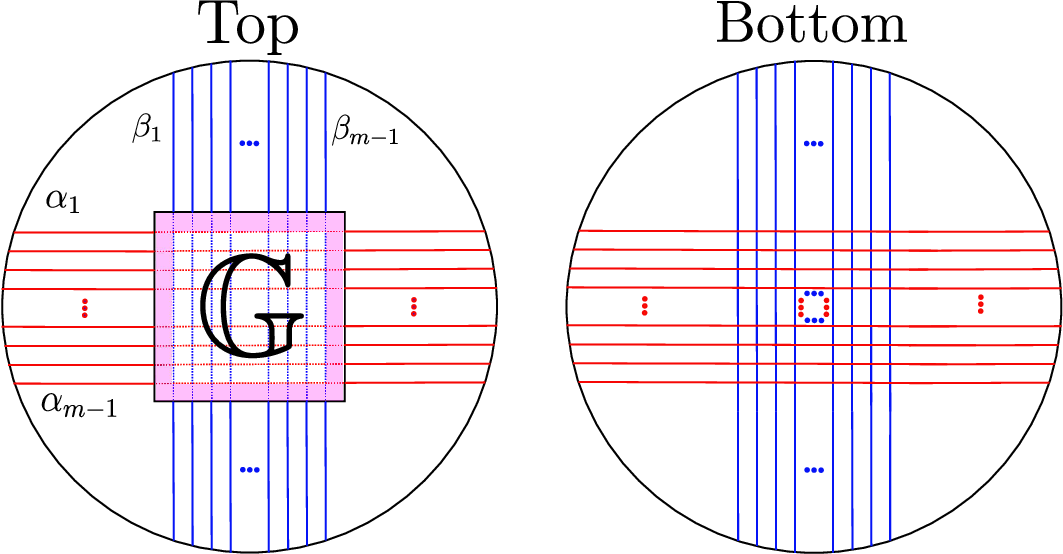}
 \caption{The genus $0$ Heegaard diagram $\mathbb{G}^{\textrm{sp}}$}\label{sphericalgriddiagram}
 \end{figure}

Say that we are given a $m\times m$ grid diagram $\mathbb{G}$, which is in particular a multi-pointed genus $1$ Heegaard diagram with  $|\bm{\alpha}|=|\bm{\beta}|=m$. In \cite[~Section 3.1]{MR2443251}, starting from $\mathbb{G}$, V\'{e}rtesi constructs a multi-pointed genus $0$ Heegaard diagram $\mathbb{G}^{\textrm{sp}}$ with $|\bm{\alpha}|=|\bm{\beta}|=m-1$ called a \emph{spherical grid diagram} as follows:

\begin{enumerate}
    \item Place $(m-1)$ $\alpha$ and $\beta$ curves that intersect in two $(m-2)\times (m-2)$ grids on the top and bottom of the sphere.
    \item  The bottom grid contains no markings.
    \item The top grid is marked with basepoints according to $\mathbb{G}$. More precisely, the markings in rows and columns $1$ and $m$ of $\mathbb{G}$ are placed to the left, right, top and bottom of $\beta_1$, $\beta_{m-1}$, $\alpha_1$ and $\alpha_{m-1}$ respectively as indicated by the pink shading in Figure \ref{sphericalgriddiagram}, and the center $(m-2)\times(m-2)$ subgrid of $\mathbb{G}$ gets placed in the top $(m-2)\times (m-2)$ grid of $\mathbb{G}^{\textrm{sp}}$.
\end{enumerate} 

Apply this construction to $\mathbb{G}_S$, $\mathbb{G}_{S_-}$ and $\mathbb{G}_{S_0}$ to obtain spherical grid diagrams $\mathbb{G}_S^\textrm{sp}$, $\mathbb{G}_{S_-}^\textrm{sp}$ and $\mathbb{G}_{S_0}^\textrm{sp}$. Let $c$ denote the center point of the top grid of $\mathbb{G}_{S_-}^{\textrm{sp}}$ and $\mathbb{G}_{S_0}^{\textrm{sp}}$. Define $I_-$ and $I_0$ to be those states containing $c$ in $\widetilde{CFL}(\mathbb{G}_{S_-}^{\textrm{sp}}$ and $\widetilde{CFL}(\mathbb{G}_{S_0}^{\textrm{sp}})$ respectively. Similarly define $N_-$ and $N_0$ to be those states not containing $c$ in $\widetilde{CFL}(\mathbb{G}_{S_-}^{\textrm{sp}}$ and $\widetilde{CFL}(\mathbb{G}_{S_0}^{\textrm{sp}})$ respectively. We reuse the notation of the proof of Proposition \ref{skeintriangles} to facilitate a direct comparison between that proof and the proof below of Proposition \ref{alexandergradingshift1}. 
\begin{proposition}\label{alexandergradingshift1}
    The map $F(\mathbf{x})=\mathbf{x}\cup c$ is a chain isomorphism 
    \begin{equation}
        F\colon \widetilde{CFL}(\mathbb{G}_S^{\textrm{sp}})[[2j,\frac{1}{2}]]\cong I_0
    \end{equation}
    for some integer $j.$
\end{proposition}

\begin{proof}
There is an isomorphism $\psi\colon N_0\to N_-$ defined as the identity map on tuples of grid intersection points, and an isomorphism $\widetilde{CFL}(\mathbb{G}_S^{\textrm{sp}})\cong I_0$ defined by taking a generator in $\widetilde{CFK}(\mathbb{G}_S^{\textrm{sp}})$ and mapping it to the corresponding generator in $I_0$ containing $c.$ Consider the composition of chain maps 
\begin{equation}
    (I_0, \partial_{I_0}^{I_0})\xrightarrow[(-1,0)]{\partial_{I_0}^{N_0}}(N_0, \partial_{N_0}^{N_0})\xrightarrow[]{\psi}(N_-, \partial_{N_-}^{N_-})\xrightarrow[(-1,0)]{\partial_{N_-}^{I_-}}(I_-, \partial_{I_-}^{I_-})
\end{equation}
    If we can show that the bidegree of $\psi$ is $(0,0)$ and that the map $\xi\colon I_0\to I_-$ defined as the identity on tuples of grid states containing $c$ has bidegree $(0,1)$, then the argument given in the proof of Proposition \ref{skeintriangles} applies verbatim to prove this proposition as well. In what follows we will make use of the gradings $M=gr_w$ and $gr_z$ which are the direct analogs of $M_{\mathbb{O}}$ and $M_{\mathbb{X}}$ for more general Heegaard diagrams. We see that $gr_w$ is preserved by both $\psi$ and $\xi$ since both spherical diagrams have the same $w$ basepoints. We follow the grading pinning procedure demonstrated in of \cite[~Lemma 3.2, Figure 8]{MOSgrids} to determine the $gr_z$ shift of $\psi$. Consider a state $\mathbf{y}\in\widetilde{CFL}(\mathbb{G}_{S_-}^{\textrm{sp}})$ and a state $\mathbf{x}\in \widetilde{CFK}(\mathbb{G}_{S_0}^{\textrm{sp}})$ that differ only in that they contain the intersection points southwest of the two central $z$'s on their respective diagrams. Then $\mathbf{x}$ and $\mathbf{y}$ must in fact have the same $gr_z$ grading since upon forgetting about $w's$ and handle-sliding all but one of the $\beta$ curves to become small loops about $z$'s, these two states correspond through empty Maslov index 1 triangles to the \emph{same} bottom-most $gr_z$ state in the complex used to pin down the $gr_z$ grading. This is illustrated in Figure \ref{grV}. So we see that \[gr_z(\mathbf{x})-gr_z(\psi(\mathbf{x}))=gr_z(\mathbf{y})-gr_z(\psi(\mathbf{x}))=-1.\]The last equality follows from applying the relative $gr_z$ grading formula, which states that $gr_z(\mathbf{y})-gr_z(\psi(\mathbf{x}))=\mu(\phi)-2n_z(\phi)$ with $\phi$ equal to the one by one square containing a $z$ and connecting $\mathbf{x}$ and $\mathbf{y}$ in $\mathbb{G}_{S_-}^{\textrm{sp}}$. The Alexander grading is $A=\frac{gr_w-gr_z}{2}-\frac{m-\ell}{2}$ where $m$ is the number of basepoints on the Heegaard diagram and $\ell$ is the number of link components. Therefore we see that

\begin{align*}A(\mathbf{x})=&\frac{gr_w(\mathbf{x})-gr_z(\mathbf{x})}{2}-\frac{(k+1)-1}{2}=\frac{gr_w(\psi(\mathbf{x}))-(gr_z(\psi(\mathbf{x}))-1))}{2}-\frac{(k+1)-1}{2}=\\&\frac{gr_w(\psi(\mathbf{x}))-gr_z(\psi(\mathbf{x})))}{2}-\frac{(k+1)-2}{2}=A(\psi(\mathbf{x}))
\end{align*}
and hence $deg(\psi)=(0,0)$. A very similar computation shows that $deg(\xi)=(0,1).$
\begin{figure}[htb!]
 \center
 \includegraphics[scale=1.5]{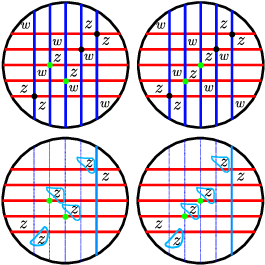}
 \caption{Top: The top halves of the diagrams $\mathbb{G}_{S_0}^{\textrm{sp}}$ and $\mathbb{G}_{S_-}^{\textrm{sp}}$. Bottom: Both diagrams are the same up to isotopy after deleting $w$ basepoints and performing the displayed handleslides amongst the $\beta$ curves.} \label{grV}
 \end{figure}
\end{proof}
Now let $*$ stand for either $0$ or $-$. Next we modify the construction of $\mathbb{G}_{S_*}^\textrm{sp}$ to accomodate the symmetry axis $A$. This construction produces more generally a spherical grid diagram compatible with an annular link. Noticing that the grid number $2k+2$ of $\mathbb{G}_{S_*}$  is even, let $\mathbb{G}_{S_*}^{\textrm{I}}$, $\mathbb{G}_{S_*}^{\textrm{II}}$, $\mathbb{G}_{S_*}^{\textrm{III}}$ and $\mathbb{G}_{S_*}^{\textrm{IV}}$ denote the four quadrants of $\mathbb{G}_{S_*}$. To construct $Ax\mathbb{G}_{S_*}^{\textrm{sp}}$, place $(2k+2)$ $\alpha$ and $\beta$ curves that intersect in two $(2k+1)\times (2k+1)$ grids on the top and bottom of the sphere, and decorate these grids as follows. The bottom grid is empty besides a $w$ marking in the center. The top grid has a $z$ marking in the center, and the remaining grid formed by removing the central row and column of the top grid is filled in with $\mathbb{G}_*^{\textrm{I}}$, $\mathbb{G}_*^{\textrm{II}}$, $\mathbb{G}_*^{\textrm{III}}$ and $\mathbb{G}_*^{\textrm{IV}}$. See Figure \ref{sphericalgriddiagramwaxis}.
 \begin{figure}[htb!]
 \center
 \includegraphics[scale=.7]{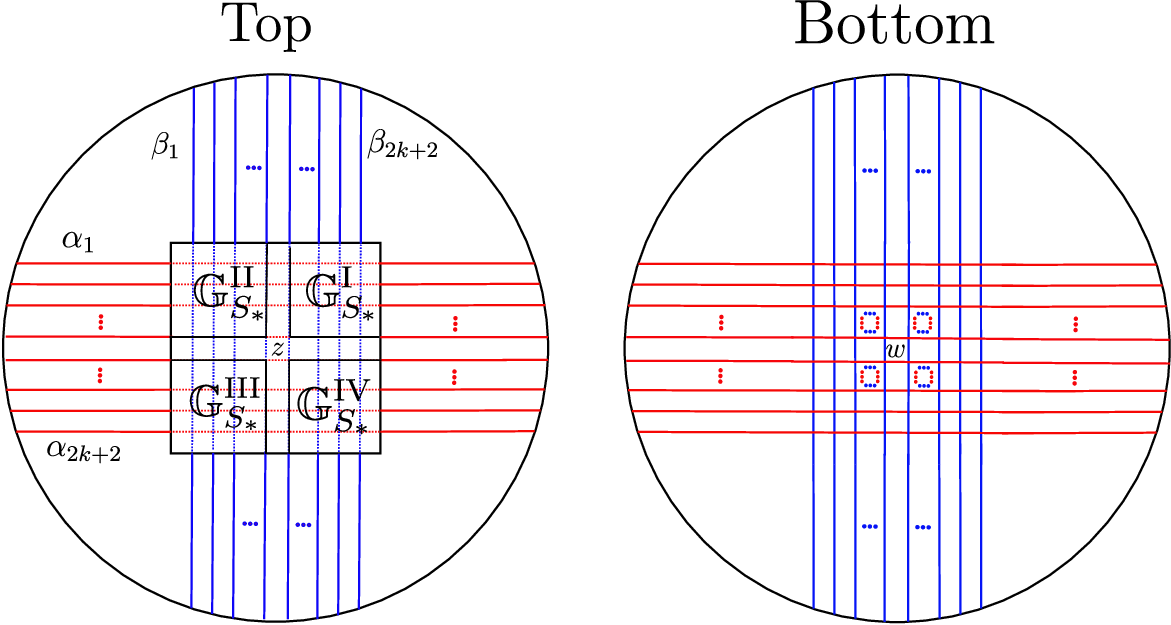}
 \caption{The Heegaard diagram $Ax\mathbb{G}_{S_*}^{\textrm{sp}}$ for the link $(S_b^n(K))_*\cup A$.}\label{sphericalgriddiagramwaxis}
 \end{figure}
 
We consider the complex $\widetilde{CFL}(Ax\mathbb{G}_{S_0}^{\textrm{sp}})$ equipped with the differential $\partial_A$ which counts disks passing over the bottom $w$ basepoint (but not the top central $z$ basepoint). See \cite[~Propostion 7.1]{OSlink} for more details; also see Proposition 2.6 and the discussion preceding it in \cite{HendricksDP}. This complex has a single Alexander grading corresponding to $L^{n+\frac{1}{2}}_b(K)$ and its homology is $\HFKhat(L_b^{n+\frac{1}{2}}(K))\otimes V^{\otimes (2k+1)}\otimes W$ where $W=\F_{(0,0)}\oplus \F_{(-1,0)}$. From now on when we write $\widetilde{CFL}(Ax\mathbb{G}_{S_0}^{\textrm{sp}})$ we mean the singly Alexander graded complex with differential $\partial_A$. Let $a$ and $b$ denote the northeast and southwest intersection points of the square containing the central $z$, and let $AxI_0$ be the quotient complex of $\widetilde{CFL}(Ax\mathbb{G}_{S_0}^{\textrm{sp}})$ generated by states containing both $a$ and $b$. The following proposition follows from a similar argument to the one given in the proof of Proposition \ref{alexandergradingshift1}. 
\begin{proposition}\label{alexandergradingstatement2}

The map $\Phi\colon I_0\to AxI_0$ sending a grid state $\bf{x}$ containing $c$ to the corresponding grid state $\Phi(\bf{x})$ containing $a$ and $b$ (i.e the rest of the constituent grid intersection points in $\bf{x}$ and $\Phi(\bf{x})$ agree in the quadrants of $\mathbb{G}_{S_*}$) is a chain isomorphism of bi-degree $(0,0)$. 
     
\end{proposition}

\begin{figure}
\centering
\begin{minipage}{.5\textwidth}
  \centering
  \includegraphics[width=.7\linewidth]{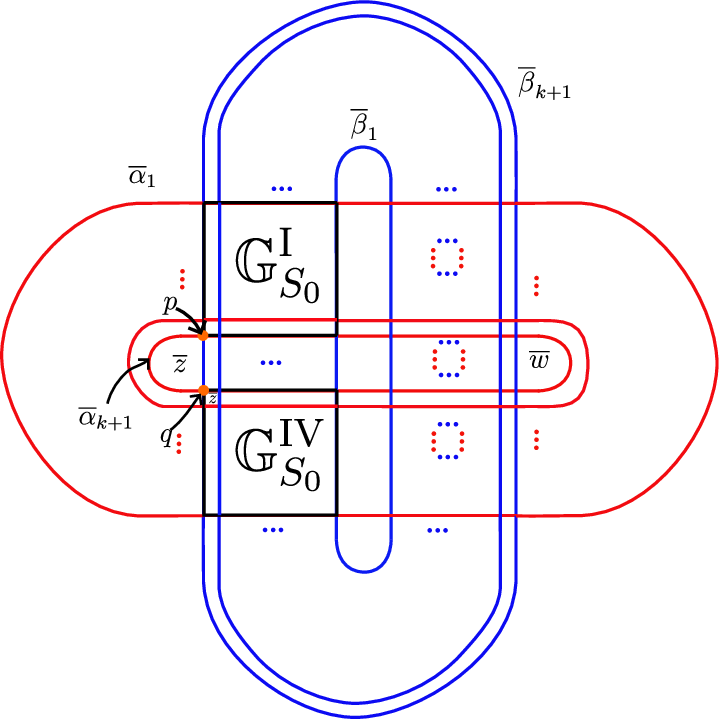}
  \captionof{figure}{The quotient diagram $Ax\mathbb{G}_{S_0}^{\textrm{sp}}/\tau$.}
  \label{sphericalquotient}
\end{minipage}%
\begin{minipage}{.5\textwidth}
  \centering
  \includegraphics[width=.7\linewidth]{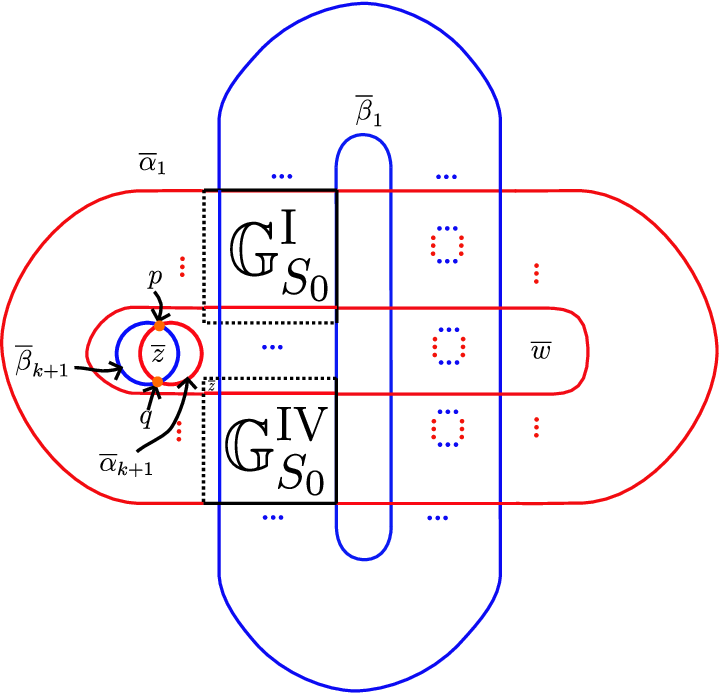}
  \captionof{figure}{The diagram $B$ obtained by isotoping $\overline{\alpha}_{k+1}$ over the rightmost $\overline{w}$ and all of the $\overline{\beta}$ curves besides $\overline{\beta}_{k+1}$ and then isotoping $\overline{\beta_{k+1}}$ over the rightmost $\overline{w}$ and over all the $\overline{\alpha}$ curves besides $\overline{\alpha}_{k+1}$.}
  \label{isotopedsphericalquotient}
\end{minipage}
\end{figure}

 In Figure \ref{sphericalquotient}, we show the quotient diagram $Ax\mathbb{G}_{S_0}^{\textrm{sp}}/\tau$, which is a Heegaard diagram for $\overline{L_b^{n+\frac{1}{2}}(K)}\cup \overline{A}$. We use the notation $\overline{L_b^{n+\frac{1}{2}}(K)}$ instead of $\overline{K}$, even though these are the same quotient knots, to emphasize that $\overline{A}$ does not intersect the quotient. Let $\partial_{\overline{A}}$ be the differential on the complex $\widetilde{CFL}(Ax\mathbb{G}_{S_*}^{\textrm{sp}}/\tau)$ which ignores the $w$ marking on the bottom grid, which is to say $\partial_{\overline{A}}$ includes disks passing over the bottom $\overline{w}$ basepoint. As in the previous paragraph, this singly Alexander graded complex has homology $H_*(\widetilde{CFL}(Ax\mathbb{G}_{S_*}^{\textrm{sp}}/\tau),\partial_A)\cong \HFKhat(\overline{K})\otimes V^{\otimes k}\otimes W$. From now on we assume that $\widetilde{CFL}(Ax\mathbb{G}_{S_*}^{\textrm{sp}}/\tau)$ is equipped with the differential $\partial_{\overline{A}}$ and correspondingly has a single Alexander grading.

Figure \ref{isotopedsphericalquotient} shows $Ax\mathbb{G}_{S_*}^{\textrm{sp}}/\tau$ after isotopies of $\alpha_{k+1}$ and $\beta_{k+1}$ that go over the rightmost $w$.  Such isotopies do not affect the homology of the complex $(\widetilde{CFL}(Ax\mathbb{G}_{S_*}^{\textrm{sp}}/\tau), \partial_{\overline{A}})$. This isotoped diagram $B$ \emph{is} the quotient diagram $\mathbb{G}_S^{\textrm{sp}}/\tau$ except for the addition of
 \begin{itemize}
     \item the isotoped $\overline{\alpha}_{k+1}$ and $\overline{\beta}_{k+1}$ curves that now meet only each other in the two points $p$ and $q$, 
     \item the $z$ basepoint contained in the central bigon formed by $\overline{\alpha}_{k+1}$ and $\overline{\beta}_{k+1}$, and
     \item the rightmost $\overline{w}$ that the differential $\partial_{\overline{A}}$ does not see.
 \end{itemize}
 With these observations and notation in place, we can state and prove the following proposition
 \begin{proposition}\label{quotientgradingshift}
    The map \[\mathcal{P}\colon  \widetilde{CFL}(\mathbb{G}_S^{\textrm{sp}}/\tau)\to  \widetilde{CFL}(Ax\mathbb{G}_{S_0}^{\textrm{sp}}/\tau)\] defined by $\mathcal{P}(\mathbf{x})=\mathbf{x}\cup p$ preserves Alexander gradings.
\end{proposition}
\begin{proof}
 Since the region of $B$ that $\alpha_{k+1}\cup \beta_{k+1}$ is contained in also contains a $z$ marking, namely the top left corner of $\mathbb{G}_{S_0}^{\textrm{IV}}$, we obtain a grading preserving isomorphism of chain complexes \[(\widetilde{CFL}(B),\partial_{A})\cong \widetilde{CFL}(\mathbb{G}_S^{\textrm{sp}}/\tau)\otimes W.\] The gradings on the tensor factor of $W=\F_{(0,0)}\oplus \F_{(1,0)}$ (the two summands correspond to those generators in $\widetilde{CFL}(B)$ containing either $p$ or $q$) are fixed by the homology of the complex $(\widetilde{CFL}(B),\partial_{A})\cong (\widetilde{CFL}(\mathbb{G}_{S_0}^{\textrm{sp}}/\tau), \partial_{A})$. In particular, 
\[
P\colon \widetilde{CFL}(\mathbb{G}_S^{\textrm{sp}}/\tau)\rightarrow \widetilde{CFL}(B, \partial_{A})\] defined on a generator $\mathbf{x}\in \widetilde{CFL}(\mathbb{G}_S^{\textrm{sp}}/\tau)$ by $P(\mathbf{x})=\mathbf{x}\cup p$ is an Alexander grading preserving map. The map $P'$ sending a generator containing $p$ in $\widetilde{CFL}(B)$ to the corresponding generator containing $p$ in $\widetilde{CFL}(\mathbb{G}_{S_0}^{\textrm{sp}}/\tau), \partial_{A})$ is also Alexander grading preserving. Finally we have that $\mathcal{P}=P'\circ P.$
\end{proof}

\begin{figure}
    \centering
    \includegraphics[width=0.25\linewidth]{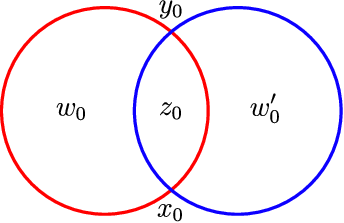}
    \caption{A type $0-1-2-3$ stabilization is the replacement of $w_0$ by the local diagram shown above.}
    \label{0123stab}
\end{figure}

We are now able to put together the results from this appendix to prove that the grading shift of Theorem \ref{Mainthm2} is as stated.
\begin{proof}
Consider a $\tau$-equivariant generator $x\in \widetilde{CFL}(\mathbb{G}_S^{\textrm{sp}})$ and its projection $\overline{x}\in \widetilde{CFL}(\mathbb{G}_S^{\textrm{sp}}/\tau)$.
Recall the chain maps $F:\widetilde{CFL}(\mathbb{G}^{\textrm{sp}})\to I_0$  and $\Phi:I_0\to AxI_0$. From Propositions \ref{alexandergradingshift1} and \ref{alexandergradingstatement2}, we see that $A(\Phi(F(x)))=A(x)-\frac{1}{2}$. Let $\overline{\Phi(F(x))}\in \widetilde{CFL}(Ax\mathbb{G}_{S_0}^{\textrm{sp}}/\tau)$ be the projection of the $\tau$-equivariant generator $\Phi(F(x))\in \widetilde{CFL}(Ax\mathbb{G}_{S_0}^{\textrm{sp}})$, and let $\lambda=lk(L^{n+\frac{1}{2}}_b(K),A)$. From \cite[~Theorem 1.3]{HendricksDP}, if $A(\Phi(F(x)))=2a+\frac{1-\lambda}{2}$, then $A(\overline{\Phi(F(x))})=a+\frac{1-
\lambda}{2}$. From the constructions of the maps $F,\Phi$ and $\mathcal{P}$ (cf. Proposition \ref{quotientgradingshift}), we see that $\overline{\Phi(F(x))}=\mathcal{P}(\overline{x})$. The map $\mathcal{P}$ has $0$ Alexander grading shift, so putting it all together we have \[A(\overline{x})=A(\overline{\Phi(F(x))})=\frac{1}{2}A(\Phi(F(x)))+\frac{1-\lambda}{4}=\frac{1}{2}A(x)-\frac{\lambda}{4}.\]Specializing to $n=0$ and applying Equation \ref{linkingnumberequation}, we see that $\lambda=2\widetilde{lk}(K)+1$. Therefore $A(\overline{x})=\frac{1}{2}A(x)-\frac{\widetilde{lk}(K)}{2}-\frac{1}{4}$. In other words, if $A(x)=2a+\widetilde{lk}(k)-\frac{1}{2}$, then $A(\overline{x})=a$.

From here on, diagram means $\tau$-equivariant Heegaard diagram. All that is left to show is that the Alexander grading shift is independent of the choice of diagram $\mathcal{H}_S$, so that the Alexander grading shift associated to $\widetilde{CFK}(\mathbb{G}_S^{\textrm{sp}})$ and $\widetilde{CFK}(\mathbb{G}_S^{\textrm{sp}}/\tau)$ is the same as the Alexander grading shift of Theorem \ref{Mainthm2}. In the proof of Theorem \ref{Mainthm2}, we showed that the spectral sequence is independent of triply pointed diagram used to compute it. In particular the Alexander grading shift is independent of triply pointed Heegaard diagram. That proof of invariance generalizes straightforwardly to show that the spectral sequence analagous to Theorem \ref{Mainthm2} computed from a $4k+3$ pointed diagram for $S^n_b(K)$ is independent of choice of $4k+3$ pointed diagram. In particular, for fixed $k$, the Alexander grading shift is independent of the choice of $4k+3$ pointed diagram for $S^n_b(K)$. So if we can show that taking a $4k+3$ pointed diagram $\mathcal{H}_S$, and replacing neighborhoods of basepoints $w_0$ and $\tau(w_0)$ (where $w_0\neq\tau(w_0) $) with the type $0-1-2-3$ stabilizations shown in Figure \ref{0123stab} to obtain a $4(k+1)+3$ pointed diagram $\mathcal{H}_S'$ has no effect on the Alexander grading shift, by induction we will be done. Notice that $\widetilde{CFK}(\mathcal{H}_S')=\widetilde{CFK}(\mathcal{H}_S)\otimes V\otimes V$ since any generator of $\widetilde{CFK}(\mathcal{H}_S')$ is $x\cup \{r_1,r_2\}$ for $r_1\in\{x_0,y_0\}$ and $r_2\in \{\tau(x_0),\tau(y_0)\}$. A simple computation shows that $A(x\cup \{x_0,\tau(x_0)\})=A(x)$, and so if $x$ is a $\tau$-equivariant generator of $\widetilde{CFK}(\mathcal{H}_S)$ then $x\cup\{x_0,\tau(x_0)\}$ is a $\tau$-equivariant generator of $\widetilde{CFK}(\mathcal{H}_S')$ of the same Alexander grading. But looking at the quotient diagrams $\mathcal{H}_S/\tau$ and $\mathcal{H}_S'/\tau$, we notice that $\mathcal{H}_S'/\tau$ is a type $0-1-2-3$ stabilization of $\mathcal{H}_S/\tau$ at $\overline{w_0}=q(w_0)=q(\tau(w_0))$, and we also see that $A(\overline{x}\cup\{\overline{x_0}\})=A(\overline{x})$.
\end{proof}

\begin{remark}
    With some more effort one can show that the spectral sequence computed from a diagram with more basepoints is the spectral sequence of Theorem \ref{Mainthm2} tensored with certain simple complexes as in Section 5.3 on equivariantly destabilizing basepoints in localization spectral sequences of \cite{Hendricks_2016}.
\end{remark}

\bibliographystyle{amsalpha}
\bibliography{bib}
\end{document}